%%%%%%%%%%%%%%%%%%%%%%%%%%%%%%%%%%%%%%%%%%%%%%%%%%%%%%%%%%%%%%%%%%%%%%%%

%%  Joseph M. Landsberg
%%  Laboratoire de Math\'ematiques,
%%  Universit\'e Paul Sabatier, UFR-MIG
%%  31062 Toulouse Cedex 4,
%%  FRANCE 

%%  jml@picard.ups-tlse.fr 
%%  05-61-55-82-84 (bureau) (from abroad: 33-5-61-55-82-84)
%%  05-61-55-82-00 (fax)    (from abroad: 33-5-61-55-82-00)
 %file snueprint.tex  
%us page lenghts
 
\input amstex

\def\BA{\Bbb A }\def\BC{\Bbb C }
\def\BO{\Bbb O }\def\BG{\Bbb G }
\def\BP{\Bbb P}

\def\a{\alpha}
\def\b{\beta}
\def\c{\gamma}
\def\fg{\frak g}
\def\xsm{X_{sm}}

\def\dtau{\delta_{\tau}}
\def\dsigma{\delta_{\sigma}}
\def\dstar{\delta_{*}}
\def\rtau{r_{\tau}}
\def\tx{\tau(X)}
\def\sx{\sigma (X)}
\def\al#1{\alpha_{#1}}
\def\aa#1{e_{#1}}
%NOTE changed def here to get std notation

%\def\al#1#2#3{a^{#1}_{{#2}{#3}}}

\def\bii#1#2{\Bbb F\Bbb F^{#1}_{#2}}
\def\bcc#1{\Bbb C^{#1}}

\def\brr#1{\Bbb R^{#1}}
\def\bpp#1{\Bbb P^{#1}}
\def\bppp{\Bbb P }\def\bccc{\Bbb C }
\def\bb#1#2#3{b^{#1}_{{#2}{#3}}}

\def\cf{\Cal F}
\def\ci{\Cal I}
\def\ee#1{e_{#1}}
\def\EE#1#2{E^{#1}_{#2}}
\def\uee#1{e^{#1}}

\def\ett#1#2{\eta^{#1}_{#2}}
\def\et#1{\eta^{#1}_0}
\def\ff#1{\Bbb F\Bbb F^{#1}}
\def\fff#1#2{\Bbb F\Bbb F^{#1}_{#2}}

\def\frp#1#2{\frac{\partial {#1}}{\partial {#2}}}
\def\gg#1#2{g^{#1}_{#2}}

\def\hd{, \hdots ,}
\def\inv{{}^{-1}}
\def\ii{| II |}

\def\La#1{\Lambda^{#1}}
\def\na{n+a}
\def\none{n+1}
\def\nl{n+l}
\def\nk{n+k}
\def\ooo#1#2{\omega^{#1}_{#2}}
\def\oo#1{\omega^{#1}_0}
\def\ot{\!\otimes\!}
\def\pii{\partial II}
\def\pp#1{\Bbb P^{#1}}
\def\ppp{\Bbb P}
\def\qq#1#2#3{q^{#1}_{{#2} {#3}}}

\def\ra{\rightarrow}

\def\rr#1#2#3#4{r^{#1}_{{#2} {#3}{#4}}}

\def\tim{\text{Image}\,}

\def\tann{\text{Ann}\,}
\def\tdim{\text{dim}\,}
\def\tcodim{\text{codim}\,}
\def\tsingloc{\text{singloc}\,}
\def\tbaseloc{\text{Base}\,}
\def\tker{\text{ker}\,}
\def\tspan{\text{span}\,}
\def\tmod{\text{ mod }}
\def\tmin{\text{ min }}
\def\tmax{\text{ max }}
\def\thom{\text{ Hom }}
\def\trank{\text{rank}\,}
\def\ttrace{\text{trace }}
\def\ud#1{\underline d^{#1}}
\def\up#1{{}^{({#1})}}
\def\upperp{{}^{\perp}}
\def\upperpq{{}^{\perp Q}}
\def\vdx{v_d(X)}

\def\vtwox{v_2(X)}

\def\ww{\wedge}
\def\xx#1{x^{#1}}
\def\xsm{X_{sm}}

\define\intprod{\mathbin{\hbox{\vrule height .5pt width 3.5pt depth 0pt %
        \vrule height 6pt width .5pt depth 0pt}}}
\documentstyle{amsppt}
\magnification = 1100
\hsize =15truecm
\hcorrection{.5truein}
\baselineskip =18truept
\vsize =21truecm
\NoBlackBoxes
%\nologo
\topmatter
\title Algebraic geometry and projective differential geometry\linebreak
  Seoul National University concentrated lecture series, 1997
\endtitle
\rightheadtext{Preliminary version}
\author
  J.M. Landsberg
\endauthor

%\abstract{ }
%\endabstract

\date {September, 1998}\enddate
\address{ Laboratoire de Math\'ematiques,
  Universit\'e Paul Sabatier, UFR-MIG
  31062 Toulouse Cedex 4,
  FRANCE}\endaddress
\email {jml\@picard.ups-tlse.fr }
\endemail
\thanks {Supported by  NSF grant DMS-9303704.}
\endthanks
\endkeywords
\subjclass{ primary 53, secondary 14}\endsubjclass
 
\endtopmatter

\document

\subheading{Contents}

0. Introduction

1. Examples of homogeneous varieties 

2. Auxilliary varieties 

3. Topology and consequences 

4. Projective differential invariants

5. Varieties with degenerate Gauss images 

6. When can a uniruled variety be smooth?

7. Dual varieties 

8. Linear systems of  bounded and constant rank 

9. Secant and tangential varieties

10. Systems of quadrics with tangential defects 

11. Recognizing uniruled varieties 

12. Recognizing intersections of quadrics 

13. Recognizing homogeneous spaces

14. Complete intersections

15. Errata and clarifications 

\subheading{ Preface}

\noindent{\bf Please note}: This is a preliminary version, the final
version
will be published in book form by Seoul National University. Please send
comments, corrections and questions to me at jml\@picard.ups-tlse.fr. I am
particularly interested in historical comments, as I have been unable to
trace
the original authorship of several results (those marked with a ??).

\smallpagebreak

This is an expanded and updated version of a lecture
series I gave at
Seoul National University  in September 1997. The series was
focus on Zak's theorem on Severi varieties and my differential-geometric
proof of the theorem at the request of my hosts. A French language version
of the introduction served as the text for my {\it habilitation
\`a diriger des recherches} (November 1997).

These notes are written with two different audiences in
mind: graduate students and  algebraic geometers
not familiar with infinitesimal methods.

For the graduate students, I have included    
exercises  and   open questions to work on. For the   algebraic
geometers not familiar with infinitesimal
methods, I have attempted to relate the
techniques used here with  standard methods in algebraic geometry.

These notes are in some
sense an  update to the   paper [GH].

   I do not discuss    
hyperdeterminants,
  Geometric Invariant
Theory, Variation of Hodge Structure,
  web geometry,  theta divisiors and other related topics,
although I would have liked to.   I strongly encourage graduate
students to study the connections with some of these other topics.

I would like to thank Jun-Muk Hwang and Seoul
National University.

\heading \S 0. Introduction\endheading

Let $X^n\subset\Bbb C\pp\na=\ppp V$ be a variety ($V=\bcc{\na+1}$). In
these
notes
I   discuss the geometry of $X$ as an algebraic variety,   the local
projective
differential geometry of the smooth points of $X$, and most importantly,
the
relations between the two.

In the spirit of F. Klein, let's consider a property of a variety {\it
geometric} if it is invariant under projective transformations
 (that is, the action of $PGL(V)$ on $\ppp V$). For example,
two geometric properties
of a variety are its {\it dimension} (the dimension of its tangent
space at a smooth point) and its {\it degree} (the number of points of
intersection with a general linear space of complementary dimension).
The first property is intrinsic, the second extrinsic.

One way to measure the pathology of $X$ is to
 construct auxilliary varieties from $X$, and to calculate the
difference between the expected and actual 
     dimensions   of these auxilliary varieties.  
 In these notes we will  study such auxilliary 
varieties using modern techniques
combined with infinitesimal methods developed by
E. Cartan and others.

 In what follows I 
give an example of an auxilliary variety,
namely the {\it secant variety} of
$X$ and discuss its study. This example is typical   and will serve as a
model
for the other cases.

\subheading{0.1 Secant varieties}

Given two points $x, y\in \ppp V$, there exists a unique
line $\BP^1_{xy}$ containing them.
Given a subvariety $X\subset\ppp V$,   define
the {\it secant variety} of $X$, 
$$
\sx:= \overline{\cup_{x,y\in X}\BP^1_{xy}},
$$
 the closure of the union of all secant lines   to
$X$. There are $2n$-dimensions   of pairs of points
on $X$ and  one parameter of points on each line, so one expects
that $\tdim\sx= 2n+1$ if $2n+1\leq n+a$, or $\sx=\pp\na$   if $2n+1>n+a$.
If not, 
we say $\sx$ is {\it degenerate} and  let $\dsigma= 2n+1-\tdim\sx$ denote
the
{\it secant defect} of $X$.

\subheading{0.1.1 Example} 
Let $V= \bcc{k+1}\ot\bcc{l+1}$ denote the space of 
$(k+1)\times (l+1)$ matrices. Let $X\subset \ppp V$ be the projectivization
of the rank one matrices.   $X\simeq \pp k\times \pp l$
because every rank one matrix is the product of a column vector with
a row vector.
$X=Seg(\pp k\times\pp l)$ is  called the {\it Segre variety}. $X$ is the
zero
locus of the   two by two minors.

The sum of two rank one matrices has rank at
most two, so   $\sigma(Seg(\pp k\times\pp
l))=\ppp (\trank\leq 2$ matrices). 

\subheading{0.1.1.1 Exercise}   $\tdim\sx = 2(k+l)-1$,
so $\sx$ is degenerate, with $\dsigma=2$.

\smallpagebreak

A general principle  is
that pathology should be rare if $X$ is smooth and $\tcodim X$ is
relatively small (see \S 3). A theorem to this effect
is the following:

\proclaim{0.1.2 Zak's theorem on linear normality [Z]} 
If $X^n\subset\pp\na$ is smooth, not contained in a hyperplane, and
$a<\frac
n2+2$, then $\sx=\pp\na$.
\endproclaim

(The name of this theorem is explained in \S 3.)

In addition,  Zak (with Lazarsfeld) classified the
varieties $X$ in the borderline case
of $a=\frac n2+2$ and $\sx$ degenerate:

\proclaim{0.1.3 Zak's theorem on Severi varieties [Z], [LV]} 
If $X^n\subset\pp\na$ is smooth, not contained in a hyperplane,   $a=\frac
n2+2$ and $\sx\neq\pp\na$,     then $X$ is one of
$$\align & i.\ \text{Veronese}\ \pp 2 \subset \pp 5\\
& ii.\ \text{Segre}\  \pp 2\times\pp 2 \subset \pp 8\\
& iii.\ \text{Pl\"ucker  embedded Grassmannian}\ G(\Bbb C^2,\Bbb
C^6)\subset
\pp {14}\\
& iv.\  E_6/P_1 \subset \pp {26}.\endalign
$$
\endproclaim
 
These four varieties are 
called the {\it Severi varieties} after F. Severi who proved the theorem in
the special case $n=2$.
They are described   in  \S 1  along with many other homogeneous
varieties. Homogeneous varieties often provide examples of extremal
pathologies so I discuss them   in \S 2.
 In [L5] I give   new proofs of Zak's theorems.  The proofs
have   five steps, which will serve as a model for
many of  the questions discussed in these notes.

\medpagebreak

\noindent {\it Step 1: Describe the pathology infinitesimally.}

The   condition  imposed on the differential invariants of $X$
  when 
$\dsigma >0$
was determined by Terracini in 1913  [T]. The condition
$\dsigma >1$   is 
essentially that the quadrics in the projective second fundamental
form satisfy a polynomial
(second fundamental forms
will be  defined shortly). See \S 9 for the precise condition. In
particular,
it is a closed condition.   For other questions, such
as  in the study of complete intersections,  this step
can be quite involved, 
see \S 12 and 0.5.

\smallpagebreak

\noindent {\it Step 2: Analyze the infinitesimal condition.}

Here one determines which
systems of quadrics satisfy the Terracini condition. This type of
question can be studied from several perspectives. To reprove
Zak's theorems,
I localize  the problem yet again and use  differential-geometric
methods, see \S 10.  In 
our study of dual varieties  
[IL], Ilic  and I  used methods from algebraic geometry
at this stage; the study of vector bundles on projective space,
see \S 7.

\smallpagebreak

\noindent {\it Step 3: Determine infinitesimal consequences of smoothness.}

In this step one studies the   additional conditions  
placed on the differential invariants
of $X$  because $X$ is   smooth (or almost smooth).
These are usually open (genericity) conditions.
Recovering global information
from  infinitesimal
invariants is   central to my research. See \S 4 and \S 7.

\smallpagebreak

\noindent {\it Step 4: Combine steps two and three.}

In the case of degenerate secant varieties,
  one  combines  the open conditions
implied by smoothness with the
restrictions on the systems of quadrics arrived at in the second step.
At  this point Zak's theorem on linear normality follows immediately.
In the Severi variety case, one is restricted to  four possible
second fundamental forms.

\smallpagebreak

\noindent {\it Step 5: Pass from infinitesimal to local (and hence global)
geometry.}

  This step of 
  passing from the infinitesimal geometry to the local geometry
(and since one is in the analytic category, the result is
a global) is made using the Cartan machinery of moving
frames and exterior differential systems.  
This step can be viewed as the
study of the  deformation theory, or rigidity, of   systems of
quadrics, cubics, etc...
Such results are discussed in \S 13.

\medpagebreak

I consider step 3 to be the most important, so let's begin there: 

\subheading{0.2 A   principle relating global smoothness to local 
projective geometry (the third step)}

In projective 
space, 
the global geometry restricts the local geometry.
One can view these restrictions as consequences of    
    the very defining property of the
projective plane:  that parallel lines meet at infinity, or, more
generally,
{\it in projective space, 
  linear spaces
(in fact arbitrary varieties) of complementary dimension  must
intersect.}

At each point $x$ of a submanifold $X\subset
 \ppp V$, 
there is a unique {\it embedded tangent space}, that is, a unique
linear space that best approximates $X$ at $x$ to first order.
 I  denote the embedded tangent space by
$\tilde T_xX$   to distinguish it from
the intrinsic holomorphic tangent space, which I denote $T_xX$.

Consider the following two surfaces
in affine space $\BA^3$
$$
\
$$

\bigpagebreak

$$ 
\text{ hyperbola }\ \ \ \ \ \ \ \ \ \ \ \   \text{ cylinder }\tag 0.2.1
$$

Both the hyperbola and the cylinder are defined by a 
quadratic equation, and both are ruled by lines. Both can be
completed to projective varieties
by considering $\Bbb A^3\subset\pp 3$. When one completes the hyperbola,
one obtains a smooth quadric surface in $\pp 3$. In contrast,
  completing the cylinder, one
obtains a singular cone. 
The cylinder
obtains a singularity because as one travels
along one of its  rulings,  the embedded tangent space $\tilde T_xX$ is
constant. This   forces the rulings to crash into each other
at infinity.
In contrast,  the embedded tangent space
of the hyperbola  rotates as one travels along
a ruling and a singularity   at infinity is thus avoided.

The contrast between these cases leads to the following principle:

\proclaim{0.2.2 Smoothness Principle [L6]}In order for $X$ to be
smooth, its embedded tangent space must \lq\lq move enough\rq\rq .
\endproclaim

How much the tangent space needs to move will depend on 
$\tdim X$ and $\tcodim X$.

 In order to make the smoothness principle precise,
we need a way to measure \lq\lq how much\rq\rq\ 
$\tilde T_xX$ is moving. We will make such measurements using
the   {\it projective second fundamental form}. 
 
\subheading{0.3 The   projective differential geometry of $X$}

Recall that in Euclidean geometry, 
the basic measure
of how a submanifold 
of Euclidean space is bending (that is, moving away from
its embedded tangent space to first order) is the Euclidean
second fundamental form. In projective
geometry, there is a  {\it projective second fundamental form} that 
can be defined the same way as its Euclidean analogue:

\subheading{0.3.1 The Gauss map $\gamma$ and a definition of $II$ via
$\gamma$}

A natural way to keep track of the motion
of $\tilde T_xX$  is the {\it Gauss map}
$$
\align
\gamma : X&\ra \Bbb G(n, \ppp V)\tag 0.3.1.1 \\
x &\mapsto \tilde T_xX.\endalign
$$
where $\BG (n, \ppp V)$ is the {\it Grassmanian} of $\pp n$'s in
$\ppp V$. A $\pp n$ in $\ppp V$ is equivalent to an $n+1$ plane passing
through the origin in $V$. I use the notation $G(n+1,V)=
\BG (n, \ppp V)$ when I want to emphasize this 
second description.

To measure how $\tilde T_xX$ moves to first order, one
calculates the derivative of $\gamma$:
$$
\gamma_{*x} : T_{x}X\ra T_{\hat T_{x}}G(n+1,V)
= \hat T^*\otimes (V/\hat T) =\hat T^*\otimes N(-1). \tag 0.3.1.2
$$

Here I use the notation $N_xX= V/\hat T\ot \hat x^*$ and
$N(-1)=V/\hat T_xX$. If you are not familiar with such notations,
there is little harm in thinking of both $N$ and $N(-1)$ as
the normal space to $X$ at $x$. (If you are
concerned, see (1.2).)   Here and in what follows, I often
omit reference to $X$ and the base point $x$.

$\gamma_{*x}$ is such that  the kernel of the endomorphism
$\gamma_*(v):\hat T\ra N(-1)$ contains $\hat x\subset\hat T$
 for all $v\in T$.  Thus  $\gamma_*$ factors to a map
$$
\gamma_*' : T\ra (\hat T/\hat x )^*\otimes N(-1)= T^*\otimes N. \tag
0.3.1.3
$$
Furthermore, $\gamma_*'$ is symmetric, 
essentially because the Gauss map is
already the derivative of a map and mixed partials commute,
(see (4.4) for a proof).  It is called the
{\it projective second fundamental form}, and denoted
$$
II:= \gamma_*'\in S^2T^*\otimes N. \tag 0.3.1.4
$$

Roughly, $II$ measures how $X$ moves away from its tangent space at each
point
to first order.  In other words, if one considers $X$ as being mapped into
$\ppp V$, $II$ is  the geometrically relevant second derivatives of the
mapping.     

\smallpagebreak

 Note  that $II$ is an {\it algebraic} object, and thus
can be used in finite characteristic etc... as long as the usual
precautions are taken.

\smallpagebreak

\subheading{0.3.2 A coordinate definition of II}

Let $x\in X$ be a smooth point. Choose local coordinates 
$(\xx 1\hd\xx\na )$ around $x$ such that
$x=(0\hd 0)$ and
$T_xX= \{ \frp {}{\xx \alpha}\}$, $1\leq\alpha,\beta\leq n$
and $N^*_xX = \{ d\xx\mu\}$, $n+1\leq\mu,\nu\leq\na$.
Write $X$ locally as a graph $\xx\mu = f^{\mu}(\xx\alpha )$. Then,
in these coordinates, 
the
projective second fundamental form of $X$ at $x$ is:
$$
II_{X,x}=  \frac{\partial^2 f^{\mu}}{\partial \xx\alpha
\partial \xx\beta}|_x
d{\xx\alpha}\circ d{\xx
\beta}\ot 
\frp {}{\xx\mu} .
\in S^2T^*_xX\ot N_xX.\tag 0.3.2.1
$$

\subheading{0.3.3 Interpretations and measurements}

It is convenient to consider $II$ as a map $II :N^*\ra S^2T^*$
(dual to the standard Euclidean perspective) and to set $\ii = \ppp
II(N^*)$.  One can think of
$\ii$ as a linear family of quadric hypersurfaces in $\ppp T$.

\smallpagebreak
 
$\ppp N^*_xX$ has the geometric interpretation as the
space of hyperplanes tangent to $X$ at $x$, i.e., 
the  hyperplanes $H$
such that $X\cap H$ is singular at $x$.  
$II$ is essentially the map that (up to scale)  sends a hyperplane
to the quadratic part of the singularity
of $X\cap H$  at $x$. 

\smallpagebreak

Step one for the secant variety problem leads to the following condition:
if  $a\leq n$, then $\dsigma >1$ implies that for all $v\in T$,
there exists $q\in \ii$ such that $[v]\in q_{sing}$. See \S 8 for
the degeneracy condition in the general case. The study of
systems of quadrics satisfying this condition 
(the second step) is complicated, so
I will wait until \S 10 to discuss it.

\medpagebreak

Returning to  the general study of $\ii$, a natural question
is: 

\smallpagebreak

{\it How much of the geometry of $X$ is determined by $\ii$ at a
general point?}

\smallpagebreak

To investigate this question, we need to
extract geometric information from $\ii$.

\smallpagebreak

In projective geometry, unlike Euclidean geometry, 
one cannot measure   how fast a submanifold is
bending, but only whether or not it is bending. 

For example, in projective
geometry,  one can measure if a line in the embedded tangent space $\tilde
T$ 
  osculates to order two, i.e.,  if $X$ appears to contain the line to
second order. Let 
$\tbaseloc\ii\subset\ppp T$ denote the variety of
  directions tangent to
the  lines that osculate
to order two at $x$. That is,
$$
\align 
\tbaseloc\ii &= \ppp \{ v\in T \mid II(v,v)=0\}\\
&= \{ [v]\in \ppp T \mid [v]\in
q\, \forall q\in\ii \}.\tag 0.3.3.1 \endalign
$$
Given $H\in \ppp N^*_xX$, the quadric hypersurface
$q_H \subset\ppp T_xX$ is
the set of tangent directions not moving away from $H$ to second order.
If $c(t)$ is a curve lying in $X$ with $c(0)=x$, then $[c'(0)]\in 
\tbaseloc\ii$ if and only if $[c''(0)]\in\tilde T_xX$ (assuming
$c''(0)\neq 0$).

A stronger condition on a tangent vector
$v$  is   that the embedded
tangent space does not move to first order in the direction of $v$. 
Let $\tsingloc\ii$ denote the set of such directions. That is,
$$
\align 
\tsingloc\ii :&=
\ppp\{  v \in   T \mid
II(v,w)=0\forall w\in T\}\\
&= \{ [v]\in \ppp T \mid [v]\in q_{sing}\forall
q\in\ii
\}.\tag 0.3.3.2 \endalign
$$
Note that   the first equality
implies that $\tsingloc\ii$ is a linear subspace
of $\ppp T$.  If $c(t)$ is a curve lying in $X$ with $c(0)=x$, then
one can interpret the condition 
$c'(0)\in \tsingloc\ii$ as follows:  let  $E_t:=
\tilde T_{c(t)}X$ denote  
the projective bundle defined over the curve $c(t)$.
 $[c'(0)]\in 
\tsingloc\ii$ if and only if $E_0'\subseteq\hat T_{c(0)}X$.

 \smallpagebreak

  In terms of our invariants of the second fundamental
form, directions tangent to the rulings of
the cylinder are in $\tsingloc\ii$.
Directions tangent to the rulings of
the hyperbola  are in $\tbaseloc\ii$,  but not $\tsingloc\ii$.

The study  of   projective second fundamental forms leads one
back to   algebraic geometry: the geometry
of systems of quadrics.  Sections 
\S 8 and \S 10 are dedicated to studying
properties of systems of quadrics.

The 
smoothness principle is 
illustrated by the following theorem of Griffiths
and Harris (??)
that generalizes the example of the cylinder above:

\proclaim{0.3.3.3 Theorem, [GH]}Let  $X^n\subset\pp{n+a}$ be a variety.   
Let
$x\in X$ be a general point. If 
 $ \tsingloc\ii_x\neq\emptyset$ then $X$ is singular.\endproclaim

A more precise version of this theorem is given in \S 5.

The following theorem is another illustration of the smoothness
principle. It states that the embedded tangent space must move 
  away from 
{\it each} tangent hyperplane in a minimum number of directions:
 
\proclaim{0.3.3.4 Theorem, (special case of) rank restrictions, [L3]  }
Let  $X^n\subset\pp{n+a}$ be a variety.   Let $b=$dim($X_{sing})$. (Set
$b=-1$
if $X$ is smooth.)  Let $x\in X$ be a general point.

1. For any quadric $q\in\ii_x$,
 $$
\text{dim(Singloc }q)\leq 2(a-1)+(b+1).
$$

2. For generic quadrics $q\in\ii_x$,
$$
\text{dim(Singloc }q)\leq a-1+(b+1).
$$
\endproclaim

The rank restriction theorem is analogous to
Bochner type formulae in differential geometry  
in which global considerations
impose pointwise  conditions on the curvature.

Varieties with the degeneracies that are discussed in these notes 
can be viewed  as solutions to systems of partial differential equations.
   From this perspective,
 the rank restriction
theorems have the effect of ruling out characteristic or
degenerate initial 
data to initial value problems.

A long term goal is to prove stronger and higher order
versions of the rank restriction theorem.  

\smallpagebreak

Before describing other problems and results, I will briefly 
motivate the method of calculation, the {\it moving frame}.

\smallpagebreak

\subheading{0.4 The moving frame}

Let $X\subset\ppp V$ be a variety and $x\in X$ a smooth point.
Then $x\in X$ determines a flag
$$
\hat x\subset\hat T_xX\subset V.
$$
One way to study the geometry of $X$ is to examine how this
flag varies as we vary $x$. For example, the second fundamental form
contains the information of how $\hat T$ moves to first order in
relation to the motion of $\hat x$.
When studying particular geometric properties, one often
refines this flag. For example, in the study of varieties
with degenerate secant varieties, it is convienent to
also fix a vector $v\in \hat T/\hat x$, to obtain a flag
$$
\hat x\subset\{\hat x,v\}\subset\hat T\subset V.
$$
This flag admits further refinements, e.g., there is
a subspace $\tann (v)\subset N_xX$ corresponding to the quadrics
  $q_H\in \ii$ annhilating $v$, and
also  a subspace $\tsingloc\tann (v) \subset
\hat T/\hat x$  consisting of all the vectors
$v\in T$ annhilated by $\tann (v)$. Thus we obtain a refined flag:
 $$
\hat x\subset\{\hat x,v\}\subset \{\hat x,\tsingloc (\tann (v) ) \}
\subset\hat T\subset \{\hat T , \tann (v)\} \subset V .
$$
One would like to keep track of the infintesimal motions
of all these spaces and their relations. Fortunately there
is a method developed by  G. Darboux, E. Cartan and others, the {\it moving
frame}, exactly designed for this purpose. It keeps track of
the relations between all the infinitesimal motions, and
is designed in such a way that as further refinements of the
flag are made, a minimal amount of additional work is necessary.
 From a geometrical perspective, it would be
natural to begin  on the manifold of complete flags of
$V$. However a computational advantage is gained if one works on
a slightly larger space, the space of bases, or {\it framings}
of $V$, as this is a Lie group ($GL(V)$ in fact) and on a
Lie group, derivatives
can be calculated algebraically.   Moving frames for subvarieties
of projective space are discussed in \S 4.

\smallpagebreak

In the remainder of this introduction, I describe additional problems
in projective geometry and results I have obtained using the methods
discussed above in their study. In the chapters I have attempted
to give an overview of what is known in general regarding these
questions.

\subheading{0.5 Complete intersections}

The least pathological algebraic varieties are the smooth
hypersurfaces. For example, the dimension is
obvious   and the degree is simply
the degree of the single polynomial
defining $X$. (Varieties are reduced and irreducible.)

A class of varieties that share many of the simple properties
of hypersurfaces is the class of {\it complete intersections}.

\subheading{0.5.1 Definition}
A variety $X^n\subset\pp\na$ is a {\it complete
intersection} if  the ideal of $X$, $I_X$,
can be generated by $a$ elements.

\smallpagebreak

 The following example shows how   varieties that are not complete
intersections can arise:

\subheading{0.5.2 Example} In $\bpp 3$, let
$Q_1,Q_2$ be quadric hypersurfaces. Consider
$X=Q_1\cap Q_2$ (the common zero locus of
two degree $2$ homogeneous polynomials).
$X$ is an algebraic set of dimension one and
degree four.
If $Q_1,Q_2$ are reasonably general, then $X$ is a curve of degree four:
$$
\ 
$$
\bigpagebreak

$$
\ \tag 0.5.3
$$

\noindent But consider the following example: let $(\xx 1\hd\xx 4 )$
be linear coordinates on $\bcc 4$, and let
$$
\align &Q_1 = \xx 1\xx 4 -\xx 2\xx 3\\
&Q_2 = (\xx 2)^2-\xx 1\xx 3\endalign
$$
Then $X$ is  a curve of degree three plus a line:
$$
\
$$
\bigpagebreak

$$
\ \tag 0.5.4
$$

$C= [s^3,s^2t,st^2,t^3]$, $[s,t]\in \bpp 1$,
$l= [0,0,u,v]$, $[u,v]\in\bpp 1$.

    Varieties are irreducible, thus we need to get
rid of one of these 
components. We pretend we understand lines, so we elimiate $l$ 
by intersecting $X$ with $Q_3=(\xx 3)^2-\xx 2\xx 4$ to be
left with the cubic curve.
Degree $(C)=3$, so  $C$ cannot be  the intersection of two hypersurfaces.
(If it were, it would have to be in a cubic hypersurface
and a hyperplane, but $C$ is not contained in any linear subspace.)  

\smallpagebreak

 The difference between the topology
  of complete intersections and non-complete intersections  
has been
studied extensively. I have attempted to understand how the projective
differential geometries of the complete intersections
and non-complete intersections differ.

If one compares the two pictures
$$
\ 
$$

\bigpagebreak

$$
\ \tag 0.5.5
$$

\noindent and tries to understand   the difference between
them, one might say that
the non-complete intersection \lq\lq bends less\rq\rq ,
or that its tangent space \lq\lq moves less\rq\rq\ 
than the complete intersection. (This idea can be made precise
 if one adds a K\"ahler metric and is willing to integrate.
Since we will work locally, this is not   what we will do.)
This idea  is central to what follows, so
I record it informally:

\proclaim{0.5.6 Complete intersection principle} If  
    $X$ \lq\lq bends
enough\rq\rq , then $X$ will be a complete intersection.
\endproclaim

The determination of how much is \lq\lq enough\rq\rq\
 will be based on
information about the degrees of hypersurfaces containing $X$. 
In our example, the cubic curve 
\lq\lq bends less\rq\rq\ than a complete intersection of quadrics
would.

A precise explanation of the phrase \lq\lq bends enough\rq\rq\
is given below.
For now, consider the curve in affine space $y= x^3$.
$$
\ 
$$

$$
\ \tag 0.5.7
$$
At the origin, the embedded tangent space \lq\lq moves less\rq\rq\
 than at
other points on the curve in the sense that the curve {\it osculates}
 to its tangent line
to order two at the origin (versus order one for all other points). 
  (0.5.7) motivates one to study bending.
However,  the actual
 type of bending we will  study  will 
concern whether or not  there is  osculation to orders 
higher than expected (in dimensions greater than one).

\smallpagebreak

Let's now discuss
what one means by orders of osculation \lq\lq higher than
expected\rq\rq .
To decide what order of osculation is expected, 
 I determine 
  {\it a priori} information about osculating hypersurfaces.
The following results are of interest, independently of the study
of complete intersections.

\smallpagebreak

If $x\in X$ is a smooth point, then there is always an
$(a-1)$-dimensional space of hyperplanes (degree one hypersurfaces)
tangent
(osculating to order one) at $x$.
The following proposition generalizes the case of hyperplanes:

 \proclaim{0.5.8 Proposition   [L4,3.16]} Let $X^n\subseteq 
 \bpp{n+a}$ be a  variety and let  $x\in \xsm$.
 For all $  p\leq d$,  
$$
\align &
\tdim
\left\{
\matrix  \text{(not necessarily irreducible)
hypersurfaces} \\\text{of degree }d\text{ osculating   to order }  p 
\text{
at }x\endmatrix \right\} \\
&
=    \binom{n+a+d}  d - \binom {n+p}p.\endalign
$$
\endproclaim

 For $k>d$, the dimensions
of the spaces of osculating hypersurfaces depend
on the geometry of $X$.
One might think that the pattern would
continue,
so that if
$a$ is relatively small one would  expect that generically   no
hypersurfaces would  osculate to order $d+1$.
However, at higher orders, a new phenomena occurs
because of hypersurfaces that are
singular at $x$. 
Independent of $X$,  
for $d+1\leq k\leq 2d-1$,  there
are {\it lower} bounds on the dimensions
of the space of hypersurfaces of
degree $d$ osculating to order $k$ at $x$. For example:

 \proclaim{0.5.9 Proposition  [L4,3.17]} Let $X^n\subseteq 
 \bpp{n+a}$ be a  variety,
 and let $x\in \xsm$.
$$
\tdim
\left\{
\matrix  \text{(not necessarily irreducible )
hypersurfaces of}\\\text{degree }d\text{  osculating   to order } 2d-1
  \text{
at }x\endmatrix \right\}
  \geq
\binom {a+d-1}d-1.
$$
\endproclaim

One possible definition of \lq\lq bending less\rq\rq\
would be that the  lower bounds for $d+1\leq k\leq 2d$
are attained (modulo the contribution of $I_{d-1}\circ V^*$).

Singular osculating hypersurfaces are the key to understanding
the projective geometry of non-complete intersections.
To explain why, for notational simplicity, assume 
$X$ is the intersection of hypersurfaces of degree $d$
(see \S 12 for the general case).

\proclaim{0.5.10 Proposition [L4, 1.1]} Let $X\subset\bppp V$ be a variety
such
that  $I_X= (I_d)$ (i.e. $I_X$ is generated by
$I_d$)
and $I_{d-1}=(0)$. Then the following are equivalent:

1.  $X$ is a complete intersection.

2. Every hypersurface of degree $d$ containing $X$ is smooth at all
$x\in\xsm$.

3.  Let $x\in \xsm$. Every
 hypersurface of degree $d$ containing $X$ is smooth at $x$.\endproclaim

Thus if $X$ is a complete interesection,
there can be no hypersurfaces
singular at any $x\in\xsm$
containing
$X$ that osculate to order infinity.

If $X$ is not a hypersurface, there are always singular hypersurfaces
in $I_X$. The proposition says that the singularities occur away from
$X$.

\subheading{0.5.11 Exercise} Verify   that the proposition
is true for the plane quartic, and to find the   quadric singular
at a point of the cubic curve.

\smallpagebreak

A precise version of
the complete intersection principle (0.5.6)
would be  that,
e.g., in the case of [L4, 1.1],
 no hypersurfaces of
degree $d$, singular at a general point $x\in X$, can osculate
to order $2d+1$ at $x$. If one understands \lq\lq bending\rq\rq\
in terms of genericity
conditions on differentials invariants, then
\lq\lq bending enough\rq\rq\ can be understood as the non-vanishing
of certain contractions
of differential invariants.
For example, in the case $d=2$,
it is sufficient that the   symmetrization map
$T^*\ot \hat{\ii}\ra S^3T^*$ is injective. (See \S 12.)

Taken together,  the smoothness principle and the complete
intersection principle indicate  that perhaps
varieties of small codimension   must be complete intersections.
Hartshorne has conjectured that if $a<\frac n2$ and $X$ is
smooth, then $X$ must be a  complete intersection. 
In fact, the two principles
were
 developed in an attempt to understand Hartshorne's conjecture
and other work motivated  by it
from the perspective of projective differential geometry.

The following result is proved using the rank restriction
theorem combined with a local study explained in \S 12:

\proclaim{0.5.11 Theorem [L6, 6.28]}
Let  $X^n\subset\bpp{n+a}$ be a variety and let 
 $x\in X$ be a general point.
  Let $b=\tdim X_{sing}$. (Set $b=-1$ if $X$ is smooth.)
If
$a<\frac{n-(b+1)+3}{3}
$
then any quadric osculating to order four at $x$ is smooth at $x$.
\endproclaim

By the discussion above, [L6, 6.28] implies

\proclaim{0.5.12 Corollary [L6, 6.29] }
Let  $X^n\subset\bpp{n+a}$ be a variety
with $I_X$ generated by quadrics.   Let $b=\tdim
X_{sing}$. (Set $b=-1$ if $X$ is smooth.)
If
$a<\frac{n-(b+1)+3}{3},
$
then $X$ is a complete intersection.
\endproclaim

\subheading{0.6 Monge equations}

[L6, 6.29] above poses the question: {\it How can one determine
if $I_X$ is generated in degree two?}. Before addressing the
question of determining if
$X$ is contained in quadric hypersurfaces, 
let's try an easier one: How can one
tell if
$X$ is contained in a hyperplane $H$?

$X\subset H$ if there exists a general point $x\in X$ such that
$n_H\in \tker\fff k{X,x}$ for some $k$, where $n_H\in N^*_xX$ is
a vector corresponding to $H$.

A corollary of the  rank restriction theorem is that if
$a<\frac {n-(b+1)}2+1$, and $x\in X$ is a general point,
then $III_{X,x}=0$. Thus:

\proclaim{0.6.1  Theorem  [L1]}Let $X^n\subset\bccc\bpp{n+a}$ be a
variety with $a<\frac
{n-(b+1)}2+1$ (where $b=\tdim X_{sing}$). Let $x\in X$ be a general point.
 If a hyperplane  $H$ osculates to
 order two at  $x$, then $X\subset H$.
\endproclaim

Thus in the situation of (0.6.1), to determine if $X$ is contained in a
hyperplane, two derivatives are sufficient.

Now consider the simplest case of a quadric hypersurface: Let
$X\subset\pp 2$ be a curve. How many derivatives does one need
to take to determine whether or not $X$ is a conic?

  To fix a plane conic,
one needs five points,
or equivalently, one point and four derivatives. Thus
to determine whether or not a given curve is a conic, it is necessary
to take five derivatives (One point and four derivatives
determines an osculating conic
$C_x(X)$, the fifth
derivative determines if $X=C_x(X)$). Plane conics
given as a graph
$y=f(x)$ are characterized by the classical {\it Monge equation}:
$
 ((y'')^{-\frac 23})'''=0.  
$
   See   [L4] for a derivation. 

A generalization of the classical Monge equation to determine if
an arbitrary variety is the intersection of quadrics would be
impossible, because no fixed number of derivatives would suffice
for all situations.  However, if 
$X$ is smooth and $a$ is small, one could hope
to have a fixed system. For example, if $X$ is a smooth hypersurface
and $n>1$,  Fubini [Fub] showed there exists a third order
system characterizing quadric hypersurfaces, so the situation here is
  better than for curves.

It turns out that if $X$ is of small codimension, but not a hypersurface,
then one needs five derivatives. In \S 12 I derive a fifth order
system of pde that I call the {\it generalized Monge system} that
characterizes intersections of quadrics when $a<\frac 13(n-(b+1)+3)$.

\subheading{0.7 Rigidity}

A general question related to step 5, of which (0.6) above is
a special case, is to know how many derivatives one needs to take
to recognize a given variety (or type of variety). In \S 13 I discuss
several
recognition questions centered around the results of [L9].

In [L9], I sharpen the result in
step five of my   proof of Zak's theorems on Severi
varieties. I show that  
   if $X$ is a variety and  $ x \in X$ a general point such
that $II_{X,x}$ is isomorphic to the second fundamental form
at a point of a Severi variety other than $v_2(\pp 2)$, then $X$ is the
corresponding Severi variety. (The result
is false when $n=2$, one must take third derivatives into account as well.)
The
case
$n=4$ had been conjectured by Griffiths and Harris in [GH].
  I also prove the following theorems:
  
\proclaim{0.7.1 Theorem [L9] }Let $X^{n+m}\subset\pp{nm+n+m-1+z}$,
$n,m\geq 2$, be an open subset of a variety   not contained
in a hyperplane with the second fundamental form of the
Segre $\pp n\times \pp m$ at   general points. Then $z=0$
and $X$ is an open
subset of the Segre $\pp n\times \pp m\subset\pp{nm+n+m-1}$.
\endproclaim

Note that the result is false if $n=m=1$.  

\proclaim{0.7.2 Theorem  [L9] }Let $X^{2(m-2)}\subset\pp{\binom m2
-1+z}$,
$m\geq 6$, be an open subset of a variety  not contained
in a hyperplane with the second fundamental form of the
Grassmanian $G(2,m)\subset\pp{\binom m2 -1}$ at   general points. Then
$z=0$
and $X$ is an open
subset of the Grassmanian.
\endproclaim

Note that the result is false if $m<5$.  

The varieties above are examples of {\it Hermitian symmetric spaces}.
A class of homogeneous varieties that resemble the Hermitian symmetric
spaces are the {\it adjoint varieties}. For example, the variety
of rank one  and traceless  $n\times n$ matrices  is an adjoint variety.
One might hope that the adjoint varieties are also determined by
their second fundamental forms. This turns out not to be the case:

\proclaim{0.7.3 Theorem [LM2]} There exist varieties
having the same second fundamental form as the adjoint varieties
at general points that  are not  adjoint varieties.\endproclaim

\smallpagebreak

\subheading{0.8 Gauss maps}

Consider again the cylinder:

$$
\ 
$$

The tangent directions to its rulings  are in $\tsingloc\ii$.
 In fact, the embedded tangent space is not
just constant to second order, but is constant all along the ruling, i.e.,
the ruling is a fiber of the   Gauss map.

  Another way to state
the result (0.3.3.3)   above is that   if  the  Gauss map of $X$
is degenerate, then $X$ is singular. A proof of this
statement and discussion of  varieties with degenerate Gauss mappings
is given in \S 5.

 Let $Y\subset\ppp V$ be a smooth variety and let $X=\tau (Y)$ be the
union of all embedded tangent lines to $Y$,
the {\it tangential variety} of $Y$. Then $X$ has a degenerate
Gauss map because its tangent space is constant along the tangent
lines of $Y$.

In [GH], Griffiths and Harris state that all varieties with degenerate
Gauss maps are constructed from cones and (generalized) tangential
varieties. They prove this statement when $\tdim X=2$. However
when $\tdim X=3$, there are already counter-examples to their announcement.
See \S 5 for some such examples.

\subheading{0.9 Dual varieties}

When $X$ is a hypersurface,   $\gamma (X)\subset\ppp V^*$ is the set
of hyperplanes tangent to $X$. A generalization  of this 
situation is as follows:

Let $X^n\subset\pp\na=\ppp V$. Define the
{\it dual variety} $X^*\subset\ppp V^*$ as the set of  hyperplanes tangent
to
$X$:
$$
X^*=\overline {\{ H\in \ppp V^*\mid \exists x\in \xsm\text{ such that
}
\tilde T_xX\subseteq H\} }.\tag 0.9.1 
$$
 One expects $X^*$ to be a hypersurface,
because there is an $(a-1)$-dimensional space of hyperplanes
tangent to each point and an $n$ dimensions   of points.
Let $\dstar = \dstar (X) = n+a-1-\tdim X^*$; the {\it dual defect} of $X$.
 Zak proved that if $X$ is smooth and
not contained in a hyperplane, then $\dstar\leq a-1$ (see \S 7).

A standard fact is that $(X^*)^*=X$, so we may think of $X^*$ as
a transform of $X$. In fact, the dual variety is the generalization
to algebraic geometry of the Legendre transform in classical mechanics.
Transforms are useful because they reorganize data 
in such a way  that information one   looks for
in an object  becomes more easily accessible
in its transform. For example,
 the   Fourier transform
exchanges global and local data. 
In my work with B. Ilic, [IL], our perspective was to view
$X^*$ as a transform of $X$. We describe ways in which the global
geometry of $X$ is reflected in the local geometry of $X^*$. 
We prove an inversion formula
that shows  that the second fundamental form at a point  
$H\in X^*$ contains information about all the points to which $H$ is
tangent.

Using the inversion formula, we show:

\proclaim{0.9.2 Theorem [IL]}Let $X^n\subset\pp\na$ be a smooth variety
with dual defect $\dstar$.
 If $H\in 
X^*_{sm}$, then $|II_{X^*,H}|$ is a system of quadrics of
projective dimension $\dstar$ and constant rank $n-\dstar$.
\endproclaim

Compare this result to an earlier result in [GH] (Bertini??)   that
if $X$ is any variety and $x\in X$ is a general point, then
$|II_{X,x}|$ is a system of bounded rank $n-\dstar$.

(0.9.2) led us to examine systems of quadrics of constant rank.
We were able
to solve an old question: What is the maximum dimension
of a system of quadrics of constant rank $r$ on $\bcc m$?  
(The answer when $r$ is odd was known classically to be one, see
\S 8). We showed:

\proclaim{0.9.3 Theorem [IL, 2.16]} If $r$ is even,   then
$$\tmax \{ \tdim(A) \mid A \subset S^2\bcc m \text{ is of constant rank } r
\}
 = m-r+1.$$
\endproclaim

[IL, 2.16] and [IL, 3.24] together furnish a new proof
of Zak's theorem
that  $\dstar\leq a-1$.  

\bigpagebreak

\subheading {0.10 Ruled and uniruled varietes}

 A variety 
$X\subset\ppp V$
is {\it uniruled} 
by $k$-planes if through
each
$x\in X$ there exists a $k$-plane contained in $X$ passing through $x$.
$X$ is {\it ruled} if it can be described as a fibration over
a base space that parametrizes the $k$ planes, e.g., 
$X$ is ruled if there is a unique   $k$-plane
through each point. The quadric hypersurface for $n>3$ is uniruled but not
ruled.

Consider again the hyperbola, a ruled surface.

$$
 \ \tag 0.10.1
$$

 The tangent directions to
the ruling are in $\tbaseloc\ii$. 
Directions in $\tbaseloc\ii$ of an arbitrary
variety $X$ are not usually tangent to linear spaces contained
in $X$. In [L7] I determine
  additional conditions that imply 
a variety
is   uniruled by $k$-planes.
In coordinates, 
the question is:  How many 
derivatives are needed  to determine if $X$ is
uniruled?  For example, given a
surface in $\Bbb P^3$,   there are always at least two tangent directions
in
$\tbaseloc\ii$ (so in this aspect, the hyperbola is
not special at all).
Such directions are called {\it
asymptotic directions} in classical differential 
geometry.  So two derivatives are   not enough to see if a
surface is ruled. A classical result [Blaschke ??] states that   three
derivatives are enough. Here is a  generalization of the classical result:

\proclaim{0.10.2 Theorem [L7]}
Let $X^n\subset\Bbb A^{\na}$ or $X^n\subset\pp\na$ be an open subset of a
variety  of a smooth (respectively analytic)
submanifold of an affine or projective space
such that at every point
(resp. at a general point) there is a line osculating to order $n+1$.
Then $X$ is uniruled by lines.

  There exist analytic 
open subsets of varieties $X^n\subset\Bbb A^{n+1}$ or $X^n\subset\pp{n+1}$ 
having a line osculating to order $n$ at every point that are not uniruled.
Over $\Bbb C$, every hypersurface has this property.   
\endproclaim

(0.10.2) also sharpens Z. Ran's dimension$+2$ secant lemma [R2]. Ran  
proves
that if a variety has lines osculating to order $n+1$ at each point, then
the union of the osculating lines is at most $n+1$ dimensional. The
above result shows that the union is in fact $n$ dimensional.

 In \S 11 I present some preliminary results 
regarding the integer $m=m(n,a,k)$
such that any open subset of a variety  $X^n\subset\pp\na$ having
$k$-planes
osculating to order $m$ at each point must be uniruled
by $k$ planes, but that there exist open subsets of   varieties with $k$
planes osculating to order $m-1$ that are not uniruled.
 
In a similar vein,   certain geometric situations when
osculation to order two is sufficient to imply containment of a linear
space are presented in \S 7.  

\subheading{0.11 Aside}

I spend alot of time studying systems of quadrics, so
people occasionally ask me questions about them. 
An amusing question was posed
to me  by M.
Kontsevich: Let $A= (a^i_j)$ be an orthogonal matrix with no
entries zero, form the {\it Hadamard inverse of $A$},
$Hinv(A)$, by $(Hinv(A))^i_j=1/{a^i_j}$. He conjectured
that the rank of $Hinv(A)$ was never two (but that perhaps all
other ranks were possible).

 A classical theorem  regarding systems of quadrics  is 
Castelnouvo's lemma, which states
that if $2n+3$ points lie on an $\binom n2$-dimensional system
of quadrics, then the points must lie on a rational normal curve.

  I resolve Kontsevich's conjecture by translating
it to a problem in algebraic geometry that turned out to be a variant
of Castelnuovo's lemma:  
If  a self-associated point set
of $2n+2$ points in $\pp n$ (which necessarily lie on
a $\binom n2$ dimensional space of quadrics) is such that
the standard
Cremona transform of $n+1$ of the points
 with respect to the simplex formed by
the other $n+1$ points is contained in a $\pp 2$, then in fact 
the $2n+2$ points
are colinear.  (see
[L8] for details).

\subheading{0.12 Self-similarity of geometric issues under change of scale}

When working with infinitesimal invariants, similar issues
often arise as when working globally.  To give three examples:
 
Ilic and I used the
  Lefschetz theorem on hyperplane sections
to prove the constant rank theorem, which
we used to obtain an infinitesimal proof of
Zak's bound on the dual defect. Zak's original proof   rests
on the Fulton-Hansen connectedness theorem, which in turn relies
on an extension of the Lefschtez theorem due to Goresky and MacPherson
(see \S 3 ).

A second example is that if $X$ is not a complete intersection, $I_X$ must
have
syzygies. (For example the quadrics defining the Segre have linear syzygies
that arise from expanding out the three by three minors.) If $X$ is not a
complete intersection, its differential invariants also must 
have syzygies. I hope  that these syzygies will be easier to study because
their existence can be related to degeneracies of auxiliary varieties.

A third example is the following proposition in [LM1],
where secant varieties arise in the study of infinitesimal geometry:

\proclaim{0.12.1 Proposition [LM1]} Let $X^n\subset\pp\na$ be an
  open subset of a variety   
and let $x\in X$ be a general point. Then
$\tbaseloc|\ff k|\supseteq\sigma_k(\tbaseloc\ii )$.\endproclaim

 $\sigma_k(X)$ denotes the closure of the union of all secant
$(k-1)$-planes
to $X$.

\subheading{0.13  Notation} 

 I  will use the following conventions for indices
$$
\align &0\leq B,C\leq n+a \\
& 1\leq \alpha ,\beta \leq n \\
%& n+1\leq \xi ,\eta \leq n+k \\
& n+1\leq \mu ,\nu \leq n+a   \endalign
 $$

Alternating products of  vectors will be denoted with a wedge ($\wedge$),
and
symmetric
products will not have any symbol (e.g. $\omega\circ\beta$ will be denoted
 $\omega\beta$).
$T_xX$ denotes the holomorphic tangent space to $X$ at $x$
and $\tilde T_xX$ the embedded tangent space. I often
 supress reference to $X$ and $x$, abbreviating the names of bundles,
e.g.
  $T$ should be read as $T_xX$,   $N$ as $N_xX$ etc...  
If $v\in \hat T_{[w]}X$, then I write
$\underline v:= v\ot w^* \tmod w\in T_{[w]}X$, where $v,w$ a part of
a basis of $V$ and $w^*$ denotes the dual basis vector to $w$.
  $\{ e_i \}$ means the span of the vectors $e_i$ over the index range $i$.
If $Y\subset\bpp m$ then $\hat Y\subset\Bbb C^{m+1}$ will
be used to  denote the cone over
$Y$
 (with the exception that  the cone over the embedded tangent space
${\tilde T}$ will be denoted
$\hat T$).
If $A\in \Bbb C^{n+a+1}$, its projection to $\bpp{n+a}$ will be denoted
$[A]$.
If $V$ is a vector space and $W$ a subspace, and $(e_1\hd e_n )$ a basis of
$V$
such that $\{e_1\hd e_p \}=W$, I write
$\{e_{p+1}\hd e_n\}$ mod$W$ to denote the space $V/W$.
For vector subspaces $W\subset V$, I  use the notation
$W\upperp  \subset V^*$ for the annihilator of $W$ in $V^*$.
I will use the summation convention throughout (i.e. repeated indices are
to
be summed over).
$\frak S_{\alpha\beta\gamma}$ denotes cylic summation over
the fixed indices $\alpha\beta\gamma$.
In general, $X$ will denote a variety,  $\xsm$ its smooth points, and
$X_{sing}$ its singular points.  $\Bbb C\bpp k$ will  be denoted $\bpp k$.
$\bii k X$ is the $k$-th fundamental form of $X$.
I often denote $\bii 2X$ by $II$ and $\bii 3X$ by $III$.
$F_k=F_k^{\mu}e_{\mu}\tmod \hat T$ is the differential invariant called
  the  {\it $(k-2)$-nd variation of $II$}.
By a {\it general point}
$x\in X$ I mean  a smooth point of $X$ such that all the discrete
information
in the
differential invariants of $X$ is locally constant.  The nongeneral points
of
$X$ are
a codimension one subset of $X$.

\heading \S 1. Examples of homogeneous varieties and their uses\endheading

A smooth projective variety $X\subset\ppp V$ is {\it homogeneous} if
$X$ is the closed orbit of a 
(complex semi-simple) Lie group $G$ acting on  $\ppp V$.
Such $X$ can be described intrinsically as $X=G/P$ where $P$
is a subgroup of $G$, called a {\it parabolic} subgroup.

One can reduce to the case of 
varieties that are orbits of simple groups
(those whose Lie algebras have no nontrivial ideals),
as others are just products of such. There are three
simple groups occuring in series,
$SL(V,\Omega)$, $Sp(V,\omega)$, $O(V,Q)$,
where respectively $\Omega\in \La m V^*,\omega\in\La 2V^*, Q\in
S^2V^*$ are nondegenerate elements and the groups are the subgroups
of $GL(V)$ preserving the forms ($\tdim V=m$). Actually $O(V,Q)$ is not
connected or simply
connected, its connected
component of the identity is called
$SO(V,Q)$,  its simply connected double cover
is called $Spin(V,Q)$. For $Sp(V,\omega)$, the
dimension of
$V$ must be even to have a nondegenerate two-form. Since the behaviour of
nondegenerate quadratic forms is quite different in even and odd dimensions
(e.g. every rotation in $\brr{2n+1}$ has a fixed axis), 
for $O(V,Q)$, the
case where  $\tdim V$ is even and odd are   considered as different groups.

In addition to these groups, there are five {\it exceptional} groups,
which are called $G_2,F_4$,
$E_6$,
$E_7,
E_8$. In what follows I will describe
all but the last two. (Actually $G_2$ is quite easy to describe,
it is the subgroup  of $GL(7,\BC)$ preserving a nondegenerate
(generic)
 element of $\La 3\bcc 7$.)

\subheading{1.1 More on Segre varieties}

Given vector spaces $W_1\hd W_r$, define an embedding
$X=Seg(\ppp W_1\times\hdots\times \ppp W_r)\subset \ppp (W_1\ot\hdots\ot
W_r)$ by $[w_1]\times \hdots \times [w_r]
\mapsto [w_1\ot\hdots\ot w_r]$. $X$ is called  the {\it Segre}
embedding of 
$\ppp W_1\times\hdots\times \ppp W_r$. (In the language of [GKZ],
these are the  \lq\lq rank one multidimensional matrices\rq\rq .)

One can form the product of varieties. Unlike
the affine case where the product of $M\subset\Bbb A^m$ and
$N\subset\Bbb A^n$ is naturally $M\times N\subset\Bbb A^{m+n}$,
the product
to two projective varieties  naturally occurs as a subvariety of 
the Segre of their ambient spaces.
If $X\subset\ppp V$, $Y\subset\ppp W$,  
$Seg(X\times Y)\subset \ppp (V\ot W)$ is the natural
  projective embedding of $X\times Y$.

\smallpagebreak
 If $W_1=W_2$, 
one can consider instead of arbitrary matrices, symmetric or
skew-symmetric matrices of minimal rank and their generalizations.
We do so in the following two examples.

\subheading{1.2 Grassmanians}

Let $V$ be a vector space, choose an identification
$V\simeq \bcc m$ so  $\La 2V$ is identified with
the $m\times m$ skew symmetric matrices. Let $G(2,V)\subset\ppp( \La 2
V)$ denote the 
projectivization of the rank two matrices (i.e. the
matrices of minimal rank, as the rank of a skew-symmetric
matrix is always even). $G(2,V)$ is called the {\it Grassmanian
of two planes in $V$}. 

Let
$\ee 1\hd\ee n$ be a basis of $V$,
then $\{\ee i\ww\ee j\mid i\leq j\}$, is a
basis of $\La 2 V$, and we may think of $\ee i\ww\ee j$ as the
skew symmetric matrix with $1$ in the $(i,j)$-th slot,
$-1$ in the $(j,i)$-th slot and zero elsewhere. Any
$E\in G(2,V)$ can be written $E= 
v\ww w=v^iw^j(\ee i\ww\ee j)$,
where $v=v^i\ee i, w= w^j\ee j$. This gives two interpretations of $G(2,V)$;
as the decomposable elements in $\La 2V$
(i.e. elements of the form $v\ww w$, with
$v,w\in V$), and as the set of two
planes through the origin in $V$ ($E$ is the two-plane spanned by $v$ and
$w$).

Generalizing, let $G(k,V)\subset\ppp \La k V$ denote the set of $k$-planes
through the origin in $V$, or equivalently the decomposable $k$-vectors
(those that can be written $v^1\ww\hdots\ww v^k$,
with  each $v^j\in V$). We may also
think of $G(k,V)$ as the space of $\pp{k-1}$'s in $\ppp V$. When
we use this perspective, we will write $G(k,V)=\Bbb G(k-1,\ppp V)$.

The tangent space to any manifold at a point is a vector space,
the tangent space to $G(k,V)$ is a vector space with additional
structure, namely $T_EG(k,V)\simeq \thom (E,V/E)= E^*\ot (V/E)$.
 To see this take a curve $E(t)= v^1(t)\ww\hdots\ww v^k(t)$
and differentiate at $t=0$. 
The tangent space to any homogeneous space is always a vector
space with additional structure, and this additional structure
can be deduced either intrinsically or extrinsically.

In particular for the case of $G(1,V)=\ppp V$, 
  $T_x\ppp V= \hat x^*\ot (V/\hat x)$. (We use the notation
that for $Y\subset\ppp V$, $\hat Y\subset V$ is the corresponding cone.)
If $X\subset\ppp V$ is a subvariety, the (intrinsic holomorphic) tangent
space to $X$ inherits this additional structure. Letting
$\tilde T_xX$ denote the embedded tangent space and
$\hat T_xX= \hat{\tilde T}_xX$,  then   $T_xX= \hat x^*\ot 
(\hat T_xX/\hat x)$.
Since we are not using a metric, the normal bundle is just a quotient
bundle, $N_xX= T_x\ppp V/T_xX = \hat x^*\ot (V/\hat T_xX)$. Traditionally
the line bundle with fiber $\hat x^{\ot k}$ is denoted $\Cal O_{\ppp
V}(-k)$
and if $E$ is any vector bundle, $E(k)=E\ot \Cal O(k)$. 

 Note that $G(k,V)=G(n-k,V^*)$ as specifying
a $k$ plane $E\subset V$ is equivalent to specifying its
annhilator $E\upperp\subset V^*$. Of particular
importance is  $G(n-1,V)=\ppp V^*$ the {\it dual
projective space}, where points of $\ppp V^*$ 
correspond to  hyperplanes in $\ppp V$.

If $G\subset SL(V)$ is a group preserving additional structure
(e.g. a quadratic form $Q$ or symplectic form $\omega$) one
can define the corresponding null Grassmanians, e.g.
$$
\align
G_{Q-null}(k,V)&:=\{ E\in G(k,V)\mid Q(v,w)=0\forall v,w\in E\}\tag 1.2.1\\
G_{\omega-null}(k,V)&:=\{ E\in G(k,V)\mid \omega(v,w)=0\forall v,w\in E\}
\endalign
$$
The $Q$-null Grassmanians are naturally embedded
in $\ppp\La kV$. 
In the  case $\tdim V=2m$
and $k=m$,   $G_{Q-null}(m,2m)$ has two isomorphic components.
The components are called the {\it Spinor varieties} $\Bbb S_m$
and each $\Bbb S_m$ embedds into a smaller projective space
which  I will
describe after explaining Clifford algebras.

\subheading{1.2.2 Exercise} Show that $G_{\omega -null}(k,V)$ naturally
embedds
into a linear subspace of $\ppp ( \La k V)$. Calculate its tangent space
at a point.

\subheading{1.3 Veroneses}

Let $S^2V$ denote the symmetric matrices and let
$v_2(\ppp V)\subset\ppp (S^2V)$ denote the
projectivization of the rank one elements.  
$v_2(\ppp V)$ is the image of $\ppp V$ under the injective mapping
$$
\align v_2 :\ppp V&\ra \ppp S^2 V\tag 1.3.1 \\
[x ]&\mapsto [x\circ x]\endalign
$$
The {\it $d$-th Veronese
embedding of $\ppp V$}, $v_d(\ppp V)\subset\ppp S^d V$
is defined 
by $v_d(x)= x^d=x\circ\hdots\circ x$. Given $X\subset\ppp V$, we can
consider
the Veronese re-embeddings of $X$, $\vdx\subset\ppp (S^dV)$,
which will turn out to be useful in our study of complete intersections.
 Note that if $Z$ is a hypersurface of degree $d$,
$v_d(Z)= H_{ Z}\cap v_d(\ppp V)$, where $H_Z $ is the
hyperplane associated to the
equation of $Z$.  
  
\subheading{1.4 Division algebras and
the spinor variety $\Bbb
S_5$}

There are four division algebras over $\brr{}$: 
 $\Bbb R$, $\Bbb C$, $\Bbb H$, and $\Bbb O$ (where
$\Bbb O$ denotes the {\it octonians}, or {\it Cayley numbers}). The
octonians
are similar to the quaternions. If one thinks of $u\in\Bbb H$
as $u= u^0+ u^1\epsilon_1+ u^2\epsilon_2+ u^3\epsilon_3 $
where
$u^i\in\Bbb R$ and
 the $\epsilon_i$ satisfy $\epsilon_i^2=-1$ and the
following multiplication table:
$$
\ \tag 1.4.1 
$$

 (to be read...), then given $u\in\Bbb O$,
write $u= u^0+ u^1\epsilon_1+\hdots + u^7\epsilon_7$ where
$\epsilon_j^2=-1$
and the $\epsilon_j$ satisfy the following multiplication table:
$$
\
$$

  Let $\Bbb A =\Bbb R,\Bbb C, \Bbb H$ or $\Bbb O$ and let $\Bbb A = \Bbb
A_{\Bbb R}\ot_{\Bbb R} \Bbb C$. 
The only spinor variety we will have immediate need of
is $\Bbb S^{10}_5\subset\pp{15}$ which may be described as follows:
Let $\bcc{16}=\Bbb O^2\ot_{\Bbb R}\Bbb C$ have octonionic coordinates
$u,v$.
$\Bbb S_5$ is defined by the equations $u\overline u=0, v\overline v=0,
u\overline v=0$ where the last is eight 
equations. 

\subheading{1.4.2 Exercise} what are the corresponding varieties
for the other division algebras?)

\smallpagebreak

Let  $Aut(\Bbb A):=
\{ g\in GL(\Bbb A) \mid (gu)(gv)=g(uv)\,
\forall u,v\in\Bbb A\}$.
$Aut(\Bbb A)\subset GL(\text{Im}\Bbb A)$ is
respectively   $\{Id\}, \Bbb Z_2, Sl_2\Bbb C, G_2$, providing a second
definition of $G_2$.

\subheading{1.5 Clifford algebras}

 Let  
$(V,Q)$ be as above. Given any linear subspace $L\subset V$,
one can define $L\upperpq\subset V$,
by $L\upperpq = \{ w\in V\mid Q(v,w)=0\, \forall v\in L\}$.
If $Q|_L$ is nondegenerate
$V= L\oplus L\upperpq$. In this case, for all $v\in V$, we may
write $v= v_1+v_2$ with $v_1\in L$, $v_2\in L\upperpq$. We may
define the {\it reflection} of $v$ in $L$ by $refl_L(v)= v_1-v_2$.

$$
\ 
$$

Recall that $O(V,Q)$ is the subgroup of $GL(V)$ preserving $Q$,
and $SO(V,Q)$ is the  component of  $O(V,Q)$
containing the identity.

\proclaim{1.5.1 Theorem, Cartan-Diedonne  (see  [Hv])}
$O(V,Q)$ is the group generated by   reflections in lines.
$SO(V,Q)$ is the group generated by even numbers of reflections.
More precisely, $O(V,Q)=
\{ refl_{l_1}\circ\hdots\circ refl_{l_k}\mid l_j\in  V\}$
and we may assume $k\leq n$. Similarly for $SO(V,Q)$, only $k$
must be even.
\endproclaim

To define $Spin(V,Q)$,
the connected and
simply connected
group corresponding to $O(V,Q)$, 
 we will need to
generalize the notion of a reflection. Let 
$  \La\bullet V= \ot V/\{ x\ot y-y\ot x\}$ be the exterior algebra.
The {\it exterior product} in
$\La\bullet V$, $(x,y)\mapsto x\ww y$ may be interpreted as follows:
Let $\hat G(i,V)\subset\La i V$ denote the cone over the Grassmanian.
If $x\in \hat G(i,V)$, $y\in \hat G(j, V)$,
  then $x\ww y\in \hat G(i+j, V)\subset \La{i+j}V$
represents the $i+j$-plane spanned by $x$ and $y$.
If $x\in V$, then $x\ww y$ is analogous to the component of $x$
in $[y]\upperpq$. If
$y$ has unit length, then $||x\ww y||_Q= ||proj_{y\upperpq}(x)||_Q$.
Note that we do  not need $Q$ to define $x\ww y$.

Let $x\intprod y$ be defined to be the $Q$-adjoint
of $x\ww y$, that is
$Q(x\ww y, z)= Q(y, x\ww z)$ for all $x,y,z\in \La\bullet V$. If $x\in V$
and $y\in\hat G(j,V)$, then $x\intprod y$ is analogous to the component
of $x$ in $\ppp \{ y\}$.
Note that  if
$y$ has unit length, then $||x\intprod y||_Q= ||proj_{y}(x)||_Q$.

Consider, for $x\in V$, $y\in\La \bullet V$,
$$
x\circ y:= x\intprod y - x\ww y
$$
which can be thought of the
\lq\lq generalized reflection\rq\rq\
of $x$ in $y$. We can, by linearity,
extend $\circ$ to a mulplication on $\La\bullet V$.

\subheading{1.5.2 Definition} Let $V$ be a vector space with a
nondegenerate
quadratic form $Q$. Let $Clifford(V,Q)=
Cl(V,Q):= (\La\bullet V,\circ)$,
the {\it Clifford algebra} of $(V,Q)$. 

As is usual with algebras formed from vector spaces, we have a 
fundamental classifying/universality lemma:

\proclaim{1.5.3 Fundamental Lemma of Clifford algebras (see   [Hv])}Let
$V$ be a vector space with a nondegenerate quadratic form $Q$ and let $\Cal
A$
be an associative algebra with unit. If $\phi : V\ra \Cal A$ is a mapping
such that for all $x,y\in V$
$$
\phi(x)\phi (y)+ \phi (y)\phi(x)= 2Q(x,y)Id_{\Cal A}
$$
then $\phi$ has a unique extension to an algebra mapping $\tilde\phi :
Cl  (V,Q)\ra\Cal A$.\endproclaim

\subheading{1.5.4  Exercise} Verify that $\Cal A = (\Lambda^\bullet V, \circ )$
satisfies the hypotheses of the fundamental lemma.

\smallpagebreak

In $Cl(V,Q)=(\Lambda^\bullet V, \circ )$,  
the degree of a form
is no longer well defined, but there is still a notion of parity. Let
$Cl^{even}(V,Q),Cl^{odd}(V,Q)\subset Cl (V,Q)$ denote the corresponding
even and odd subspaces. As vector spaces, $Cl^{even}(V,Q)=\La{even}V$,
$Cl^{odd}(V,Q)=\La{odd}V$.

Let $Cl^*(V,Q)\subset Clifford(V,Q)$ denote the invertible elements.

\subheading{1.5.5 Definition}
$$
\align
Pin(V,Q)&:= \{ a\in Cl^*(V,Q)\mid a= u_1\circ\hdots\circ u_r,
u_j\in V, Q(u_j,u_j)=1\}\\
Spin(V,Q)&:=\{ a\in Pin (V,Q)\mid r\text{ is even }\}.\endalign
$$

 From now on, assume $\tdim V=2m$.

We have a description of $Spin (V,Q)$ as an abstract group
and as a group acting on a vector space, namely $\La\bullet V$. This action
is clearly reducible because it preserves parity, but
in fact it   reduces further. Fix $U^m\subset V$, a
maximal null subspace. It is easy to verify that in fact $Spin(V,Q)$
preserves
$\La{even}U\subset \La\bullet V$ and, slightly more difficult to verify,
but
also true, that this action is irreducible. (Traditionally, $\La{even}U$
is denoted $\Cal S^+$ and called the space of {\it positive spinors}.)

Given $a= u_1\circ\hdots\circ u_r\in Clifford(V,Q)$, let
$\tilde a= (-1)^ru_r\circ\hdots\circ u_1$. This is a well defined
involution.
Using it, we may define an action  $\rho$ of $Spin(V,Q)$ on $V$, by
$\rho(a)v:=av\tilde a$, and it is easy to verify that $\rho$ is  
a $2:1$ map $Spin(V,Q)\ra SO(V,Q)$.

One may describe    $\Bbb
S_m\subset\ppp
\Cal S^+$,   as   the null
$m$-planes 
$E$ such that $\tdim(E\cap U)-n\equiv 1\tmod 2$  as
  follows:

Define a map
$$
\align
V\times\La{even}U&\ra \La{odd}U\tag 1.5.6 \\
(v,\alpha)&
\mapsto
proj_{\La{odd}U}(v\circ \alpha)\endalign
$$
($v\circ \alpha \in \La{odd}V$ and then we project it to $\La{odd}U$).
 In other words, for $\alpha\in\La{even}U$, we have a
map $L_{\alpha}: V\ra \La{odd} U$. In particular, for $1\in\La 0 U$,
$\tker L_1= U'$,
where $U'\subset V$ is the unique null $m$-plane such that
$V= U\oplus U'$. Since the map is
$Spin(V,Q)$ equivariant, for $g\in Spin (V,Q)$,
$\tker L_g= \tker L_{g\circ 1}= \rho (g)(\tker L_1)= gU'\tilde g$.
In summary:
\proclaim{1.5.7 Proposition/Definition}
$$\align
\Bbb S_m(V,Q,[U])&= \{ 
\pp{m-1}\subset Q\mid \tdim  (\pp{m-1}\cap \ppp U)-m\equiv 1\tmod 2
 \}\\
&=\ppp \{\text{the }Spin(V,Q)\text{ orbit of }1\in \La{even}U\}
\subset\ppp (\La{even}U). \endalign
$$
\endproclaim
In \S 4 I show that $  T_{[1]}\Bbb S_m\simeq   \La 2U $;
in particular, $\tdim \Bbb S_m=\binom n2$.  

\subheading{1.6 The Severi varieties}

Here is a different generalization of the Veronese $v_2(\ppp V)$:
  Let
$\Bbb A$ denote the complexification of a division algebra as above, and
let
$\Cal H_{\Bbb R}$ denote the
$\Bbb A_{\Bbb R}$-Hermitian  forms on $\Bbb A_{\Bbb R}^3$,  i.e., the
$3\times 3$
$\Bbb A_{\Bbb R}$-Hermitian matrices.   If $x\in\Cal H_{\Bbb R}$, then
we may write
$$
x = \pmatrix r_1 & \bar u_1 & \bar u_2 \\
u_1 & r_2 & \bar u_3 \\ u_2 & u_3 & r_3 \endpmatrix \ \ r_i\in \Bbb R , \
u_i\in \Bbb A_{\Bbb R}.\tag 1.6.1
$$
Let $\Cal H = \Cal H_{\Bbb R}\ot_{\Bbb R} \Bbb C$.

\subheading{1.6.2 Exercise}Verify that the notion of $x^2$ and $x^3$ make
sense (one needs to use the Moufang identites  in the case of the
octonians,
see [Hv]).

\smallpagebreak

Define a cubic form $det$ on $\Cal H$ by
$$ det(x) := \frac{1}{6}((\text{trace}(x))^3 + 2\text{trace}(x^3) - 3
\text{trace}(x)\text{trace}(x^2)).
\tag 1.6.3
$$
  $det$
is just the usual determinant of a $3\times 3$ matrix when $\Bbb A = \Bbb
C$.
When
$\Bbb A_{\Bbb R} = \Bbb O$ one cannot define $det$ for $4\times 4$ or
larger
matrices.   

Now, considering $\Cal H$ as a vector space over $\Bbb C$, let $G$ be the
subgroup of $GL(\Cal H,\Bbb C )$ preserving $det$, i.e., define
$$ G:= \{g\in GL(\Cal H,\Bbb C) | det(gx) = det(x)\ 
\forall x\in\Cal H\}.
$$
The   respective groups are:
$$\align \Bbb A_{\Bbb R} =& \ \ \ \ \ \ \ \ G = \tag 1.6.4\\
 &\Bbb R \ \ \ \ \ \ \ \ \ \  \ Sl(3,\Bbb R)^{\Bbb C}= Sl(3,\Bbb C)\\
  &\Bbb C \ \ \ \ \ \ \ \ \ \  \ Sl(3,\Bbb C)^{\Bbb C}= Sl(3,\Bbb C)\times
Sl(3,\Bbb C)\\
   &\Bbb H \ \ \ \ \ \ \ \ \ \  \ Sl(3,\Bbb H)^{\Bbb C}= Sl(6,\Bbb C)\\
    &\Bbb O \ \  \ \ \ \ \ \ \ \  \ \lq Sl(3,\Bbb O)^{\Bbb C}\text{\rq} =
E_6
\endalign
$$
 where I write    \lq$Sl(3,\Bbb O)^{\Bbb C}$\rq\  merely
 to be suggestive. We take the above as the
definition of $E_6$. The group $F_4$
is defined to be the subgroup of $E_6$ preserving the quadratic
form $Q(x,x)= \ttrace (x^2)$ where $xy$ is the usual matrix multiplication
and one must again use the Moufang identities to be sure $Q$ is well
defined.

\subheading{1.6.5 Exercise} Show that the action of $F_4$ on $\Cal H$
preserves the line $\BC\{ Id\}$ and $\Cal H_0:=\{ x\in \Cal H\mid
\ttrace (x)=0\}$.
(In fact, $F_4$ acts irreducibly on both factors.)

\smallpagebreak

$det$ tells us which elements of $\Cal H$ are of less than full rank. One
can
also unambiguously define a notion of being rank one; either by taking
$2\times
2$ minors or by noting  that under the $G$ action each $x\in \Cal H$ is
diagonalizable and one can take as the rank of $x$ the number of nonzero
elements in the diagonalization of $x$.

Let
$$
X := \pp {}\{\text{ rank one elements of }\Cal H \} = \pp {}\{ G\text{
orbit of
any rank one matrix}\}.
$$
  $X = (\Bbb A_{\Bbb R}\pp 2)^{\Bbb C}$,
that is, the complexification of the space of $\Bbb A_{\Bbb R}$-lines
in
$\Bbb A_{\Bbb R}^3$.
The four varieties $X\subset\ppp\Cal H$ are called the {\it Severi
varieties}.
We have already seen the first three, they are
$v_2(\pp 2)\subset\pp 5, Seg(\pp 2\times\pp 2)\subset \pp 8$,
and $G(2,6)\subset\pp {14}$.  

\subheading{1.6.6 Exercise} What are the analogous groups and varieties
if one takes instead the $2\times 2$ $\Bbb A$-Hermitian forms?

\subheading{1.6.7 Exercise} Let $\BO\pp 2_0=\Cal H_0\cap\BO\pp 2$.
Show that
$$
\BO\pp 2_0=\ppp \{ x\in\Cal H_0\mid x^2=0\}
$$
and deduce that it is a homogeneous space of $F_4$.

\heading{\S 2. Constructing new varieties from old}\endheading

As mentioned in the introduction,
one way to study   subtle properties of a
variety
$X^n\subset\ppp V$ is to study coarse properties of auxilliary varieties
one
constructs from $X$. We will be primarily concerned with studying the
dimensions of auxilliary varieties. Here are some constructions:

\subheading{2.1 The Gauss images of $X$}

We continue the notations $V=\bcc{\na +1}$ and $X^n\subset\ppp V$
a subvariety.

The {\it Gauss map}
of $X\subset\ppp V$ is defined at smooth points by:
$$
\align \gamma :\xsm &\ra \Bbb G(n ,\ppp V)\\
x&\mapsto \tilde T_xX.\endalign
$$
$\gamma$ can be completed to a rational map
$X\dashrightarrow\Bbb G(n ,\ppp V)$ whose image
we denote $\gamma (X)$, the {\it Gauss image} of $X$.
We will say $\gamma$ is {\it degenerate} if $\tdim \gamma (X)<\tdim X$.

More generally, for each $x\in \xsm$ we can consider the
$(n+k)$-dimensional projective spaces that contain $\tilde T_xX$,
and the resulting submanifold of $\Bbb G(n+k,\ppp V)$. The
closure of these is a variety which we denote $\gamma_k(X)$
and call the {\it $k$-th Gauss image of $X$}. The notation is such that
$\gamma_0(X)=\gamma(X)$.
Of particular interest is the  {\it dual variety}
$X^*=\gamma_{a-1}(X)\subset\Bbb G(n+a-1,\ppp V)=\ppp V^*$.  One expects 
$\tdim \gamma_k(X)$ to be $n+k(a-k)$ because there are $n$ dimensions of
points on $X$ and a  $k(a-k)$-dimensional
space of $(n+k)$-planes tangent to each smooth point. We  
say $\gamma_k(X)$ is {\it degenerate} if it fails to be of
the expected dimension. We will be particularly interested in
the degeneracy of the dual variety, we
  let $\dstar= n+a-1 -\tdim X^*$ denote the {\it dual defect} of $X$.

To better understand the higher Gauss images, it is useful to
use   a standard construction  in algebraic geometry,
the {\it incidence correspondence}.

Let 
$$
\ci= \ci (X,\Bbb G(n+k,\ppp V))= 
\overline{\{ (x,L) \mid
\tilde T_xX\subseteq L\} }.\tag 2.1.1 
$$
 We have the following picture
$$
\align &\ci_k \tag
2.1.2\\
& \\
& \\
X \ \ \ \ \ & \ \  \ \ \ \Bbb G(n+k,\ppp V)  \endalign
$$
 **maps are $\pi_k,\rho_k$***

By definition $\gamma_k(X)=\rho_k\pi_k\inv (X)$.
We will often use the notation $Y_L= \pi(\rho\inv L)$.
In  the most important case of $\ci_{a-1}$ we write
$\ci=\ci_{a-1}$ etc...

\smallpagebreak

With an eye towards Morse theory
(which is used in the proof of the
Lefschetz theorems, see \S 3), one can define $X^*$ as follows:
Consider the pairing
$<,> :\hat X\times V^*\ra \Bbb C$. For $p\in \hat X$, consider the function
$$
\align 
<p,\cdot >:V^* &\ra \Bbb C\tag 2.1.3
\\q&\mapsto <p,q>.\endalign
$$
Then $\pi\inv (x)=\ppp (\text{critical point of}<p_x,\cdot >)$,
where $p_x\in \hat x$ and
$$
X^*=\overline{
\cup_{x\in\xsm}\ppp (\text{critical points of } <p_x,\cdot >) }\tag 2.1.4
$$

$\tilde T_xX\subseteq H$ if and only if $X\cap H$ is singular at $x$, so
  another definition of 
$X^*$ is the union of the hyperplanes $H$ such that $X\cap H$ is
singular.

\subheading{2.2 Secant varieties and joins}

 Let $Y,Z\subset\ppp V$ be two varieties.
We define the {\it join} of $Y$ and $Z$,
$$
S(Y,Z)=\overline {\cup_{y\in Y, z\in Z}\BP^1_{yz}}
$$
where the closure is not necessary if the two varieties do not
intersect. This generalizes the secant variety which is the
case $Y=Z$. We can similarly form the join of $k$ varieties
$Y_1\hd Y_k$, 
$$
S(Y_1\hd Y_k)=\overline{\cup_{y_j\in Y_j }\BP^{k-1}_{y_1\hd y_k}} \tag
2.2.1
$$

A special case of this construction is if $Z=L$ is a linear
space. Then $S(Y,L)$ is   a cone over $Y$.  

An important result about joins is the following:

\proclaim{2.2.2 Terracini Lemma} Let $Y,Z\in\ppp V$ be
 varieties and let $x\in \BP^1_{yz}$. Then
$$
\hat T_xS(Y,Z)\supseteq \hat T_yY +\hat T_zZ.
$$ 
Moreover,  if $x$ is a general point of $S(Y,Z)$, then equality
holds.  
\endproclaim

\subheading{2.2.3 Exercise} Prove the inclusion part of
Terracini's lemma. Hint: consider the two curves
$p+q(t)$ and $p(t)+q$, where $p\in Y$, $q\in Z$ and differentiate.

\subheading{2.2.4 Special cases}

 If $Z=L$ is a linear space,
then $\tilde T_xS(Y,L)$ contains $L$ for all $x\in S(Y,L)_{sm}$ and
$S(Y,L)$ is the cone over $L$.

If 
   $Y_j=Y$ for all $j$, we call $S(Y\hd Y)$ the
the {\it $k$-th secant variety of $Y$} and use
the notation  $ \sigma_k(Y)$.   The notation is such that 
$\sigma_1(Y)=Y$.
We often denote $\sigma_2(Y)$ by $\sigma (Y)$. Note that one expects $\tdim
\sigma_k(Y)=
\tmin\{ nk +k-1, n+a\}$
as there are $nk$ dimensions worth of picking $k$ points on $Y$
and their span is a $\pp{k-1}$.
  
\smallpagebreak

A related notion  is
the {\it $k$-th multi-secant variety of $X$}, $MS_k(X)$ is defined
to be the closure of the union of all lines in $\ppp V$ containing
$k$ points of $X$.   An essential
observation regarding them, due to Severi, is
that if $X$ is contained in a hypersurface of degree $d$, then
$MS_k(X)\subset X$ for all $k>d$ (and thus is empty if $X$ contains
no lines). See [R2] for one use of this observation.

\subheading{2.3 Tangential varieties}

If $X$ is smooth, the tangential variety of $X$,
$\tau (X)\subset\ppp V$, is simply the union
of all the embedded tangent spaces. When $X$ is not smooth, there
are several possible notions one could use to define tangent
spaces (see [L5..]).
The notion that turns out to be useful is the 
  {\it tangent star},
$T^{\star}_xX\subset\ppp V$. Intuitively,
$T^{\star}_xX$ is the limit of secant lines.
More precisely, let $x\in
X$,
$\Bbb P^1_*$ is a line  
in $T^{\star}_xX$ if there exist smooth
curves
$p(t),q(t)$ on $X$ such that $p(0)=q(0)=x$ and
$\Bbb P^1_* = \text{lim}_{t\rightarrow 0} \Bbb P^1_{pq}$. 
$T^{\star}_xX$ is the union of all $\pp 1_*$'s at $x$ and we define
the {\it tangential variety of $X$}, by
$\tx= \cup_{x\in X}T^{\star}_xX$. One can define higher
tangential varieties to be the union of higher osculating spaces,
 i.e.,
$\tau_k(X):=\cup T^{\star (k)}_xX$, where at smooth points
$T^{\star (k)}_xX=\tilde T^{  (k)}_xX$ is the usual osculating space
and at singular points one takes the union of $\pp{k-1}$'s that are
limits of $k$ points on $X$ moving towards $x$.
Tangential varieties are sometimes called {\it developpable varieties}.

\heading\S 3.   Topology and consequences\endheading

Three basic theorems about the topology of subvarieties of projective
space are Bertini's theorem, Lefshetz's theorem
and Bezout's theorem. 

  Bezout's theorem states
that  varieties of complementary dimension
must intersect (generalizing the fundamental
theorem of algebra). We have already seen that
Bezout's theorem  has consequences for
the inifinitesimal   geometry of subvarieties. We will see the same
is true of the Bertini and Lefschetz theorems. The Bertini and
Lefschetz theorems and  their   variants 
 are also quite useful at the infinitesemal
level to prove results about systems of quadrics.

\proclaim{3.1 Bertini's Theorem} Let $A\subset S^dV^*$
be a linear subspace. Let $P\in A$
be a general element and let $Z_P\subset\ppp V$ be the hypersurface
of degree $d$ that
it determines. Then $(Z_P)_{sing}\subseteq\tbaseloc (A)$.
\endproclaim

Bertini's theorem can be understood as a quantitative version
of Sard's theorem for polynomials. Its proof is elementary,
see e.g. [GH2].

One of the most important topological theorems in algebraic geometry
is the Lefschetz theorem. Its simplest form is:

\proclaim{3.2 Lefschetz theorem, version 1}
Let $X^n\subset\pp\na$ be
a smooth variety and let $H\subset\pp\na$
be a hyperplane, then the restriction map on cohomology:
$$
H^i(X,\Bbb Z)\ra H^i(X\cap H, \Bbb Z)
$$
is an isomorphism for $i<n-1$ and injective for $i=n-1$.
\endproclaim

There are two standard proofs of the Lefschetz theorem, one using
Morse theory (see [Milnor]), and the other using harmonic differential
forms
and the Hodge theorem (see [GH2]). In both proofs, an essential point
is that when one writes a $n\times n$ 
Hermitian matrix with complex
entries as a $2n\times 2n$ matrix with real entries, the
eigenvalues of the new
matrix will occur in pairs $\lambda, -\lambda$.

Note that if $X$ is smooth, then by Poincar\'e duality  the cohomology
of $X\cap H$ is determined
 in all but the three middle dimensions.

Now let $Z\subset\pp\na$ be a smooth  hypersurface of degree $d$
and consider
$X\cap Z$. By re-embedding $\pp \na$ by the $d$-th Veronese, we may
linearize
the equation of $Z$, so the same conclusion holds for the map 
$
H^i(X,\Bbb Z)\ra H^i(X\cap Z, \Bbb Z)
$.
Call this extension version 2.

If we take $X=\pp\na$ and cut by hypersurfaces, as long as we have a
smooth variety, i.e. a smooth complete intersection at each step,
we can continue iteratively to get:

\proclaim{3.3 Lefschetz theorem version 3  }
Let $X^n\subset\pp\na$ be
a smooth complete intersection.  Then the restriction map on cohomology:
$$
H^i(\pp \na ,\Bbb Z)\ra H^i(X, \Bbb Z)
$$
is an isomorphism for  $i\neq n$  and injective for $i=n$.
\endproclaim

Thus for smooth complete intersections, almost all the cohomology
 is inherited from the ambient
projective space. The Lefschetz theorem  extends
to fundamental groups,
see, e.g. [Ful], so   a good deal of the topology
of a complete intersection
is inherited from the ambient projective space.

The hypotheses on the Lefschetz theorem can be relaxed   in several
different
ways. Here is a generalization due to Goresky and MacPherson:

\proclaim{3.4 Theorem [GM]} Let $X^n$ be an algebraic variety and 
suppose $\pi: X\ra\pp\na$ is a (not necessarily proper)
algebraic map with finite fibers. Let $L\subset\pp\na$ be a linear subspace
of codimension $c$.
Let $L_{\epsilon}$ denote an $\epsilon$-neighborhood of $L$ with respect to
some Riemannian metric on $\pp\na$.

If $\epsilon>0$ is sufficiently small, then the inclusion
$\pi\inv (L_{\epsilon})\ra X$ induces an isomorphism
on intersection homology groups
$IH_i(\pi\inv (L_{\epsilon})\Bbb Z)\simeq IH_i(X,\Bbb Z)$ for all
$i<n-c$ and a surjection 
$IH_{n-c}(\pi\inv (L_{\epsilon})\Bbb Z)\ra IH_{n-c}(X,\Bbb Z)$.

Furthermore, if $L$ is generic, then $L_{\epsilon}$ may be replaced
by $L$ in the above formulae.
\endproclaim

Using the generalized version
of the Lefschetz theorem (actually   other results
that are essentially equivalent, e.g.
Deligne's extension of Bertini's theorem) Fulton and Hansen
proved the following theorem which has striking consequences:

\proclaim{3.5 Connectedness Theorem [FH]} Let $Z^n$ be a projective
variety, let $f: Z\ra\pp m\times\pp m$ be a finite morphism, and
let $\Delta\subset\pp m\times\pp m$ denote the diagonal. If
$n>m$, then $f\inv (\Delta )$ is connected.
\endproclaim

As explained
by Fulton, see [Ful],
 Bertini, Bezout and Lefschetz's theorems
can all be understood in terms of  
connectedness.

An important application of the connectedness theorem is the following,
which reduces the studies of degenerate secant and tangential
varieties to the same problem. 
 
If $Y\subseteq X\subseteq\pp\na$, define $\tau (Y,X)=\cup_{y\in
Y}T^{\star}_yX$.

\proclaim{3.6 Theorem [Z], [FH]} Let
$X^n,Y^y\subset\ppp  V$ be varieties, respectively  of dimensions
$n,y$. Assume $Y\subseteq X$.
Then either

$\tdim \sigma (Y,X)= n+y+1$ and $\tdim\tau (Y,X)= n+y$

or

$\sigma(Y,X)=\tau (Y,X)$.\endproclaim

Taking $Y=X$, one sees that if either $\sx$ or $\tx$ is
not of the expected dimension, then they must be equal, which 
is the version proved by Fulton and Hansen. The extension is
due to Zak.

\demo{Proof}
Assume $\tdim\tau (Y,X)=t< n+y$. We need to show 
$\sigma(Y,X)=\tau (Y,X)$.
Project $\pp\na$ 
from some linear space $L^{n+a-t-1}$ avoiding 
$\tau (Y,X)$ to a $\pp t$ and
consider the map $F:  X\times Y\ra \pp t\times \pp t$. Since
$F$ restriced to each factor is finite,
$\tdim F(X\times  Y)=n+y>t$, so $F\inv (\Delta )$ is connected. Assume
 $\sigma(Y,X)\neq\tau (Y,X)$, so
$\tdim \sigma(Y,X)>\tau (X,Y)$ 
and thus $\sigma (X,Y)$ 
 intersects $L$. 
Thus there exists $x\in X\backslash Y$, $y\in Y$, such that
the line $\overline{xy}\cap L\neq \emptyset$. Thus $F(x,y)=F(y,y)$.
By the connectedness theorem, there exists an arc 
 $(x_t,y_t)$ in $F\inv (\Delta)$ such that 
$(x_0,y_0)=(x,y)$ and $(x_1,y_1)=(y,y)$. By continuity there exists
some minimal $t'$ such that $x_{t'}=y_{t'}$ and thus the limiting
line
as $t\ra t'$ is a tangent line that intersects $L$, giving a contradiction.
\qed\enddemo
\smallpagebreak

 The following important application of the connectedness theorem is due
to Zak: 

\proclaim{3.7 Zak's theorem on tangencies,
  [Z]} Let $b=\tdim X_{sing}$
(Set
$b=-1$ if
$X$ is smooth.) Let $L\in \gamma_k(X)$ be any point, then
$\tdim \{
x\in X\mid \tilde T_xX\subseteq L\}
\leq k+ (b+1)$.
 \endproclaim

\demo{proof} Assume $X$ is smooth. Let $Y=Y_L$ and let $y=\tdim Y$.
We have $\tau (Y,X)\subset L$, but since $X$ is not contained in
a hyperplane $\sigma (Y,X)\not\subseteq L$. Thus by (3.6), $\tdim\tau(Y,X)=
y+n$. Thus $y+n\leq n+k$, i.e. $y\leq k$.\qed\enddemo

\subheading{3.7.1 Exercise} Prove the general case. Hint: consider
the intersection of $X$ with a linear space.

\smallpagebreak

In comparison to   Lefschetz theorem version 3, Barth and Larsen
proved the following theorem about an arbitrary smooth subvariety
of projective space:

\proclaim{3.8 Theorem, see [B], [BL], [Hart]} Let
$X^n\subset\pp\na$ be a smooth variety, then the restriction map on
cohomology
$$
H^i(\pp \na ,\Bbb Z)\ra H^i(X, \Bbb Z)
$$
is an isomorphism for $i\leq n-a$.
\endproclaim

Thus, smooth varieties of small codimension cohomologically \lq\lq look
like\rq\rq\ complete intersections.

A basic question in geometry is: To what extent is the geometry of 
a manifold   determined by its topology? For example, in Riemannian
geometry,   certain averages
of the Riemann curvature tensor are invariants of the differentiable
topology.

R. Hartshorne posed several questions in algebraic geometry to the
effect of asking when something is true on the level of cohomology 
implies the corresponding result on the level of geometry. For example, 
there is the classical question: under
what circumstances are cohomology classes represented by an algebraic
cycle (or variety)? Among his questions is the following
famous conjecture:

\proclaim{3.9 Hartshorne's conjecture on complete intersections}([Hr],
1974)
Let $X^n\subset\bccc\bpp\na$ be a smooth variety. If 
$a<\frac n2$ then $X$ is a complete intersection.
\endproclaim

If $X$ is a complete intersection, then all polynomials of degree
$d$ on $X$ are the restriction of polynomials of degree
$d$ on $\ppp V$. In other words, the map
$H^0(\pp\na, \Cal O_{\pp\na}(d))\ra H^0(X,\Cal O_X (d))$
is a surjection.  If  all up to 
 the $d$-th maps are surjective, one
says that
$X$ is $d$-{\it normal}, and if all maps are surjective,
one says that $X$ is {\it projectively normal}. $1$-normal
varieties are called {\it linearly normal}.  Another way of saying
that $X$ is linearly normal is that $X$ cannot be
realized as the linear projection of some $\tilde X\subset\pp{\na+1}$.
More generally, $X$ is $d$-normal if $v_d(X)$ is not the linear projection
of some $Y\subset\pp{\binom{\na +d}d}$. It
has been shown that under certain
circumstances, if $X$ is projectively normal, it must be
a complete intersection **ref?**. Thus we may think of $X$ being $d$
normal   as an approximation to being a complete intersection.

 Thus a first approximation to Hartshorne's conjecture
on complete intersections
(also conjectured by Hartshorne) is that if $X^n\subset\pp\na$ is
a smooth variety, it must be linearly normal if 
$a<\frac n2$. Now say
$Y\subset\pp { N+1}$ is a smooth variety. We can project $Y$ 
to a smooth subvariety of some $\pp N$ iff there exists
a $p\in \pp{N+1}$ such that $p$ does not lie on any secant
or tangent line of $Y$, i.e. iff $\sigma (Y)\neq\pp{N+1}$.
Zak proved that if $N+1-n\leq\frac n2$ then $\sigma (Y)=\pp{N+1}$
and this is his theorem on linear normality quoted in the introduction.

\heading\S  4.  Projective differential invariants\endheading

\subheading{4.1 The moving frame}

As discussed in the introduction, we will study the local geometry
of a subvariety $X^n\subset\ppp V$ by studying the infinitesimal motions
of flags in $V$, in fact of bases of $V$ adapted to
the geometry of $X$.

Let $N=n+a$.
To begin, let  $\cf$ denote the space of bases of $V$. $\cf$  is
isomorphic to $GL(V)$. We write $f\in \cf$ as
$f= (\ee 0\hd\ee N )$,  where we think of the basis
vectors $\ee B\in V$ as column vectors.

We will use the projection 
$$
\align \cf &\ra\ppp V\tag 4.1.1 \\ 
 f= (\ee 0\hd\ee N )&\mapsto [\ee 0].
\endalign
$$

\subheading{4.1.2 Exercise} Define a projection $\cf\ra \BG (k, \ppp V)$.
What is the group preserving the fiber?

\smallpagebreak

We will often consider $f\in\cf$ as a matrix, e.g. $f= (g^A_B)$,
writing
$$
\ee B = \pmatrix g^0_B\\ \vdots\\ g^{N}_B\endpmatrix .
$$

One any Lie group $G$, there is a canonical left-invariant 
$\fg$-valued one-form, the {\it Maurer-Cartan form}
$\Omega\in \Omega^1(G,\fg )$.

In the case of a matrix Lie
group (our situation),
we think of the matrix elements $g^A_B$ as coordinate functions,
and take their derivative to get a mapping
$dg_f:T_f\cf\ra M_{N+1\times N+1}(\Bbb C)$.
($M_{k\times l}(\Bbb C)$ denotes the vector space of
$k\times l$ matrices.) In this situation, the Maurer-Cartan
form is
$$
\Omega := g\inv dg\in \Omega^1(\cf,
M_{N+1\times N+1}(\Bbb C)).
$$

 Utilizing the fact that the differential
of a constant map  is zero, we calculate
$0=d(gg\inv)=gdg\inv + (dg)g\inv$ and find that $d\Omega = d(g\inv dg)=
d(g\inv )\ww dg = -g\inv dg g\inv\ww dg$. I.e.
$$
d\Omega =
-\Omega\ww\Omega ,\tag 4.1.3 
$$
 which is called the {\it Maurer-Cartan equation}. In 
indicies, 
$$
d\ooo AB = -\ooo AC\ww\ooo CB\tag 4.1.4
$$
(here and throughout, we use the convention that repeated indicies
occuring up and down are to be summed over).

$\ooo AB$ has the geometric interpretation of measuring the
infinitesimal motion of $\ee A$ towards $\ee B$,
as $d\ee A = \ooo 0A\ee 0 +\hdots+ \ooo N A\ee N$.
In particular, the infinitesimal motion of $\ee 0$ is measured 
by the forms $\ooo B0$, so it is not suprising that  
$\{ \oo 1\hd \oo N\}= \pi^* (T^*_{[\ee 0]}\ppp V)$.

We now adapt to the geometry of our situation, the flag
$\hat x\subset \hat T_xX\subset V$. We let
$\cf^1_X=\cf^1\subset\cf$ be the subbundle   respecting
this flag. 
 
 Write  $f=(\ee 0 ,\ee\alpha ,\ee\mu )$ for an element
of $\cf^1$, where $1\leq\alpha , \beta\leq n,\  n+1\leq\mu , \nu\leq n+a$
where $[\ee 0 ]=x$, and  $\{\ee 0\hd\ee n\}=\hat T$.
$\cf^1$ is not a Lie group, but it is a $G_1$-principal bundle where
$$
G_1=\left\{
g\in GL(V) |\,  g = \pmatrix g^0_0 & g^0_{\beta} & g^0_{\nu} \\
 0 &  g^{\alpha}_{\beta} & g^{\alpha}_{\nu} \\
0 &  0 & g^{\mu}_{\nu} \endpmatrix \right\} \tag 4.1.5
$$

Let $i: \cf^1\subset \cf$ denote the inclusion.

 $i^*(\oo\mu )=0$ because
$d\ee 0\equiv \oo\a\ee\a\tmod \ee 0$ 
and for the same reason,
$i^*(\oo 1\ww\hdots\ww \oo n)$ is nonvanishing. 

   From now on, I
commit a standard abuse of notation,  omitting the $i^*$ in
the notation,   the pullback being  understood from the context.

Expanding out $0=d\oo\mu=-\ooo\mu\beta\ww\oo\beta$ using the Maurer-Cartan
equation   (4.1.4)
 yields
 $$
\ooo\mu\beta = \qq\mu\alpha\beta\oo\alpha\tag 4.1.6
$$
for some functions $\qq\mu\alpha\beta = \qq\mu\beta\alpha$
defined on $\cf^1$. 

 To prove the
symmetry, one uses the {\it Cartan lemma}: Let $W$ be a vector space and
let $w_1\hd w_k$ be independent elements. If $v_1\hd v_k\in W$
are such that $w_1\ww v_1+\hdots w_k\ww v_k=0$,
then $v_j=\Sigma_ih_{ij}w_i$ with $h_{ij}=h_{ji}$. 

\subheading{4.1.7 Exercise} prove 
the Cartan Lemma.

\smallpagebreak

 Consider 
$$
\tilde{II}=\qq\mu\alpha\beta\oo\alpha\oo\beta\ot (\aa\mu\tmod\hat T)
\tag 4.1.8
$$
  
\subheading{4.1.9 Exercises} 

1. Show that  $\tilde{II}\ot \ee 0^*$
is invariant under motions in the fiber over $[\ee 0]$ and thus descends to
be
a well defined section of $S^2T^*X\ot NX$, which,  as
expected, is $II$.

2. Show that this definition is equivalent to
the  definitions of $II$ given in the introduction. Hint: use
the projection $\cf^1\ra G(k+1,V)$ to show equivalence
with the Gauss map definition and use a section of 
$\cf^1$ given in coordinates to prove the second equivalence,
e.g., take $\ee\a = \frp{}{\xx\a} - f^{\mu}_{\xx\a}\frp{}{\xx\mu}$,
$\ee\mu=\frp{}{\xx\mu}$.

\bigpagebreak

\subheading{4.2 Higher fundamental forms}

\smallpagebreak

We can define the higher fundamental forms $\ff k$ as maps
$\ff k : \tker \ff{k-1}\ra S^kT^*$, where we begin
by considering $II=\ff 2 : N^*\ra S^2T^*$. Note that in projective
geometry this mapping is more natural
than the dual map
$S^2T\ra N$ as   $\ppp N^*_xX$ has the geometric
interpretation of the  space of hyperplanes tangent to $X$ at $x$.  

In coordinates, the $k$-th  fundamental form 
corresponds to the first nonzero
terms in the Taylor series expansion using  
coordinates adapted to
the filtration of $N^*$. 

Fix a general point $x\in X$ and let $a_1=\tdim II(T,T)$ and
let 
$$\hat T\up 2_xX= \hat T+ II(T,T)(-1)\subseteq V
$$ 
denote
the {\it second osculating space} of $X$ at $x$.

One way to define $III$ is as the derivative of
the second Gauss map
$$
\align \gamma\up 2 : X &\dashrightarrow G(n+a_1, V)\tag 4.2.1 \\
x&\mapsto \hat T\up 2  \endalign
$$

\subheading{4.2.2 Exercise} Write out this definition.

\smallpagebreak

To define $\ff 3$ in frames, consider the quantity
$$
\tilde{\ff 3}:=
 \ooo\nu\mu\ooo\mu\beta\oo\beta \ot e_{\nu}\text{mod }\hat T\up 2 \tag
4.2.3
$$
$\tilde{\ff 3}$
descends (after twisting) to  a well
defined  element $\ff 3\in\Gamma ( S^3T^*\ot (N/ II(S^2T)))$.
$\ff 3$ has the geometric interpretation of measuring how $X$ is 
moving away from its second osculating space to first order.
We sometimes use $III$ to denote $\ff 3$.

Adapt frames to the filtration of $N$, letting
$\{\ee\xi\} =II(S^2T)$ and $\{\ee\phi\}= N\tmod II(S^2T)$,
with index ranges $n+1\leq \xi,\eta\leq n+a_1$,
$n+a_1+1\leq\phi,\psi\leq\na$. Let $\cf^2$ denote
the resulting bundle.
Then
$$
\tilde{\ff 3}:=
 \ooo\phi\xi\ooo\xi\beta\oo\beta \ot e_{\phi}\text{mod }\hat T\up 2. \tag
4.2.4
$$

\subheading{4.2.5 Exercise} Using this definition, show that
$\ff 3$ is a symmetric cubic form. 

\smallpagebreak 

There is a further restriction on $\ff 3$ at general points.
(4.2.3) implies so
  $|\ff 3|\subseteq (\ii\ot T^*)$.

\subheading{4.2.6 Definition} Let $T $ be a vector 
space  and let  
$A\subset S^dT^*$ be a linear subspace. 
Define the $k$-th {\it prolongation} of
$A$ by
$A\up k:= S^{d+k}T^*\cap (A\ot S^kT^*)$. Geometrically, $A\up k$
is the space of 
  polynomials of degree $d+k$ on $T$  having the property that all
their $k$-th derivatives lie in $A$.

\proclaim{4.2.7 Theorem (Cartan) [C vIII.1  p 377]}Let $X^n\subset\pp\na$
be a
variety and let $x\in X$ be a general point. Then
$$
|\ff k|_{X,x}\subseteq \ii_{X,x}\up{k-1}. \tag 4.2.8
$$
\endproclaim

We will call (4.2.8) the {\it prolongation property}. It can fail
to hold at special points where $II$ degenerates.

The two remarks above prove Cartan's theorem for $\ff 3$,  
the general case is left  as an exercise.

\bigpagebreak

\subheading{4.3 Relative differential invariants}

\smallpagebreak

You may wish to skip this section on a first reading. Life would
be easier if one could just skip it altogether, but unfortunately
the cubic form plays an important role any time the geometry  of $X$ is
not completely determined by $\ii$.

\smallpagebreak

In addition to  $III$, there is
another third order invariant  called the {\it cubic form}, which we will
denote
$F_3$.   The geometric interpretation of $F_3$ is that it measures  how the
second fundamental form varies infinitesimally, or in other words, how $X$
is
leaving its embedded tangent space to second order.

Namely, differentiating $0= d(\ooo\mu\alpha - \qq\mu\alpha\beta\oo\beta)$  
implies there exist functions $\rr\mu\alpha\beta\gamma$ defined 
on $\cf^1$ that satisfy
$$
 \rr\mu\alpha\beta\gamma\oo\gamma  =
-d\qq\mu\alpha\beta - \qq\mu\alpha\beta\ooo 00 -\qq\nu\alpha\beta\ooo\mu\nu
+\qq\mu\alpha\delta\ooo\delta\beta + \qq\mu\beta\delta\ooo\delta\alpha .
\tag 4.3.1
$$
For $f\in\cf^1$, let $F_3=(F_{3})_f\in \pi^*(S^3T^*\ot N)$ be
$$
F_3= \rr\mu\alpha\beta\gamma\oo\alpha\oo\beta\oo\gamma\ot\ee\mu .\tag 4.3.2
$$ 
Note that $F_3$ does {\it not} descend to be a well defined section
of $S^3T^*\ot N$ (even after twisting). In fact,
if $(\tilde \aa 0,\tilde\aa\alpha,\tilde\aa\mu )$ is a new frame with
$$
\align
&\tilde\aa\mu = \aa\mu + g^0_{\mu}\aa 0 + g^{\alpha}_{\mu}\aa\alpha
\tag 4.3.3\\
& \tilde\aa\alpha = \aa\alpha + g^0_{\alpha}\aa 0\endalign
$$
then
$$
\tilde\rr\mu\alpha\beta\gamma = \rr\mu\alpha\beta\gamma +
\frak S_{\alpha\beta\gamma}g^0_{\alpha}\qq\mu\beta\gamma +
\frak
S_{\alpha\beta\gamma}g^{\delta}_{\nu}\qq\nu\alpha\beta\qq\mu\gamma\delta
\tag 4.3.4
$$

Motions by $g^{\mu}_{\nu},g^{\alpha}_{\beta}, g^0_0$
also vary the $\rr\mu\alpha\beta\gamma$ but the change they effect
is cancelled by the corresponding changes in the $\oo\alpha$ and $\ee\mu$.
$F_3$ is an example of what is called a {\it relative invariant}.
We will use the notation $\Delta\rr\mu\alpha\beta\gamma$ to denote
the change in $\rr\mu\alpha\beta\gamma$  by a  fiber motion of
the type in (4.3.3).
By (4.3.4),
$$
\Delta\rr\mu\alpha\beta\gamma =
\frak S_{\alpha\beta\gamma} (\gg 0\alpha\qq\mu\beta\gamma +
\gg\delta\nu\qq\nu\alpha\beta\qq\mu\gamma\delta ). \tag 4.3.5
$$

\smallpagebreak 

It is possible to define $F_3$ as a section of a bundle
well defined over $X$, namely
$$
F_3\in\Gamma (X,\frac{ S^3T^*\ot N}{T^*\circ II + <II, T\ot N^*, II>})
\tag 4.3.6
$$
where $<\cdot, \cdot,\cdot>$ is the natural contraction.

\smallpagebreak 

One can compare $F_3$ to the covariant derivative of the second
fundamental form of a submanifold
of a Riemannian manifold, $\nabla^{riem} II $, which is a well defined
tensor. In fact one may think of $F_3$ as an the equivalence class of
$\nabla^{riem} II$'s, where they range over the
(holomorphic)  metrics compatible with the
projective  structure. 

Differentiating
$F_3$ one obtains a fourth order invariant $F_4\in \pi^*(S^4T^*\ot N)$
whose
coefficients
$\rr\mu\alpha\beta{\gamma\delta}$  are defined by
$$
\rr\mu\alpha\beta{\gamma\delta}\oo\delta =
-d\rr\mu\alpha\beta\gamma -\rr\mu\alpha\beta\gamma\ooo 00
-\rr\nu\alpha\beta\gamma\ooo\mu\nu +
\frak S_{\alpha\beta\gamma}(
\rr\mu\alpha\beta\epsilon\ooo\epsilon\gamma -\qq\mu\alpha\beta\ooo 0\gamma
+\qq\mu\alpha\epsilon\qq\nu\beta\gamma\ooo\epsilon\nu ).\tag 4.3.7
$$
(If $III=0$, $F_4$ is the only fourth order invariant.)
The geometric interpretation of $F_4$ is that it measures how $X$ leaves
its
embedded  tangent
space to third order.

 Under a change of frame (4.3.4) the
coeffecients of $F_4$ vary as follows:
$$
\tilde\rr\mu\alpha\beta{\gamma\delta} = 
\rr\mu\alpha\beta{\gamma\delta} +
\frak S_{\alpha\beta\gamma\delta}g^0_{\alpha}\rr\mu\beta\gamma\delta +
\frak S_{\alpha\beta\gamma\delta}g^{\epsilon}_{\nu}
(\rr\nu\alpha\beta\gamma\qq\mu\delta\epsilon+
\qq\nu\alpha\beta\rr\mu\gamma\delta\epsilon) 
+ g^0_{\nu}\qq\mu\alpha\beta\qq\nu\gamma\delta .
\tag 4.3.8
$$
  One can
continue, defining forms $F_k$ for all $k$.

In coordinates adapted to a point $x$, 
the  coefficients of the $F_k$ at $x$ 
are the $k$-th derivatives of the embedding.

\subheading{4.4  Yet
another definition of fundamental forms}

Here is yet another definition
of the fundamental forms that will be particularly
useful for calculating fundamental forms of re-embeddings and homogeneous
spaces. If you are more at home with spectral sequences than Gauss maps,
coordinates, or moving frames, then this definition is for you.

\smallpagebreak

Define inductively a series of maps (following [L6]):
$$
\ud k\aa 0 : (T\cf^1)^{\ot k}\ra V/\tim
(\ud 0\hd \ud{k-1}) \tag 4.4.1
$$
as follows:
Let $d$ denote exterior differentiation, let $\ud 0\aa 0 = \aa 0$
and let $\ud 1\aa 0 = d\aa 0 \tmod \aa 0$.
 If $v_1\hd v_k\in T_f\cf^1$, extend
$v_1\hd v_k$ to holomorphic vector fields in some neighborhood of $f$
which we   denote $\tilde v_1\hd \tilde v_k$.
Let
$$
 \ud k\aa 0 (v_1\hd v_k ) :=  v_1(\intprod d (
\tilde v_2 \intprod
\hdots d(\tilde v_k\intprod d\aa 0 )
 \text{ mod }\pi_k\inv (\text{Image}\ud{k-1})
\tag 4.4.2
$$
where
$\pi_k: V\ra V/(\text{Image}\{\ud 0 \hd\ud{k-1}\} )$ is the projection,
and $\intprod$ denotes the contraction
$T\times  T^{*\ot l}\ra  T^{*\ot l-1}$. (4.4.2) is independent of
the extension of $v_1\hd v_k$ to holomorphic vector fields.
(The proof that (4.4.2) is independent of
 the choice of extension to holomorphic
vector fields is the same as the standard argument in the real
case, see e.g. [S].)
   $\ud k\aa 0\ot (e_0^*)^{k-1} =\pi^*(\ff k)$.
The maps $\underline d^k$ can be defined more algebraically
as follows;

The quotient  map
$$
  V^* \ra  V^*/ \hat x\upperp=\Cal O_{\bppp V }(1)_x \tag 4.4.3
$$
gives rise to a spectral sequence of a filtered complex by letting
$$
\align F^0K^0 &= V^*\ \ \  \ \ \  F^0K^1 = \Cal O_{X }(1)_x \tag 4.4.4\\
F^1K^0 &= 0 \ \ \  \ \ \ \ \ F^p = F^pK^1 =   \frak m^p_x(1).
\endalign
$$
The maps are
$$
\align
& \ud 0 : V^*\ra F^0/F^1 = \Cal O_{X,x}(1)/\frak m_x(1)\simeq\bcc{}\tag
4.4.5
\\
& \ud 1: \text{ker}\ud 0  \ra F^1/F^2 = \frak m_x (1)/\frak m^2_x (1)
\simeq T^*(1) \\
& \ud 2: \text{ker}\ud 1 \ra F^2/F^3 = \frak m^2_x (1)/\frak m^3_x (1)
\simeq (S^2T^*)(1) \\
& \vdots \endalign
$$

For example,
 the first two terms of   (4.4.5) expressed in frames are
 $$
\align
\ud 1\aa 0 &= \oo\alpha\ot\aa\alpha \text{ mod } \hat x\tag 4.4.6\\
\ud 2\aa 0 &= \oo\alpha\ooo\mu\alpha\ot\aa\mu\text{ mod}\hat T
 \tag 4.4.7
\endalign
$$
  
We will use the notation $F^1=\ff 1= Id_T$, $F^0=\ff 0= \hat x\ot \hat
x^*$.

\smallpagebreak

This definition  is particulary useful in computing fundamental forms
of homogeneous spaces. Say $G/P\subset\ppp V$ and
$V$ is a vector space formed from a vector space
$W$ (e.g. an exterior power) and originally $G\subseteq GL(W)\subset
GL(V)$.
Then we can use the smaller $G$-frame bundle
to do our computations.

\subheading{4.5 Fundamental forms of Veroneses $v_p(\ppp V)\subset\ppp
S^pV$}

Let $V=\bcc{n+1}$ have basis $\{\ee 0,\ee\alpha\}$,
$1\leq\alpha\leq n$ and let $x=[(\ee 0)^p]$. Using
the bundle of $\rho_p (GL(V))\subset GL(S^pV)$ frames, we have
$$
\align
&\underline de_0^d\equiv p\oo\alpha\ee\alpha e_0^{p-1}\tag 4.5.1\\
&\underline d^2e_0^d\equiv p(p-1)\oo\alpha\oo\beta\ee\alpha
\ee\beta e_0^{p-2}\\
&\underline d^ke_0^d\equiv p(p-1)\hdots (p-k+1)\oo{\alpha_1}
\hdots \oo{\alpha_k}\ee{\alpha_1}\hdots
\ee {\alpha_k} e_0^{p-k}\endalign
$$
Thus
$$
\align
&|\ff k_{ v_p(\ppp V)}|= \ppp S^kT^*\ k\leq p\tag 4.5.2\\
&|\ff k_{ v_p(\ppp V)}|= \emptyset\ k>p. \endalign
$$

\subheading{4.6 Fundamental forms of Veronese re-embeddings
$v_d(X)\subset\ppp
S^dV$}

\subheading{4.6.1 Fundamental forms of  $v_2(X)\subset\ppp S^2V$}
Assume $X\subset\ppp V$ is such that $II_X$ is surjective
so the only differential invariants are $F_k=(F_{k})_{ X,x}$. Write
$x=[\ee 0]$, $v_2(x)=[\ee 0\circ\ee 0]$ and
use the Leibnietz rule applied to 
$\ee 0\circ\ee 0$ to compute (see [L6, 3.2] for details):
$$
\align
\Bbb F\Bbb F^1_{ \vtwox ,x } &= 2F_1F_0|_{\hat x^2\upperp} \tag 4.6.2\\
\Bbb F\Bbb F^2_{ \vtwox ,x} &= 2(F_2F_0 + F_1F_1)|_{ (\hat x\hat T)\upperp}
 \\
\Bbb F\Bbb F^3_{\vtwox ,x} &= 2(F_3F_0 + 3F_2F_1) |_{\tker\bii 2\vtwox} \\
\Bbb F\Bbb F^4_{\vtwox ,x}&= 2(F_4F_0 + 4F_3F_1 + 3F_2F_2)|_{\tker\bii
3\vtwox}\\
\Bbb F\Bbb F^5_{\vtwox ,x}&= 2(F_5F_0 + 5F_4F_1 + 10F_3F_2)|_{\tker\bii
4\vtwox}\\ &\vdots
\endalign
$$
 
\proclaim{4.6.3 Proposition [L6, 3.10]}
The fundamental forms of $v_d(X)$ are
$$
\bii k \vdx =
\Sigma_{l_1+\hdots +l_d=k}c_{l_1\hdots l_d}F_{l_1}\hdots
F_{l_d} \tmod (\Sigma_{l<k} \bii l\vdx (S^lT) )|_{\tker\bii {k-1}\vdx} \tag
4.6.4 
$$
where the $c_{l_1\hdots l_d}$ are nonzero constants.
\endproclaim

For example
$$
\bii 4{v_3(X)}= c_{400}F_4F_0F_0+ c_{310}F_{3}F_1F_0 +
c_{220}F_2F_2F_0+c_{211}F_2F_1F_1. \tag 4.6.5
$$
For the proof, see [L6].

Note that if $Z$ is a hypersurface of degree $d$, then
$Z$ osculates to order $p$ at $x\in X$ if and only if
$H_Z\in \tker\fff k{\vdx, x^d}$.

\subheading{4.6.4 Proofs of  (0.5.8), (0.5.9)}
 
\demo{Proof of [L6, 3.16]} Observe that $| \ff k_{\vdx}|= \ppp S^kT^*$ for
all
$k\leq d$ and sum up.\qed\enddemo

\demo{Proof of [L6, 3.17]}
The first term in any $\bii k\vdx$ for which $S^dN^*$ is not in the kernel
  is
$(F_2)^d$, which appears in $\bii {2d}\vdx$. Thus
$S^dN^{*}\not\subseteq\tker\bii{2d-1}\vdx$. Finally, note  that
 dim$S^dN^{*}=\binom {a+d-1}d$.\qed
\enddemo

\subheading{4.7 Fundamental forms of Grassmanians}

Let $W=\bcc n$ and let $G(k,W)\subset\ppp (\Lambda^k W)$.
 Write $V=\Lambda^kW$. Use index ranges $1\leq i,j\leq k$, $k+1\leq s,t\leq
n$. Write the Maurer-Cartan form of $GL(W)$ as
$$
\Omega = \pmatrix \ooo ij &\ooo it\\ \ooo sj&\ooo st\endpmatrix
\tag 4.7.1 
$$
Using the embedding $\rho : GL(W)\ra GL(V)$, let $E=\ee 1\ww\hdots\ww \ee k
\in G(k,W)$, 
$E_s^j=\ee 1\ww\hdots\ww\ee{j-1}\ww \ee s\ww\ee{j+1}\ww\hdots\ee k$
where $\ee j$ has been replaced by $\ee s$, let $E_{st}^{ij}$ be $E$ with
$\ee s$ replacing $\ee i$ and $\ee t$ replacing $\ee j$
and so on.
Then the $E_{s_1\hdots s_p}^{j_1\hdots j_p}$,
$1\leq p\leq k$, give a basis of $V$, 
in fact adapted to the representation $\rho$. Now we take derivatives at
$E$
using (4.4) and the Leibinitz rule:
$$
\align
&\underline dE\equiv \ooo siE_s^i \tmod E\\ &\underline d^2E\equiv
 \Sigma_{i<j,s<t}(\ooo si\ooo tj -\ooo ti\ooo sj)E_{st}^{ij}\tag 4.7.2\\
&\underline d^3E\equiv \Sigma_{i<j<l, s<t<u}
det\pmatrix\ooo si &\ooo ti &\ooo u i\\
\ooo sj&\ooo tj&\ooo u j\\
\ooo sl&\ooo tl&\ooo u l\endpmatrix
E_{stu}^{ijl}
\endalign
$$
 The first line can be used to recover that $T_EG(k,W)= E^*\ot W/E$.
Continuing, we see that for any $x\in G(k,W)$,
$$
|\ff p | = \{ p\times p\text{ minors of }(\ooo s j )\}=I_p\sigma_{p-1}(
Seg(\ppp E^*\times\ppp (V/E)).\tag 4.7.3
$$
In particular, the last nonzero fundamental form is the $\tmin(k,n-k)$-th.

\subheading{4.8 Fundamental forms of spinor varieties}

Choose a basis of $V=\bcc{2m}$ so that $Q=\pmatrix 0&I\\ I&0\endpmatrix$.
With respect to this basis $SO(V,Q)$ has Maurer-Cartan form
$$
\pmatrix 
\ooo ij & \ooo i{n+k}\\
\ooo{n+l}j &\ooo{n+k}{n+l}\endpmatrix\tag 4.8.1 
$$
where $\ooo{n+k}{n+l}= \ooo{ k}{ l}, \ooo i{n+k}=-
\ooo k{n+i}, \ooo{n+l}j=- \ooo{n+j}l$. 

Let $E=\ee 1\ww\hdots\ww\ee m\in \Bbb S_m= \Bbb S_m(V,Q)$.
  Computing as for the Grassmanians, 
$$
\align
dE &\equiv
  \ooo {n+i}j\EE j{n+i}  \tmod E\tag 4.8.2\\
d^2E 
&\equiv 
\ooo {n+i}j\ooo{n+k}l\EE {jl}{n+i,n+k}  \tmod \tim dE\\
&\vdots \endalign
$$
Note that $Q$ allows us to identify $V/E$ with $E^*$. The additional
structure on the tangent space induced from that of the Grassmanian
becomes
$$
T=\La 2 E^*.
$$
Moreover,
$$
 |\ff k|\simeq \La{2k} E^* = I_k\sigma_{k-1}G(2,E)
$$
 in particular
$|II|= I_2G(2,E)$.

Note that one can construct the  space of positive spinors
$\La{even} E^*$ from the fundamental forms, without any knowledge
of Clifford algebras.

\medpagebreak

\subheading{4.9 Hermitian symmetric spaces}

\subheading{4.9.1 Remark}
Note that in the examples of Grassmanians 
and Spinor varieties above, $|\ff k| = |\ff 2|\up{k-1}$
where for a general variety the prolongation property only implies
$|\ff k |\subseteq |\ff 2|\up{k-1}$.  Among homogeneous varieties,
there is a preferred subclass, the {\it Hermitian symmetric spaces}.

\subheading{4.9.2 Definition} Let $X=G/P\subset \ppp V$ be a homogeneous
variety. Let $H$ be the
semi-simple part of $P$.
 We say $X$ is a  {\it Hermitian symmetric space} if $T_xX$
is irreducible as an $H$-module. (This is a nonstandard definition,
normally
the definition is made using a Riemannian metric.)

\proclaim{4.9.3 Theorem [LM1]}
Let $X = G/P \subset \BP   V $ be a 
Hermitian
symmetric space in its fundamental embedding and let $x\in X$. 
Then for $k\ge 2$,
$$
|\ff {k+1}_{X,x}|= |\ff 2_{X,x}|\up {k-1}
$$
\endproclaim

For the proof, see [LM1]. This strict prolongation property fails
for a general homogeneous space, e.g.,  it does not hold for
$G_{Q-null}(k,m)$ for $2\leq k < [\frac m2]$. This property
has the following geometric  consequence (see [LM1] for the proof):

\proclaim{4.9.4 Corollary [LM1]}
Let $X$ be a Hermitian symmetric space in its fundamental embedding, and
let $x\in X$. Then
$$\tbaseloc |\ff k_{X,x}| = \sigma_{k-1}(\tbaseloc |\ff 2_{X,x}|).
$$ 
\endproclaim

\subheading{4.10 Fundamental forms of Segres}

Let $U=\bcc{a+1},W=\bcc{b+1}$ and consider the Segre,
$Seg(\ppp U\times \ppp W)\subset\ppp (U\ot W)$. Here we use
$G=GL (U)\times GL(W)\subset GL (U\ot W)$ frames. Let
$(\ee 0\hd\ee a)$, $(f_0\hd f_b)$, $\oo\a,\ett j0$, $1\leq\a\leq a$,
$1\leq j\leq b$ denotes bases and dual bases. We compute
$$
\align 
\underline d (\ee 0\ot f_0 )&\equiv \oo\a \ee \a\ot f_0 + \ett j0\ee 0
\ot f_j \tmod \ee 0\ot f_0\\
\underline d^2 (\ee 0\ot f_0) &\equiv \oo\a \ett j 0 \ee \a\ot f_j\endalign
$$
We see $II$   surjects onto $N$ and 
$$
|\ff 2|= I_2(\ppp (\ee 0\ot (W/\{f_0\}))\sqcup \ppp ((W/\{\ee 0\}\ot f_0))
\tag 4.10.1
$$

More generally,
consider    $X= Seg (\pp {r_1}\times
\hdots
\times\pp{r_d}) =
Seg(\ppp W_1\times\hdots\times \ppp W_d)$,
$\tdim W_i= r_i$. We use the
$GL(W_1)\times\hdots\times GL(W_d)\subset
GL(W_1\ot\hdots\ot W_d)$ frame bundle. 

Let $(\ee{0j}\hd \ee{r_jj})$ denote an
  element of the frame bundle for $\ppp W_j$ and write the 
$j$-th block diagonal element
of the Maurer-Cartan form as $\Omega_j= (\ooo{A_j}{Bj})$

\subheading{4.10.2 Exercises} 

1. Show $\ff 2$ is the complete system with the
base locus
 the disjoint union of the linear spaces
$\ppp (\ee{01}\ot\hdots\ee{0({j }-1)}\ot ( W_{j }/ \ee{0j})
\ot \ee{0({j }+1)}\ot \hdots\ot\ee{0r})
$.

2. Show
$$
\align
|\ff k|=&I_k  ( \cup_j 
 \ppp (\ee{01}\ot\hdots\ee{0({j_1}-1)}\ot (W_{j_1}/
\ee{0j_1})\ot\ee{0({j_1}+1)}\ot
\hdots
\hdots\ee{0({j_2}-1)}
\tag 4.10.2.1 \\
&\ot (W_{j_2}/\{\ee{0j_2}\})\ot\ee{0({j_2}+1)}\ot 
\hdots\ot \ee{0({j_k}-1)}\ot (W_{j_k}/\{\ee{0j_k}\})\ot\ee{0({j_k}+1)}\ot
\hdots
\ot \ee{0r})\endalign
$$
In particular, the last nonzero fundamental form is $\ff r$.

\proclaim{4.10.3 Proposition } Let $X_j\subset\ppp W_j$,  $1\leq j\leq r$,
be varieties and let
  $y=[x_1\ot \hdots \ot x_r]  \in Y:=Seg(X_1\times\hdots\times X_r)$.Then
$$
\align
\hat T_yY&=
\Sigma_j x_1\ot \hdots\ot x_{j-1}\ot \hat T_{x_j}X_j\ot x_{j+1}
\ot\hdots\ot x_r\tag 4.10.3.1 \\
II_{Y,y}&=
\Sigma_{j,k} x_1\ot \hdots\ot x_{j-1}\ot \hat T_{x_j}X_j\ot x_{j+1}
\hdots\ot x_{k-1}\ot \hat T_{x_k}X_k\ot x_{k+1}
\ot\hdots\ot x_r\\
& +
\Sigma_j x_1\ot \hdots\ot x_{j-1}\ot \hat II_{X_j,x_j} \ot x_{j+1}
\ot\hdots\ot x_r\endalign
$$
\endproclaim

\subheading{4.11 Fundamental forms of the Severi varieties}

Here,
it is easier to use the coordinate definition. Choose affine coordinates
based
at $[p]$ where
$$
p = \pmatrix 1& 0 & 0\\
 0 & 0& 0\\
   0 & 0& 0\endpmatrix
$$
and following the notation of (3.1), denote the affine coordinates
$u_1,u_2,u_3\in\Bbb A ,r_2,r_3\in \Bbb C$ where the tangent space to $p$ is
$\{u_1,u_2\}$ (the span is taken over $\Bbb C$).   In these coordinates:
$$
r_2(u_1,u_2)= u_1\bar {u_1}\ \text{as} \  \ \ \text{det}\pmatrix
1&\bar {u_1}\\ u_1 & r_2\endpmatrix =0
$$
$$r_3(u_1,u_2)=  u_2\bar {u_2}\ \text{as} \  \  \ \text{det}\pmatrix
1&\bar {u_2}\\ u_2 & r_3\endpmatrix =0
$$
$$u_3(u_1,u_2)= \bar {u_2}u_1\ \text{as} \  \  \ \text{det}\pmatrix
1&\bar {u_1}\\ u_2 & u_3\endpmatrix =0
$$
where the last equation gives us one, two, four or eight quadratic forms.  
The
determinants come from the vanishing of $2\times 2$ minors that must be
zero to
make the Severi variety consist only of rank one elements. 
In division algebra notation the second fundamental forms
are
$$
  |II| =  \ppp \{ u_1\bar{u_1}, u_2\bar {u_2},\bar {u_2}u_1\}. \tag 4.11.1
$$

Note that since this is all of the Taylor series, there are no
other differential invariants.

\smallpagebreak

\subheading{4.11.2 Remark}The alert reader might have observed that, in all
cases above, the base locus of the second fundamental form corresponds to
the
homogeneous space one obtains by marking the nodes in the dynkin diagram
of $H=P^{ss}$ adjacent
to the nodes removed that correspond to $P$.  This is part of
a general pattern, see [LM1] for details.

\medpagebreak

\subheading{4.12  Proof of the equality in Terracini's lemma}  

Consider the mapping
$$
\align 
\cf^1_Y\times\cf^1_Z& \ra \ppp V\\
((\ee 0\hd \ee N), (f_0\hd f_N)&\mapsto [\ee 0 +f_0]\endalign
$$
The image is $S(Y,Z)$. We compute its tangent space.
Let $\ooo AB, \ett AB$, denote the
entries of the Maurer Cartan forms over $Y$ and $Z$ respectively.
Let $1\leq \a\leq
\tdim Y$,
$1\leq j
\leq
\tdim Z$. Note that
$\hat T_{[\ee 0]}Y= \{ \ee 0,\ee\a \}$, $\hat T_{[f_0]}Z=\{ f_0, f_j\}$.
 Since
we know the line $\ppp \{\ee 0, \eta_0\} =\BP^1_{[\ee 0][f_0]}$ is
contained in
the tangent space to $[\ee 0+f_0]$ we calculate modulo
$\{\ee 0, \eta_0\}$.
$$
d(\ee 0 +f_0) \equiv \oo\a\ee\a + \et jf_j\tmod \{\ee 0, f_0\}
$$
Thus
$\tdim S(X,Z)\leq  \tdim \{\ee 0,f_0, \ee \a , f_j\}$ and
by the inequality case they are equal.  Moreover 
$\tilde T_{[\ee 0+f_0]}S(Y,Z)= \ppp\{\hat  T_{[\ee 0]}Y +
\hat  T_{[f_0]}Z\}$. (One can prove the inequality case
by observing that the forms $\{\oo\a, \eta^j_0\}$ are all
independent.)

\heading \S 5. Varieties with degenerate Gauss images\endheading

The new results in this section are from [AGL], currently under
preparation.

\subheading{5.1 Examples}

\subheading{5.1.1 Joins}    Joins have degenerate
Gauss maps with at least $(k-1)$-dimensional fibers
because
  Terracini's lemma implies that 
the tangent space to $S(Y_1\hd Y_k)$ is constant
along each $\BP^{k-1}_{y_1\hd y_k}$.

Note the special cases of cones, where one of the factors
in the join is linear, and of secant varieties, where
all factors are the same.

\subheading{5.1.2 Tangential varieties}  Let $X=\tau (Y)$ be
a tangential variety. Let $x\in X$ and write $x=[\ee 0+\ee 1]$
with $[\ee 0]\in Y$ and $[\ee 1]\in\tilde T_yY$.
Then $\tilde T_xX=
\tilde T_{[ \ee 0+t\ee 1]}X$ for all $t\neq 0$.
 
 More generally, if $y\in \tau ^k(X):=\cup_{x\in X}\tilde T_x\up k X$,
write $y=[\ee 0+\ee 1 +\ee{\mu_1}+\hdots+\ee{\mu_{k-1}}]$,
then
$\tilde T_y\tau^k(X)$ is constant on the  $\pp k\subset\tau^k(X)$
spanned by $\{ \ee 0,\ee 1 ,\ee{\mu_1}\hd \ee{\mu_{k-1}}\}$.

\subheading{5.1.3 Unions of conjugate lines}  
 Let $Y^{n-1}\subset\pp{n+1}$ be
a variety with a generic second fundamental form. In this case there exist
$n-1$ simultaneous eigen-directions for the second fundamental form 
(To make the notion of eigen-direction
  precise, choose a nondegenerate quadric in $II$ to identify
$T$ with $T^*$ and consider the quadrics as endomorphisms of $T$.
The result is independent of the choices.) 
Let $X^{n}\subset\pp{n+1}$ be the union of one of these families
of embedded tangent lines.  Such lines are called {\it conjugate lines}.

 In this case  $\gamma (X)$ is degenerate with one dimensional
fibers. In fact, the embedded tangent space of $X$ is constant along
the conjugate lines. To see this we calculate the tangent
space of $X$ using frames over $Y$. We restrict to frames of
$Y$ such that the basis of $T_{[\ee 0]}Y$ consists of conjugate directions
and let $1\leq i,j\leq n-1=\tdim Y$.
Adapt frames such that $II (\ee i, \ee j)=\delta_{ij} (\ee n +
\lambda_i\ee{n+1})$. We calculate ($2\leq \rho\leq n-1$)
$$
d(\ee 0+t\ee  1)\equiv
(\oo  \rho + t\ooo\rho 1)\ee\rho +t\oo 1(\ee n+ \lambda_1\ee{n+1})
\tmod \{\ee 0,\ee 1\}
$$
Thus  for generic values of $t$,
$$
\hat T_{[\ee 0+t\ee 1]}X=\{\ee 0,\ee 1,\ee\rho, \ee n+\lambda_1\ee{n+1}\}.
$$

and thus $\hat T_{[\ee 0 + t\ee 1]}X$ is independent of $t$.

Note that if $\ee 1$ is not a conjugate direction and one attempts
the same construction, the resulting variety will not in general
have a degenerate Gauss map.

\proclaim{5.2 Results on Gauss maps}Let $X^n\subset\pp\na=\ppp V$ be a
variety.

i. Consider the distribution 
$\Delta\subset TX$ defined by $\tsingloc\ii_x$.
This distribution is integrable and its integal manifolds
are the fibers of the Gauss map $\gamma$. 

ii. [??] The generic fibers of $\gamma$ are linear spaces.
 
iii. 
[??]
If $\gamma (X)$ is degenerate, then
$X$ is not smooth. More precisely, let $F$ be a general fiber,
then $X$ is singular along a codimension one subset in $F$.

iv. [Zak] If $X$ is smooth, then $\gamma$ is finite.
(iii. only implies $\gamma$ is generically finite.)

v. [Ran] If $\gamma$ is generically finite, then it is finite.
\endproclaim

Although Zak's theorem subsumes iii.   I include it
fopr historical reasons and 
  because the methods of proof are quite different and hold
the possiblity of generalizations in various directions.

\subheading{5.3 Remark: differences between affine and projective spaces}
Ran's arguments are valid in a more general ambient space
than a projective space because he
does not use the linearity of the fiber, but he does
use compactness in an essential way as can be seen by
the  smooth surface in affine space $z= xy^2+ (1-x)y$,
which has the $x$-axis as an isolated fiber.
See [R] for his argument.
$$
\ 
$$

\subheading{5.4 Proofs}

\demo{Proof    of i}
 From the definition of $II$ as the derivative of the 
Gauss map, the first   assertion is   immediate. A proof
is included anyway, as it  will serve as a model
for more complicated proofs and I do not prove the equivalence of the
definitions of $II$.

Adapt frames   to the flag
$$
\hat x\subset \{ \hat x , \tsingloc
\ii_x\}\subset\hat T\subset V \tag 5.4.1
$$
by letting $\{\ee 1\hd\ee f\}\simeq \tsingloc\ii$.
Use additional index ranges $1\leq s,t\leq f$,
$f+1\leq i,j,k\leq n$.
Our adaptations have the effect that $\ooo\mu s=0$.
The first assertion is equivalent to showing the distribution
$\{\oo j\}$ is integrable. Using the Maurer-Cartan equation,
we see 
$$
d\oo j\equiv -\ooo js\ww\oo s \tmod \{\oo i\} \tag 5.4.2
$$
We   examine the forms $\ooo js$. Since $\ooo\mu s=0$,
$$
0=d\ooo\mu s = -\ooo\mu j\ww\ooo js\tag 5.4.3
$$
Since we are assuming $\{\ee s\}$ is the entire singular locus,
for each $j$, there exists some $\mu$ with $\ooo\mu j\neq 0$
which implies
  $\ooo js = C^j_{sk}\oo k$ for some functions $C^j_{sk}$.
Thus $d\oo j\equiv 0 \tmod \{\oo i\}$ and the Frobenius theorem
implies the distribution is   integrable.\qed\enddemo

\demo{Proof of ii}
To see the fiber is a linear space, it will suffice to show that 
$II_{F,[\ee 0]}=0$. On $F$, $\oo j=0$ so we have
$$
\align
d\ee 0 &\equiv \oo s\ee s\tmod \ee 0  \tag 5.4.4 \\
\underline d^2\ee 0 &\equiv \ooo js\oo s\ee j\tmod\{\ee 0,\ee s\}\\
&\equiv 0\tmod\{\ee 0,\ee s\}\text{ and }\tmod\{\oo j\} \qed \endalign
$$
\enddemo

 \demo{Proof of iii} Fix a fiber $F$  of $\gamma$ and let
  $p=u^0\ee 0+u^s\ee s\in F$. We calculate $\tilde T_pX$. Since
$F\subset\tilde T_pX$ for all $p\in F$ we can work modulo
$F=\ppp \{\ee 0,\ee s\}$.
$$
\align
dp\equiv & (u^0\oo j + u^s\ooo js )\ee j
\tmod \{\ee 0,\ee s\}\tag 5.4.5\\
\equiv & (u^0\delta^k_j + u^sC^k_{sj})\oo j\ee k\endalign
$$
We may think of $[u^0\hd u^f]$ as parametrizing a
$(f-1)$-dimensional (projective dimension) family 
of   matrices, which
must drop rank at least along a codimension one subset, which corresponds
to
$(X\cap F)_{sing}\subseteq X_{sing}$.
\qed\enddemo

 An important feature of this proof is that we were able to
calculate the tangent space to $F$ using first and second order
information at just one point of $F$, namely $[\ee 0]$. We will
see this type of calculation reappear when studying $\tx$ and $X^*$.

 \demo{Second proof of iii}  Let $F$ denote a typical  fiber of the
Gauss map of $X$. 
We have
$$
\CD
F @>>>
X \\ 
 @. @VVV \\ 
  @.
\gamma (X) 
\endCD\tag 5.4.6
$$

which determines a map
$$
\align
\gamma (X)&\ra \Bbb G(k,\ppp V)\tag 5.4.7\\
\gamma (x) &\mapsto F_{\gamma (x)}=\gamma\inv \gamma (x)\endalign
$$
The differential of this mapping  is a linear
map $T_{\gamma (x)}\gamma (X)\ra \hat F^*_{\gamma (x)}\ot (\hat T_xX/\hat
F_{\gamma (x)})$ which can be thought of as a map
$\hat F_{\gamma (x)}\ra T_{\gamma (x)}{\gamma (X)}^*\ot T_xX/\hat F_{\gamma
(x)}$, from a $k$-dimensional vector space to the space of
$(n-k)\times (n-k)$ matrices. If $T_xX$ is to have 
the proper dimension, the image matrix must be of full rank
for all $f\in F_{\gamma (x)}$, but this is impossible,
one cannot have a linear space of square matrices such that all
are of maximal rank, in fact the matrices must drop rank along a
codimension one subset.\qed\enddemo

\subheading{5.5 Exercise} Prove that a smooth variety $X$ cannot be ruled
by $\pp k$'s for $k\geq a$ by considering the reduced Gauss map
$B^{n-k}\ra\BG (k,\ppp V)$
where $B$ is the base space. A generalization of this result is given
in \S 6.

\heading \S  6. Smoothness of ruled and uniruled varieties\endheading

 In this section I discuss codimension restrictions on ruled
and uniruled varieties. Ruled varieties are special cases of
varieties that can be described as fibrations,
   and the bound
on their codimension is the same for that of an
arbitrary  fibered variety:

\proclaim{6.1 Theorem, Remmert and Van de Van [RV]}
Let $X^n\subset\pp\na$ be a smooth variety that is a fibration
with fibers
of dimension $f$, then $f< a$.\endproclaim

The proof of Remmert and Van de Van's theorem is a simple application
of the intersection property of projective space. Let $p\in B$ be
a point and let $Z\subset B$ be a hypersurface not containing $p$.
$\tdim \pi\inv (p)=f$ and $\tdim\pi\inv (Z) = (n-f-1)+f=n-1$. Thus
if $f+n-1> n+a$, i.e. $f\geq a$, $\pi\inv (p)\cap \pi\inv (Z)\neq\emptyset$
which is a contradiction.

If $X^n$ is only uniruled by $k$-planes, it is easier to be smooth,
for example the smooth quadric hypersurface $Q^n$ is uniruled
by $[\frac n2]$-planes. Nevertheless
there are still restrictions on the codimension:

\proclaim{6.2 Theorem [Ran, R3]}
Let $X^n\subset\pp\na$ be a smooth variety
such that through a
 general point
  $x\in X$ there is at least one $\pp k\subset X$ containing $x$.
Then $a\geq k/(n-k)$.
\endproclaim

Ran's proof follows from his discussion of
generalized Gauss maps. He also remarks that the
results can be derived from the Barth-Larsen theorems.
Here is an alternative
proof that is local in character, in that the only way global information
is taken into account is that the singular locus of
the second fundamental form must be empty:
 
\demo{Proof}
We are given that at each point there is a $\pp {k-1}$
contained in $\tbaseloc\ii$. The proof will follow from
the following lemma:\enddemo

\proclaim{6.3 Lemma} Let $A$ be an $a$-dimensional system of quadrics
  on an $n$-dimensional vector space $V$
such that there is  a linear space
$W$  of dimension $k$ in
the base locus of $A$. If
$a<k/(n-k)$, then $A$ has a singular locus.
\endproclaim

\demo{Proof}
Any quadric $q\in A$ can be written
$$
q= v^1w^1+\hdots + v^kw^k + q'
$$
where,   $q'\in S^2W\upperp$,
and $v^j\in W\upperp$. Since $k>n-k$, at most 
$n-k$ of the $v^j$ are independent.
Thus each quadric has at least an
$k-(n-k)=2k-n$ dimensional singular locus in $W$, so
if $a(n-k)<k$ then $\tsingloc (A)\neq 0$.\qed\enddemo 

\subheading{6.4 Problem} Use the rank restriction theorem to get a
better bound. (Although among bounds of the form above, Ran's bound is
optimal as equality cases occur.)

\heading\S 7. Varieties with degerate dual varieties\endheading

Let $X^n\subset\pp\na=\ppp V$ be a variety,
let $X^*\subset\ppp V^*$ denote its dual variety, and let
$\dstar = \dstar (X) = n+a-1-\tdim X^*$ denote the dual defect,
as discussed in the introduction.

\subheading{7.1 Examples}

i. The smooth quadric hypersurface is self dual.

iii. Let $X= Seg(\ppp W_1\times\ppp W_2)=
Seg(\pp k\times\pp l)$, with $k\geq l$,
 then $X^*=\sigma_{l}(Seg(\ppp
W_1^*\times\ppp W_2^*))$ and thus $X^*$ is degenerate iff $k\neq l$
with defect $\dstar = k-l$

iii'. More generally, if $X=Seg(\pp{k_1}\times\hdots\times\pp{k_r})$
where $k_1\geq k_2\geq\hdots\geq k_r$, then $X^*$ is degenerate
iff $k_1> k_2+\hdots +k_r$
with defect $\dstar = k_1-( k_2+\hdots +k_r)$.

iv. If $X$ is a scroll, that is a linear fibration
with base a curve, then $\dstar= n-2$.

vi. If $X=G(2,W)\subset\ppp\La 2 W$, then 
$X^*=\sigma_p(G(2,W^*)$, where
$p=\frac 12(k-2)$ if $k$ is even and $\frac 12(k-3)$ if $k$ is odd, thus
$\dstar= 0$ if $k$ is even
and $2$ if $k$ is odd.  

vii. $\Bbb S_5(V,Q)^*=\Bbb S_5(V^*,Q)$,  $\dstar =4$. (Higher
Spinor varieties have nondegenerate duals).

Note that the example of the Segre shows that there can be
no absolute bound on $\dstar$. 

\proclaim{7.2 Theorems on dual varietes} Let $X^n\subset\pp\na=\ppp V$ be a

  variety and let $X^*\subset\ppp V^*$
denote its dual variety. Let $\dstar= n+a-1-\tdim X^*$ denote
the dual defect of $X$.   

i. [Bertini] $(X^*)^*=X$  ({\text the reflexivity theorem})

ii.  $\dstar\geq 
\delta_{\gamma}:= \tdim X- \tdim \gamma (X)$.

Assume $X$ is not contained
in a hyperplane. 

iii.  

\ \ a. If $X$ is a smooth hypersurface, $X^*$ is nondegenerate.

\ \ b. If $X$ is any curve, then $X^*$ is nondegenerate.

\ \ c. If $X$ is any surface, then $X^*$ is nondegenerate.

iv. If $X$ is ruled by $k$-planes,
i.e. if 
$X$ is a linear  fibration, with
$k$-dimensional fibers, then $\dstar\geq 2k-n$.

v. If $X$ is a cone over a linear space $L$, then $X^*\subset L\upperp$.
 
vi. [Bertini] If $H\in 
X^*_{sm}$, then $ 
\{ x\in X\mid \tilde T_xX\subset H\}$ is a linearly embedded
$\pp\delta\subset\ppp V$.

vii. [Zak, Z] If $X$ is smooth, then
 $\tdim X^*\geq \tdim X$, i.e. $\dstar\leq
a-1$.

viii. [Landman, E1]If $X$ is smooth, then $n-\dstar$ is even. In particular

$\dstar = n-1$ is impossible.

ix. [Ein, E1] If $X$ is smooth, $\tdim X=\tdim X^*$,
$X\not\subset \pp{\na -1}$,  and $a\geq\frac n2$, then
$X$ is one of the following:

\ \ \ \ a. $Seg(\pp 1\times\pp{n-1})\subset\pp{2n-1}$.

\ \ \ \ b. $G(2,5)\subset\pp 9$.

\ \ \ \ c. $\Bbb S^{10}_5\subset\pp {15}$.

Moreover, in these cases  $X^*$ is isomorphic to $X$.

x. [Ein, E1] [BC]. If $a=2$, $X\not\subset \pp{\na -1}$, and $X$ is
smooth, then $X^*$ is a hypersurface 
unless $X=Seg(\pp 1\times\pp 2)\subset\pp 5$.

xi. [Ein, E2]. If $\dstar= n-2$ then $X$ is a scroll.  

xii. [Ein, E2]. If $\dstar\geq\frac n2$ then $X$ is ruled by
$\pp{\frac{\dstar
+n}2}$'s.

xiii. [Ein, E2, 1.3]. 
If $X=Y\cap Z$ with $Y$ a smooth hypersurface of degree greater
than one, then $X^*$ is
nondegenerate. 

xiv. If $H$ is a hyperplane,
then $\dstar (X\cap H) = \dstar (X) -1$.
\endproclaim

Aluffi   has proved a partial converse to vi, see [Aluffi].

There are also classification results in small dimension see
[E2] and [BFLS].

Landman's original (unpublished) proof of viii. used Lefshetz pencils.  
 
Ein's   proofs rely on studying 
the conormal bundle of $Y_H:=\{x\in X\mid \tilde T_xX\subset H\}$
and its deformations.   
 
I do not know of a direct 
geometric proof of vi, both Ein and Ballico and Chianti use facts
about rank two vector bundles on projective space.

\proclaim{7.3 Theorems on the projective
differential geometry of $X^*$}
Let $X^n\subset\pp\na$ be a variety and let
$x\in X$ be a general point.  Let $\dstar= n+a-1-\tdim X^*$ denote the dual
defect of $X$. 
Given $H\in X^*$, let $Y_H=\{ x\in X\mid \tilde T_xX\subset H\}$.

i. (??) [GH], [L3]   $|II_{X,x}|$
is a system of quadrics of projective dimension $a-1$
 and   bounded rank $n-\dstar$ on an $n$-dimensional vector space.
I.e. If $r$ is the rank of a generic quadric in $|II_{X,x}|$,
then $\dstar= n-r$.  

ii. [IL, 3.4] If $X$ is smooth and $H\in 
X^*_{sm}$, then $|II_{X^*,H}|$ is a system of quadrics of
projective dimension $\dstar$ and constant rank $n-\dstar$
on an $(n-\dstar + a-1)$-dimensional vector space.
I.e. If $r$ is the rank of the quadrics in $|II_{X^*,H}|$,
then $\dstar= n-r$.

iii. [IL]  If $H\in X^*$ is a smooth point,
then $\ii_{X^*,H}$ can be recovered
from
$II_{X,x}$ as $x$ ranges over $Y_H$.
(The precise inversion formula is given below.)

iv. [IL]  If $H\in X^*$ is a smooth point,
then $\ii_{X^*,H}$ can be recovered
from
$II_{X,x}$ and $F_{3 X,x}$ where $x\in Y_H$ is any smooth point.
(The precise inversion formula is given below.)

v. [L7] Let   $X^n\subset\pp\na$
 be an analytic open subset of a variety . 
 If $a\geq 2$,
$X$ is not contained in a hyperplane, and
for    general 
$x\in X$ there exists  a linear space $L^{n-1}_x$
osculating to order two at $x$, then   $L_x\subset \overline X$.
In particular, any smooth open subset of a variety  having the second
fundamental form
of   a scroll is a scroll.

vi. [L7]  
 Let   $X^n\subset\pp\na$
 be an analytic open subset of a variety 
having the property that
through  a general 
point $x\in X$, there exists a linear space
$L_x$, of dimension $k$,
osculating to order two at $x$. Then there exists
a linear subspace $M_x\subset L_x$, of dimension $2k-n$, such that  
$M_x\subset \overline X$. In particular $\dstar  \geq 2k-n$.

vii. [IL]  If  there exists $x\in\xsm$ such that $Y_x\cap X^*\subset
X^*_{sm}$. Then $\dstar\geq a-1$.

viii.[IL]   If  there exists $x\in\xsm$ such that $Y_x\cap X^*\subset
X^*_{sm}$ and   $H\in X^*_{sm}$ such that  $Y_H\cap X\subset\xsm$,
   Then $\tdim X=\tdim
X^*$.
\endproclaim

iv  may be considered as an infinitesimal version  of  7.2.vii.
7.3.v is related to, but different from 7.2.viii.

\subheading{ 7.4 Frames for $X$ and $X^*$}

 Let $H\in X^*_{sm}$ and let 
$$
Y_H=\{ x\in X\mid \tilde T_xX\subseteq H\}
$$

Let 
$$
\align 
\ci^0 & = \{ (x,H)\mid x\in \xsm, H\in X^*_{sm}, \tilde T_xX\subseteq H\}\\
&=\{ (x,H)\mid x\in \xsm, H\in X^*_{sm}, \tilde T_HX^*\subseteq x\}\\
&=\{ (x,H)\mid x\in \xsm, H\in X^*_{sm}, x\in Y_H\}\endalign
$$
and let
  $\cf^*\ra \ci^0$ be the
 frame bundle
of bases
of $V$ over $(x,H)\in \ci^0$ adapted to the flag:
$$
0\subset \hat x\subset
 \hat T_x Y_H\subset\hat T_xX\subset \hat H\upperp\subset V.\tag 7.4.1 
$$

%%%%%NOTE change in definition here!!
\def\nl{n+\lambda}
\def\nk{n+\kappa}

To describe $\cf^*$ using indices,
let $1\leq \alpha,\beta\leq n$,
$1\leq s,t,u\leq n-r$,
 $n-r+1\leq i,j,k \leq n$, $ 1\leq\kappa,\lambda \leq  a-1$.
Require that
$$
\align   \hat x& = \{\ee 0\}\\
\hat T_xY_H &=\{ \ee 0,\ee s\}\\
\hat T_xX&= \{\ee 0,\ee s,\ee j\}=\{ \ee 0,\ee \a\}\\ 
\hat H&= \{ \ee 0,\ee \a\ee\nl\}\endalign
$$

In indices, the  flag (7.4.1) is
$$
\{\ee 0\}\subset\{\ee 0,\ee s\}\subset\{\ee 0,\ee s,\ee j\}
=\{\ee 0,\ee\alpha \}\subset \{\ee 0,\ee\alpha,\ee{n+\lambda}\}\subset
\{\ee 0,\ee\alpha,\ee{n+\lambda},\ee\na\} =V\tag 7.4.2
$$

Write the pullback  of the Maurer-Cartan form
to $\cf^*$ as:
$$
\Omega = \pmatrix \ooo 00&\ooo 0 t&\ooo 0 k&\ooo 0 \nl
&\ooo 0\na\\
\ooo s0 & \ooo sk & \ooo st& \ooo s \nl& \ooo s\na\\
\ooo j0 & \ooo jk & \ooo jt& \ooo j \nl& \ooo j\na\\
0 & \ooo  \nk t & \ooo  \nk k& \ooo \nk \nl & \ooo \nk\na\\
0 & 0 & \ooo \na k& \ooo\na \nl& \ooo\na\na \endpmatrix
\tag 7.4.3
$$
where
$
\ooo\na t=0  
$ 
because 
we have adapted such that $\{\underline e_t\} = \tsingloc ( q^{\na})$,
which proves (7.3.i).

We first prove $Y_H$ is a linear space, proceeding as we did
in \S 5 to show that the fibers of the Gauss map are linear spaces.
(Alternatively one can argue that if $H\in X^*_{sm}$, then
$\rho\inv (H)=\ppp N^*_HX^*(1)=\pp\delta$ and then that
$\pi$ linearly embedds the $\pp\delta$.) First note that
$$
0=d\ooo\na s =-\ooo\na k\ww\ooo ks\tag 7.4.4 
$$

The forms $\ooo \na k=\qq\na ki\oo i$ are all independent
 since   $q^{\na}$ descends to be a nondengenerate quadratic form
on $T_xX/T_xY_H$,
and $\hat T_xY_H\simeq N^*_HX^*(1)$ is of constant dimension.
 Thus
$\ooo ks\equiv 0\tmod\{\oo j\}$.

We  calculate $II_{Y_H, [\ee 0]}$:
$$
\align 
d\ee 0&\equiv \ooo s0 \tmod \ee 0\text{ and }\tmod\{\oo j\}\tag 7.4.5 \\
\underline d^2\ee 0 &\equiv \ooo js\oo s
\tmod\{\ee 0,\ee s\} \text{ and }\tmod\{\oo j\}\\
&\equiv 0 \tmod\{\ee 0,\ee s\} \text{ and }\tmod\{\oo j\}
\endalign
$$

 Pulled back to $Y_H$, $\oo j=0$, which implies    
 $II_{Y_H,[\ee 0]}=0$ andthus  $Y_H$ is a linear space.
 
\medpagebreak 

We may also consider $\cf^*$ as an adapted frame bundle over  $X^*_{sm}$.
It is clearly zero-th  order adapted, but is in fact as adapted to $X^*$
as it is to $X$. To see this, note that
  $\cf^*$ is also adapted to the dual flag:
$$
0\subset  \hat H  \subset (\hat T_xX)\upperp
\subset (\hat T_xY_H)\upperp\subset
\hat x\upperp
\subset
V^*.\tag 7.4.6 
$$
But this flag is just
$$
0\subset \hat H\subset\hat Y_x\subset \hat T_HX^*
\subset \hat x\upperp \subset V^*.\tag 7.4.7
$$

To consider $\cf^*$ as an adapted frame bundle over
  $X^*_{sm}$, let
$\{\uee B\}$ denote the dual basis to $\{\ee B\}$.
Let $<,>$ denote the pairing $V\times V^*\ra \bccc$,
so $<\ee A, \uee B>=\delta^B_A$. We calculate:
$$
0=d<\ee A, \uee B> =<\ooo AC \ee A , \uee B>+ <\ee A, d\uee B >\tag 7.4.8
$$
which implies
$$
d\uee B = -\ooo BC\uee C.   \tag 7.4.9
$$

We may specify a point of $\cf^*$ using  the dual basis 
$\uee B$ and it will
be convenient to do so when studying frames over $X^*$.
 From this
dual  perspective,
 $H=[\uee\na ]$, $  Y_x=\ppp\{\uee\na, \uee\nl\}$, $\tilde
T_HX^*=\ppp\{\uee\na,\uee 
\nl,\uee k\}$,
$x=\ppp\{\uee\na,\uee\nl,\uee \alpha\}$.
  
To see  $\tilde T_HX^*=\ppp\{\uee\na,\uee
k,\uee \nl\}$, note that
$$
d\uee\na \equiv -\ooo\na k\uee k -\ooo\na \nl\uee \nl \tmod \uee\na .
   \tag 7.4.10
$$
 
(7.4.10)  shows that $\{\ooo\na k,\ooo\na \nl\}$ form a basis
of the semi-basic forms for the projection to $X^*$, note that
 $\tdim \{\ooo\na
k \} =r$. It also shows that on $\cf^*$, the forms
$\ooo\na{n+\lambda}$ are independent and independent of the forms $\oo j$
because the dimension of the tangent space to $X^*$ is
constant at all smooth points of $X^*$.
 If $X$ is smooth, then (7.4.10) is valid for {\it any}
$H\in X^*_{sm}$.

\smallpagebreak

Using (4.4),  we   compute $II_{X^*,H}$:
$$
\align 
\ud 2 \uee\na
&\equiv (\ooo\na k\ooo k s + \ooo \na \nl\ooo \nl s )\uee s
+\ooo\na k\oo k\uee 0 \tmod\hat T_HX^*    \tag 7.4.11\\
&\equiv
(\rr\na s ij\oo j  + 2\qq \nl si\ooo\na \nl)\oo i\uee s +
\qq\na jk\oo j\oo k\uee 0  \tmod\hat T_HX^* \endalign
$$
where to derive the second line we utilized: 
$$
\align
&
\rr\na st\beta\oo \beta = 0 \tag 7.4.12\\
&
\rr\na sij\oo j = -\qq\nl si\ooo\na\nl +\qq\na kj\ooo ks. \tag
7.4.13\endalign
$$

The second line of (7.4.11) proves (7.3.iii), (7.3.iv) and
(7.3.ii); the first
because
$Q_x$ is the quadric corresponding to $\uee 0$
and the second because all the other invariants appearing are
coefficients either of $II_{X,x}$ or $F_{3X,x}$. The first inversion
formula in turn proves [IL, 3.4] because  the rank of the quadric
corresponding to $\uee 0$ is of rank $r$, but as one moves in $N^*_HX^*$,
every conormal direction gets a chance to be $\uee 0$ (under the hypothesis
$H\cap X_{sing}=\emptyset$).

Note that (7.4.10), (7.4.11) imply that  the dual
variety of a smooth variety has the property that an
arbitrary smooth point of $X^*$ is \lq\lq general\rq\rq\ to order two.

By a classical theorem (see the next section), it is not
possible to have a system of quadrics of constant rank if $r$ is odd,
which, combined with
(7.3.ii) provides a new proof of the Landman parity theorem (7.2.viii).

[IL, 3.4] led us to study systems of quadrics of  bounded
and constant rank. Our main result 
on systems of quadrics is the resolution of the constant even rank
problem:

\proclaim{7.4.14 Theorem [IL, 2.16]} If $r$ is even   then
$$\tmax \{ \tdim(A) \mid A \subset S^2\bcc m \text{ is of constant rank } r
\}
 = m-r+1.$$
\endproclaim

[IL, 2.16] immediately provides a new proof of Zak's theorem
that $\dstar\leq a-1$,
without using the Fulton-Hansen connectedness theorem
but still using the   Lefschetz theorem as 
will be explained in \S 8.

The discussion above proves the rank restriction theorem (0.3.3.4)
for the case of a generic quadric, and the case for an arbitrary quadric
follows from counting how dimensions drop, see [L3].

 We may state (7.3.iv) in the more
precise form:

\proclaim{ 7.4.15 Inversion formula, [IL]}With respect to the basis
$\{\ooo\na\nl,\oo j \}$ of $T^*_HX^*$,
and $\uee 0,\uee s$ of $N_HX^*$,
 we have
$$
II_{X^*,H}=\left\{ Q_0 =\pmatrix 0&0\\ 0&\qq\na jk\endpmatrix ,\ \ 
Q_s=\pmatrix 0&\qq\nl sj \\ \qq\nl sk &\rr\na sjk\endpmatrix \right\} 
$$
where the blocking is $(a-1,r)\times (a-1,r)$.
\endproclaim

 Regarding (7.4.10), we made the following observation:

\proclaim{7.4.16 Observation } If $\frac r2<a-1$,
i.e. if $\dstar <n-2(a-1)$,  then 
$II_{X^*,H}$ cannot be of constant rank $r$ unless the $\delta +1$
matrices
$(\qq\na jk ),(\rr\na 1jk )\hd
(\rr\na \delta jk )$ are all linearly independent.
\endproclaim

Were the matrices dependent,  there would be
a quadric  in $II_{X^*,H}$ whose lower right $r\times r$ block was
zero, but if $\frac r2<a-1$, such a quadric could not have rank $r$.

Observation (7.4.16) motivated the  following conjecture:

\proclaim{7.4.17 Conjecture} If
$X^n\subset\pp\na$ is a smooth variety with $\dstar >0$, then $\dstar \geq 
n-2(a-1)$.
\endproclaim

\proclaim{7.4.18 Corollary of Conjecture} If
$X^n\subset\pp\na$ is a smooth variety with $a-1<\frac n3$, then $\dstar
=0$.
\endproclaim

 The motivation for the conjecture
is that if $\dstar <n-2(a-1)$, there are {\it
genericity} conditions placed on the coefficients of $F_3$
that appear to be incompatible with the closed conditions
arising from the degeneracy of $X^*$.

\heading\S 8. Systems of bounded and constant rank\endheading
 
\subheading{8.1 Examples in coordinates}  $A$ refers to symmetric
systems, $B$ to linear systems, and $C$ to skew symmetric systems.
In all cases
$[\ee 0\hd\ee\delta]\in\pp\delta$.
Consider 
$$
B_I.\
\pmatrix
    \aa 0 & 0\\
  \aa 1 & \aa 0\\
  \aa 2 & \aa 1\\
   0 & \aa 2\endpmatrix
\ \ \ \ \ 
C_{II}. \
\pmatrix 
   0 &\aa 0 & \aa 1\\
 -\aa 0 &0 & \aa 2\\
 -\aa 1 &-\aa 2 & 0\endpmatrix
$$

$$
A_I.\
\pmatrix 
 0&0&0 &0 &\aa 0 & 0\\
0&0&0 &0 &\aa 1 & \aa 0\\
0&0&0 &0 &\aa 2 & \aa 1\\
0&0&0 &0 & 0 & \aa 2\\
\aa 0&\aa 1&\aa 2& 0&0&0\\
0&\aa 0&\aa 1& \aa 2&0&0\endpmatrix
\ \ \ \ \ 
A_{II}. \
\pmatrix 
 0&0&0 &0 &\aa 0 & \aa 1\\
0&0&0 &-\aa 0 &0 & \aa 2\\
0&0&0 &-\aa 1 &-\aa 2 & 0\\
0&-\aa 0&-\aa 1 &0 &0 & 0\\
\aa 0&0&-\aa 2& 0&0&0\\
\aa 1&\aa 2&0 & 0&0&0\endpmatrix
$$

$$
A_{III}. \
\pmatrix
 0     &0      &0      &0      & 0    &0      &0      &\aa 0&\aa 1&\aa 2\\
 0     &0      &0      &0      & 0    &-\aa 0 &-\aa 1 &0&0&\aa 3\\
 0     &0      &0      &0      &-\aa 0&0      &-\aa 2 &0&-\aa 3 &0\\
 0     &0      &0      &0      &-\aa 1&-\aa 2 &0      &-\aa 3 & 0&0\\
0      &   0   &-\aa 0 &-\aa 1 & 0    &0      &0      &0&   0&\aa 4\\
0      &-\aa 0 &0      &-\aa 2 &     0&0      &0      &0&-\aa 4&0\\
0      &-\aa 1 &-\aa 2 &0      &     0&      0&0      &-\aa 4&0&0\\
\aa 0  & 0     &0      &-\aa 3 &    0 & 0     &-\aa 4      &0&0&0\\
 \aa 1 &0      &-\aa 3 & 0     &     0& -\aa 4&0      &0&0&0\\
 \aa 2 &\aa 3  &0      & 0     & \aa 4&      0&0      &0&0&0\endpmatrix
$$

$$
C_{IV}.\
\pmatrix 
0& 0& 0& 0& 0& 0& 0& \aa 0& \aa 1& 0\\
 0& 0& 0& 0& 0& 0& \aa 0& \aa 1& 0& \aa 2\\
 0& 0& 0& 0& 0&-\aa 0& \aa 1& 0& \aa 2& \aa 3\\ 
 0& 0& 0& 0& \aa 0& \aa 1& 0& \aa 2& \aa 3& 0\\ 
 0& 0& 0&-\aa 0& 0& 0& \aa 2&-\aa 3& 0& 0\\ 
 0& 0& \aa 0&-\aa 1& 0& 0& \aa 3& 0& 0& 0\\ 
 0&-\aa 0&-\aa 1& 0&-\aa 2&-\aa 3& 0& 0& 0& 0\\ 
-\aa 0&-\aa 1& 0&-\aa 2& \aa 3& 0& 0& 0& 0& 0\\ 
-\aa 1& 0&-\aa 2&-\aa 3& 0& 0& 0& 0& 0& 0\\
 0&-\aa 2&-\aa 3& 0& 0& 0& 0& 0& 0& 0\endpmatrix
\ \ A_{IV}. \ \pmatrix 0&C_{IV}\\ {}^tC_{IV}&0\endpmatrix
$$

\medpagebreak

\subheading{8.2 Examples via standard constructions}

\subheading{8.2.1 Doubled systems}Given a linear subspace
$B\subset \bcc k\ot \bcc l$
  of constant rank
$p$,
  one can form a system of   symmetric or skew,
maps of rank $2p$ in  
$S^2\bcc{k+l}$ (or $\La 2 \bcc{k+l}$), namely
$$
\pmatrix 0 &B\\ {}^tB& 0\endpmatrix\ \ \ 
\pmatrix 0 &B\\ -{}^tB& 0\endpmatrix .
$$
We call such systems, {\it doubled} systems.
Examples $A_{I},A_{II},A_{IV}$ are doubled systems.

\subheading{8.2.2 Systems of split type}
A suitably generic subspace $B\subset \bcc l\ot \bcc k$ of dimension
at most $l-k-1$
  ($k\leq l$) is a 
system of linear maps of constant rank $k-1$. 
We will refer to  the doublings of such systems  as  
{\it systems  of split type}. The name will be explained below.

\subheading{8.2.3 Systems of graded algebra type}
Consider a vector space $V$ and the
inclusion
$$
\align
V&\ra \thom (\La k V,\La {k+1}V)\tag 8.2.3.1 \\
v&\mapsto (E\mapsto v\ww E)\endalign
$$

We will call these systems, {\it systems of graded algebra type}.
One can also use Clifford multiplication in the construction.

One can double these systems to 
get symmetric or skew systems. However, it is not always
necessary to double:
If $\tdim V=2k+1$, then $\tdim \La k V=\tdim \La{k+1}V$ and
we may identify them using a volume form. If
$k$ is even, the maps are symmetric, and they are skew symmetric if $k$ is
odd (see [IL]).
Examples $B_{II}$ and $A_{III}$ arise  in this way.

We can also consider systems of bounded rank via
$$
\La p V\ra \thom (\La k V, \La {k+p}V)
$$
and one can take subsystems of these 
to obtain systems of constant rank. 

\subheading{8.2.4 A mystery}Example $C_{IV}$ is due to Westwick
[W2]. I do not have a
geometric explanation for it, although   it might be related
to the Horroks-Mumford bundle.

\subheading{8.2.5 Examples arising from second fundamental forms of
dual varieties} The   varieties 
$Seg(\pp 1\times\pp n), G(2,5), \Bbb 
S_5$ are self-dual so
their second fundamental forms are of constant rank: 
$$
\align 
&\ii_{Seg(\pp 1\times\pp n)^*, x\ot y}= I_2(\ppp (\bcc {}\ot y )\sqcup
\ppp (x\ot\bcc n ))= I_2(\pp 0\sqcup \pp{n-1})\\
&\ii_{G(2,5), E}= I_2(Seg (E^*\ot V/E))=
I_2( Seg(\pp 1\times\pp 2)) \\
&\ii_{\Bbb S_5, E }=I_2(G(2,E))=I_2(G(2,5))\endalign
$$

 From the pattern, one might guess that $I_2( \Bbb S_5)$ might be
of constant rank. This is not the case, although it is nearly
of constant rank. We have already seen that it 
occurs as the system associated to the second fundamental
form of $\Bbb O\pp 2\subset\ppp\Cal H$
 and that only two ranks are possible
for the quadrics.

I discuss the problems of arbitrary linear systems  and skew-symmetric
systems
of constant rank  
as well as symmetric systems of constant rank because the
three are related.

\proclaim{8.3 Theorems on systems of constant and bounded rank} Let
$$
\align
 &l(r,m,n)=\tmax\{ \tdim (A) \ | \ A\subset \bcc n\ot \bcc m
\text{ is of constant rank }r\}\\
&\underline l(r,m,n)=\tmax\{ \tdim (A) \ | \ A\subset \bcc n\ot \bcc m
\text{ is of bounded rank  }r\}\\
&\overline l(r,m,n)=\tmax\{ \tdim (A) \ | \ A\subset \bcc n\ot \bcc m
\text{ is of   rank bounded below by }r\}\\
\endalign
$$
and similarly let $c(r,m)$ (resp. $\lambda(r,m)$) etc... denote the
corresponding numbers for symmetric (resp. skew symmetric) matrices.

i.  \cite{W} Suppose $2 \le r \le m \le n$. Then
\roster
\item $l(r,m,n) \le m + n - 2r +1$
\item $l(r,m,n) = n - r + 1$ if $n-r+1$ does not divide $(m-1)! / (r-1)!$
\item $l(r,r+1,2r-1) = r +1$
\endroster

ii. (Classical, see [IL, 2.8])
\roster
\item $\overline l(r,m,n) = (m-r)(n-r)$.
\item $\overline c(r,m) = \binom{m-r+1}2$.
\item $\overline \lambda(r,m) = \binom{m-r}{2}$ ($r$ even).
\endroster

iii. (Classical, see [IL, 2.10])
\roster
\item If $0 <r \le m \le n$ then $l(r,m,n) \ge n - r + 1$.
\item If $r$ is even  
then $c(r,m) \ge m-r+1$.
\item If $r$ is even  
then  $\lambda(r,m) \ge m-r+1$.
\endroster

iv. (Classical, see [IL, 2.15]) If $r$ is odd then $c(r,m) = 1$.

v. [IL, 2.16] If $r$ is even  then $c(r,m) = m - r + 1$.
\endproclaim

To prove iii, 
take  $A^{n-r+1}\subset \bcc n\ot\bcc m$
generic to show iii.1, and double symmetrically and skew-symmetrically
to get iii.2 and iii.3.

  I recently solved a conjecture of C. Pauly, characterizing
the maximal linear spaces of skew symmetric
matrices of bounded rank:

\proclaim{8.3.1 Proposition} If $L\subset\La 2\bcc m$ is of bounded rank
$m-1$ and dimension $\binom{m-1}2$, then $L=\La 2\bcc{m-1}$.\endproclaim

If the dimension of $L$ is not maximal, then the conclusion does not
necessarily hold.

\subheading{8.4 Direct geometric methods}
 There is a natural geometric way to 
study
 systems of bounded and constant rank.  
Consider $\sigma_r(Seg(\pp {m-1}\times\pp{n-1}))$, 
$\sigma_r(v_2(\pp{m-1}))$, and $\sigma_{\frac r2}G(2,m)$.
The problem
of finding systems of bounded rank is to find linear spaces
on these varieties, the problem of finding constant rank systems is to find
linear spaces contained in the smooth locus of these varieties,
and the problem of finding systems with rank
bounded below  is to find linear spaces on the ambient projective
space avoiding these varieties. To prove ii,   one 
 calculates the dimensions of the varieties and use Bezout's theorem.

\subheading{8.4.1 Problem} Explicitly identify the examples above as
linear spaces on\linebreak 
$
\sigma_r(v_2(\ppp T))\backslash
\sigma_{r-1}(v_2(\ppp T))$, 
$\sigma_{\frac r2 }(G(2,T))\backslash \sigma_{\frac r2-1}(G(2,T))$, 
$
\sigma_r(Seg(\pp m\times \pp n))\backslash
\sigma_{r-1}(Seg(\pp m\times \pp n))$.

For example, the doubling of a system $B\subset U\ot W$
with $\tdim U=k, \tdim W= m-k$ to $V= U\oplus W$, must
be contained in 
$$
T_{[x_1^2+\hdots x_k^2]}\sigma_k (v_2(\ppp V))\cap 
T_{[x_{k+1}^2+\hdots x_{m}^2]}\sigma_{m-k} (v_2(\ppp V))
$$
where $x_1\hd x_m$ is a basis of $V$
such that $U=\{ x_1\hd x_k\}, W=\{ x_{k+1}\hd x_m\}$.

\subheading{8.5 Vector bundle methods}

Let $A\subseteq V^*\ot W$ be a linear subspace. Over the projective
space 
$\ppp A$ are  the vector bundles $ V\ot \Cal O_{\ppp A}$,
$  W\ot \Cal O_{\ppp A}(1)$ and a vector bundle map $\phi$ between
them, 
$$
\align
\phi_{[a]} : V\ot \Cal O_{\ppp A}&\ra W\ot \Cal O_{\ppp A}(1)\tag 8.5.1 \\
v&\mapsto a(v)\ot a^*.\endalign
$$
 The image and kernel of $\phi$ are in
general sheaves, but if $A$ is a system of constant rank, they are
vector bundles. Let $E:=\phi (\Cal O_{\ppp A}\ot V)$ denote
the image vector bundle in this case. If $V^*\simeq W$ and the
system is symmetric or skew, then one deduces that $E\simeq E^*(1)$
by considering the dual sequence.

By Grothendiek's theorem (a proof of which is already in [HP]), any vector
bundle
  $E\ra \pp 1$ splits as a direct sum of line bundles,
$E=\Cal O(a_1)\oplus\hdots\oplus\Cal O(a_r)$. Given a vector
bundle $E\ra\pp m$, we let the {\it splitting type} of $E$
denote the sequence of integers $(a_1\hd a_r)$ on obtains by
restricting $E$ to a general $\pp 1\subset\pp m$. We say $E$ is
{\it uniform}  if the sequence is independent of the $\pp 1\subset\pp m$.

$E:=\phi (\Cal O_{\ppp A}\ot V)$ is uniform of splitting type $E\mid_{\pp
1}=\Cal O^{\frac r2}\oplus\Cal O^{\frac r2}(1)$,
because $E$ is globally generated and $E\simeq E^*(1)$.
This proves   8.3..iv. (8.3.iv. can
also be proved using a normal form, see [HP].)

 \subheading{8.5.2 Examples} 
$$
\align
&E_{G(2, 2k+1)^*}=T\pp 2(-1)^{\oplus k}\\
&E_{Seg(\pp 1\times\pp n)^*}= split\\
&E_{\Bbb S_5}= \La 2(T\pp 4(-1))\endalign
$$

If $A$ is the doubling of  $B$, then 
$E_A= E_B\oplus E^*_B(1)$. (This makes
sense because $\ppp A\simeq \ppp B$.) 

If $A$ is the doubling of 
a suitably generic $B^{\frac r2}\subset\bcc{\frac r2}
\ot\bcc{m-\frac r2}$, then   $E_A= \Cal O^{\frac r2}\oplus  \Cal O^{\frac
r2}(1)$, justifying the terminology \lq\lq split type\rq\rq .  

\smallpagebreak

Relating the systems to dual
varieties, we have $\ppp A= \ppp N^*_HX^*=\pp\dstar$,
$E= N^*_{X/Y_H}$. Ein observed that $E\simeq E^*(1)$, but he failed
to observe that the identification arises from a symmetric map, although
it is implicit in his proof.

We actually proved a more general result than
[IL, 2.16], namely:

\proclaim{8.5.3 Theorem [IL,1.2]} Let $Z$ be a nonsingular simply connected
projective variety of dimension $\delta$, $E$ a rank $m$ vector bundle on
$Z$,
and
$L$ a line bundle on $Z$. Suppose that $S^2(E^{*}) \otimes L$ is an ample 
vector bundle and that there is a constant even rank $r \ge 2$
symmetric bundle map $E \to E^{*} \otimes L$. Then $\delta \le m-r$.
\endproclaim

\demo{Idea of proof}
Consider the case that $Z\subset\ppp S^2 V^*$ is a 
variety of
quadrics. 
Consider the incidence correspondence
$$
\ci = \{ (x,P)\mid x\in P\}\subset Seg(\ppp V\times Z)
 \subset\ppp (V\ot S^2V^*)\tag 8.5.4
$$
$\ci$ is the intersection of the Segre with a cubic hypersurface,
and thus we can use the Lefschetz theorem to conclude that  $h^i(\ci )= h^i
(\pp m\times Z )$ for $i<m+\delta -1$. By the Kunneth formula,
$h^i (\pp m\times Z )= \oplus_{l+s=i}h^l(\pp m)h^s(Z )$.
$\ci$ can be considered as a fibration $Q\ra\ci\ra Z$.
Since $\pi_1 (Z) =0$ there is no monodromy and were the fiber smooth,
the Leray spectral sequence would say that
we could still compute the cohomology
as if it were a product,
i.e. $h^i (\ci) =\oplus_{l+s=i}h^l(Q)h^s(Z )$.
If we had the above results, we could conclude  $h^{r-2}(Q)=2
\neq 1=h^{r-2}(\pp m)$ and  unless
$\delta \le m-r$ there would be a contradiction.

We can easily modify things to make the fiber smooth, we 
let $\ci '$ be the fiber bundle over $Z$ with fiber the smooth quadric
obtained by quotienting $Q$ by its singular locus,
and note that $\ci '$ is contained in a   projective bundle  whose
fibers are $\ppp (V/Q_{z,sing})$. Now  Leray applies to $\ci '$ but
unfortunately Lefshetz in its standard form no longer applies 
when comparing $\ci '$ and this vector bundle. Fortunately  
  a  version of the Lefschetz theorem due to Lazarsfeld
  does apply in this situation and one obtains
the result. \qed\enddemo

\proclaim{8.5.5 Problem} Classify the boundary cases $\delta
=m-r$.\endproclaim

[IL,1.2] is an analogue for symmetric matrices of the following result
of Lazarsfeld:

\proclaim{8.5.6 Theorem   [Laz]}
 Let $X$ be a projective variety of dimension $m$. Let
$E$ and $F$ be vector bundles on $X$ of ranks $e$ and $f$ respectively.
Suppose that $E^* \otimes F$ is ample and that there is a constant rank
$r$ vector bundle map $E \to F$. Then $m \le e + f - 2r$.
\endproclaim

8.3.i.1 follows from  (8.5.6).
What follows is a   simplified 
version of Westwick's   proofs of 8.3.i.1
and   8.3.i.2, due to B. Ilic:
\demo{Proof of i.1 and the remainder of of i.2} 
let $K$ denote the kernel bundle   and
$N$ the cokernel bundle of $\phi$.
 From the resulting
 exact sequences, $c(K)c(E) = 1$ and $c(E)c(N) = (1+h)^n$.
Thus, $c(K)(1+h)^n = c(N)$. If $n-r+1 \le i \le l$ then $c_i(N) = 0$ and
looking at the coefficient of $h^i$ we get
$\sum_{j = 0}^{m-r} \binom{n}{i-j} c_j(K) = 0$ where we use the convention
that 
$\binom{n}{j} = 0$ if $j<0$ or $j > n$.  The coefficient matrix of this
collection of linear equations is 
$M = ( \binom{n}{i-j})_{0 \le j \le m-r, \  n-r+1 \le i \le l}$.
If $l = m + n -2r +1$ then this is a square invertible matrix with
determinant $\prod_{j=0}^{m-r} j!$. Thus $k_0 = 0$ which is a contradiction
since $k_0 = 1$. This proves i.1.
Westwick refers to a privately published manuscript of Muir and
 Metzler as a reference for evaluating
this determinant however one can also refer to e.g.  \cite{ACGH, pg.
93-95}.
i.2 follows directly from considering the linear equation with $i = n-r+1$.
\qed
\enddemo

\subheading{ 8.6 Systems of split type}

 In this section we restrict
attention to systems of quadrics, the modifications for the
skew case being more or less clear.

\proclaim{8.6.1 Theorem  (Sato) [Sato]} Let $E^r\ra\pp\delta$
be a uniform vector bundle. If $\delta >r$, then $E$ splits as a direct sum
of line bundles. If $\delta =r$
then either $E$ splits or $E$ is isomorphic to $T\pp\delta$ or
$T\pp{\delta *}$.\endproclaim

The proof of Saito's theorem relys on Tango's result about maps
of projective spaces to Grassmanians (which follows from an elementary
Chern class calculation) and sheafy gymnastics.

We will call systems of quadrics 
whose associated vector bundle   $E$ splits, 
{\it of split type}.

\proclaim{8.6.2 Conjecture}
If $E^r\ra\pp\delta$ is a vector bundle
arising  from a system of quadrics
and if
$\delta >\frac r2 +1$, then $E$ splits as a direct sum of line bundles.
\endproclaim

 Note
that this is the best bound possible due to the system 
  $I_2(G(2,5))$.  
I see no reason for the corresponding
conjecture to be true in the skew case,
as the corresponding Lagrangian Grassmanian is much larger than
the Spinor variety. 

 To study this conjecture let's examine maps of projective
spaces into spinor varieties as follows
(Chern classes are not enough): 

Let $A\subset\ppp S^2T^*$ be a system of constant rank.
Fix   $Q_0\in A$. Consider $\ppp (A/ A_0 )=\pp{\delta -1}=
\Bbb F (Q_0)$, the space of lines in $A$ through $Q_0$.
There is a natural vector bundle $L$ over $\Bbb F(Q_0)$,
where
$$
L_{[Q]}= \tspan_{\{[s,t]\in \pp 1\}}\tker (sQ+tQ_0)  \tmod \tker Q_0.\tag
8.6.3 
$$

 Since all pencils $\{ Q,Q_0\}$ are of  the same  type, $L$ is
indeed a vector bundle and not just a sheaf.

Note that $L$ is a subbundle of the trivial bundle  
$\Cal O^r = \Cal O\ot(V/\tker Q_0)$.
It is of rank $\frac r2$  (see [HP]). Label the resulting exact
sequence:
$$
0\ra L^{\frac r2}\ra \Cal O^{r}_{\Bbb F (Q_0)}\ra \Lambda^{\frac r2}\ra 0
\tag 8.6.4
$$
 $L_{[Q]}$ is
a $Q_0$-isotropic subspace of $(V/\tker Q_0)$
since $L_{[Q]}+\tker Q_0= W_{\{Q,Q_0\}}$.
Thus
 the 
image of the map to the 
Grassmanian $G(\frac r2, V/\tker Q_0)$ defined by
$L$ actually lies in the 
$Q_0$-spinor variety, $\Bbb S_{Q_0}$, a connected component of
$\{ E\in G(\frac r2,(V/\tker Q_0)) | E\subset  Q_0\}$.
Label the resulting map
$$
\align
\phi : \pp{\delta -1}&\ra   \Bbb S_{Q_0}\tag 8.6.5 \\
[Q]&\mapsto L_{[Q]}\endalign
$$

  $L= \Lambda^*$, as $L=\phi^*(\Cal S)$, $\Lambda=\phi^*(\Cal Q)$
where $\Cal S,\Cal Q$ denote the pullbacks of the universal sub-
and quotient bundles
on the Grassmanian
to  $\Bbb S$,  and   $\Cal S = \Cal Q^*$ when pulled back to $\Bbb S$.
This proves the assertion $\trank L=\frac r2$.
(The isomorphism is skew symmetric because $T_E\Bbb S\simeq
\Lambda^2E$.)

To prove the conjecture, one would need   to show that 
if $\frac r2<\delta -1$, then $L$ is trivial, i.e. $\phi$ maps
to a point (as then $W_{\{Q,Q_0\}}$ is independent of $Q$).

\proclaim{8.6.6 Proposition} If $k>1$, there are no nonconstant maps  
  $\pp{k}\ra \Bbb S$ or
$\pp k\ra \Bbb G_{Lag}$ such that 
$\phi^*(\Cal S)$ is  a split vector bundle,
where $\Cal S$ denotes the universal subbundle on the Grassmanian
restricted to $\Bbb S$ (resp. $\Bbb G_{Lag}$).
\endproclaim

Thus it would be sufficient to show that $L$ must be uniform,
because  if $L$ were uniform in this range, it would have to  be split.

\demo{Proof of proposition}
Given a map $\phi$, one gets a sequence as (8.6.4).
If $L$ splits we may write
$$
\align
&c(L)=\Pi_{i=1}^{\frac r2}(1+\alpha_it ) \ \ \alpha_i\leq 0\tag 8.6.7 \\
&c(\Lambda)=\Pi_{i=1}^{\frac r2}(1-\alpha_it ) \endalign
$$
We have
$$
1=c(L)c(\Lambda ) = 1-\{\Sigma_i(\alpha_i)^2\}t^2 + \hdots
$$
which implies $\alpha_i=0$ for all $i$ and thus $L$ is trivial.\qed
\enddemo

The conjecture would also follow from the following
stronger conjececture:

\proclaim{8.6.8 Conjecture} If $p<q$ then there are no nonconstant
maps, $\pp p\ra \Bbb S_{q}$\endproclaim

For the equality case $p=q$ there are linear examples. One might
hope to prove the conjecture by a degeneration argument.

\subheading{8.7 Split type systems and dual varieties}

\proclaim{8.7.1 Proposition} If $X^n\subset\pp\na$ is a linear
fibration with $f>\frac n2$ dimensional fibers, then
$\dstar\geq 2f-n$\endproclaim

It seems reasonable to conjecture that split type systems
correspond to fibrations with linear fibers. This cannot
be quite correct, as if $\dstar =1$, then $|II_{X^*,H}|$ is
necessarily of split type.  

Fix $H= [\uee\na ]$. Let $N'\subset N_xX$ be $(\uee\na )\upperp$
and let $L=\hat q^{\na}_{sing}\subset T$.
Consider the linear map
$$
\align
T/L &\ra L^*\ot N' \tag 8.7.2\\
v&\mapsto ( w\mapsto II(v,w) )\endalign
$$
which is well defined, as since $q^{\na}$ is generic, its
singular locus is in the base locus of the system.

\proclaim{8.7.3 Theorem [L10]} If $II_{X^*,H}$ is of split type
and (8.7.2) is injective, then $X$ is a linear fibration with
$\frac{n+\dstar}2$-dimensional fibers.\endproclaim

For the proof, see [L10]. Note that if $\dstar$ is large and
$X$ is smooth, then (8.7.2) is automatically injective.

Note that Ein's fibration theorem agrees with this conjecture, using
the standard bound. That is, in the geometric situation,
$r=n-\dstar$ (and $\dstar=\delta$), so Ein proves the variety
must be a fibration in the range $\ii_{X^*,H}$ is known to be split
type. If our conjecture is correct, then combined with the other
results, this would say that
if $\dstar>\frac{n+1}3$, then $X$ must be a fibration.

The fibration conjecture combined
with  [RVdV] would imply that $\dstar >
\frac{n+2}3$ implies $a>\frac 23 n+1$. 
 
\subheading{8.8   \lq\lq Dual\rq\rq\  systems of quadrics}

This section contains some observations that   should
have geometric consequences, although I have not yet been able
to determine them.

Let $T^n,N^a$ be vector spaces and let $II\in S^2T^*\ot N$ be
a system of quadrics on $T$ parametrized by $N$ of bounded rank $r$.
Given a generic $\alpha\in N^*$, let $W=W_{\alpha}:=\tsingloc
II^*(\alpha)$.
Bertini's theorem implies that $W\subset\tbaseloc \ii$ so
for each $\beta\in N^*$, the map
$$
II^*_{\beta}: T\ra T^* \tag 8.8.1
$$
descends to a well defined mapping
$$
{II^*_{\beta}}': W\ra W\upperp \tag 8.8.2
$$
i.e. each generic $\alpha\in N^*$ gives a tensor
$$
L_{\alpha}'\in (N^*/\{\alpha\})\ot W_{\alpha}\ot W_{\alpha}\upperp\tag
8.8.3
$$
which we consider as a system of linear forms  $L_{\alpha}$ on
$(N^*/\{\alpha\}) \ot W_{\alpha}\upperp$ parametrized
by $W_{\alpha}$.
  Call this system
the system {\it induced by $\alpha$}.

What is remarkable is that the doubling of the
induced  system occurs
naturally in the geometric setting. More precisely,
examining the inversion formula ([IL], 3.9), one has the following
proposition:

\proclaim{8.8.4 Proposition} Let $X^n\subset\pp\na$ be a variety
with degenerate dual variety $X^*$ of dimension $r+a-1$,
let $x\in \xsm$, let $H\in X^*_{sm}$ be a hyperplane tangent
at $x$ and let $n_H\in N^*_xX$ be a vector representing $H$.
Then  $II_{X^*,H}$ modulo the contributions
of $\pii_{X,x}$ is the doubling of the   system induced
by $n_H$. (After one makes identifications as in [IL], 3.9).
\endproclaim

This places further, to my knowledge so far unstudied, restrictions
on the second fundamental forms of varieties with degenerate dual
varieties. For example the system induced by $n_H$
must be of bounded rank $\frac r2$.

 An interesting example is when $X$ is self-dual, as in that case,
for all $\alpha\in N$, the doubling of the induced system is the
original system. One might hope to characterize all such systems
of quadrics.
 
\heading \S 9. Varieties with degenerate 
secant and tangential varieties\endheading

\subheading{9.1 Examples} The only examples of smooth varieties with
degenerate
secant varietes that I am aware of essentially come from 
varieties of  rank one
matricies. Its clear that the secant variety will be the set of rank
less than or equal to two matrices, which often  has
a secant defect. 
Consider the basic examples

i. Segres, $\pp m\times \pp n$, $n\geq m$, $\dsigma = 2$ if $n\geq 2$.

ii. Veroneses (symmetric matrices) $v_2(\pp n)$,  $\dsigma = 1$ if $n\geq
2$.

iii. Grassmanians of two-planes  (skew symmetric matrices) $G(2,n)$,
$\dsigma =
4$ if
$n\geq 6$.

iv. Severi varieties, $\Bbb A$-Hermitian symmetric
$3\times 3$  matrices. $\BA\pp 2$,  $\dsigma = \tdim_{\Bbb C}\Bbb A$

There is no known example of a variety with $\dsigma >8$.

In what follows, I only mention secant or tangential varieties in a
hypothesis, with the corresponding statement for the other understood.

\proclaim{9.2  Theorems on Secant varieties}Let $X^n\subset
 \pp{n+a}$ be a smooth  variety not contained in a hyperplane.

i. (Zak's Theorem on Linear Normality, [FL],[Z])  If  $\sigma
(X)\neq  \pp{n+a}$, then $a\geq \frac{n}{2}+2$.

ii. (Zak's Theorem on Severi Varieties,   [LV],[Z])   If $a=\frac{n}{2}+2$,
and    $\sigma
(X)\neq  \pp{n+a}$,
then $X$ is one of
the four Severi varieties $\Bbb A\pp 2\subset\ppp \Cal H$.

iii. (Zak-Fantecci superadditivity Theorem [Z],[Fan]
 [L5])   
$$
\text{dim}\sigma_k(X)\leq n+ (k-1)(n+1-\dsigma).
$$

iv.  ([Roberts], [L5])
Let $Y\subset \ppp V$ be a variety and let
$X=v_d(Y)\subset\ppp S^dV$ be the Veronese re-embedding.
If $d>2$ or
$d=2$ and $Y$ is not a linear subspace of $\ppp V$,
then $\sx$ is nondegenerate.

v. [Zak, V.2.3] Let $M(n,\delta )$ denote the maximum ambient dimension
of a smooth variey of dimension $n$, secant defect $\delta$. Then
$M(n,\delta )\leq \frac{n+\delta+2 +\epsilon (\delta -\epsilon
-2)}{2\delta}$
where $\epsilon$ is the remainder of the division of $n$ by $\delta$.

vi. [GH]    $\delta_{\gamma}(\tau (X))\geq \dtau (X)+1$. (Note
that $\delta_{\gamma}(\tau (X))\geq 1$ for any tangential variety.)

vii. [Zak] $\dstar (\tau (X))\geq \delta_{\sigma}(X)+1$
 
\endproclaim

The     Severi varieties   have many other special properties. For example,
they classify the quadro-quadro Cremona transforms
(see [ESB]).  These four varieties turn up in several areas of
geometry to construct examples of varieties
 exhibiting extremal pathology (isoparametric submanifolds [C, III.1
p1447],
tight embeddings [Kuiper]).
 
The $d>2$ case of iv. is due to Roberts, the $d=2$ to myself.

\subheading{9.3 The refined third fundamental form}

For $v\in T$ and $A\subset S^2T^*$, recall the notions
$$\align
\tann (v) &:= \{ q\in A\mid [v]\in q_{sing}\}\tag 9.3.1\\
\tsingloc (A)&:=\{ [v]\in\ppp T \mid  [v]\in q_{sing}
\forall q\in A\}\endalign
$$

Recall the map $\gamma\up 2 : X\ra G(n+a_1,V)$ defined in (*).
If $II(v,T)\subsetneq II(T,T)$, $\gamma\up 2$ can be refined to a map
$$
\align  \tilde\gamma\up 2 : TX &\dashrightarrow G(n+a_0,V)\tag 9.3.2 \\
(x,v)&\mapsto \hat T+ \widehat{II}(v,  T)\endalign
$$

Let $v\in T=T_xX$ and let $II_v$   denote the mapping:
$$
\align
II_v: T&\ra N\tag 9.3.3 \\
w&\mapsto II(v,w)\endalign
$$

Taking the derivative of $\tilde\gamma\up 2$, we obtain a refinement of
the third fundamental form:
In place of $V/\hat T\up 2$ we get the larger space $V/(\hat T
+ II_v(T))\simeq \tann (v)^*$. In place of $T$, we must
restrict to $SA(v):=
\tsingloc\tann (v)$, in order to have a mapping symmetric
in the three factors.

Fix a
$II$-generic vector $v\in T$, define
the {\it   third fundamental form refined with respect to $v$},
$$
  III^v  \in S^3 (\tsingloc (\tann (v)))^*\ot N/II_v(T)\tag 9.3.4
$$

   $III^v$ is well defined as if $w\in
\tsingloc (\tann (v))$, then $II_w(T)\subseteq II_v(T)$.
 
\smallpagebreak

To study
$III^v$, it is natural to consider the frame bundle over general points
$(x,v)\in TX$, and for each such, we get a refinement of the
flag $\hat x\subset\hat T\subset V$ to
$$
\hat x\subset\{\hat x,v\}\subset \hat T\subset
\{\hat T+ II(v,v)\} \subset \{\hat T+ II(v,T)\}\subset V
$$
In fact we may refine to a flag
$$
\hat x\subset\{\hat x,v\}\subset \{\hat x, SA(v)\}\subset \hat T\subset
\{\hat T+ II(v,v)\} \subset \{\hat T+ II(v,T)\}\subset V
$$
in fact to a partially ordered
flag, because we also have the flag
$$
\{\hat x\}\subset \{\hat x,\tker II_v\}\subset \{\hat x, SA(v)\}\subset
\hat T\subset
\{\hat T+ II(v,v)\} \subset \{\hat T+ II(v,T)\}\subset V
$$

 One can
recover $III$ from knowing $III^v$ for all $v\in T$.

\proclaim{9.4 Theorems on the local differential geometry
of varieties with degenerate secant varieties}
Let $X^n\subset\pp {n+a}$
be a
   a open subset of a variety   
 with degenerate tangential   manifold of
dimension
$n+a_0$.

i. (Terracini ??) [GH] [L5] Let $X^n\subset\ppp V$ be a smooth
variety. let $x\in X$ be a general point and $v\in T_xX$ a
generic tangent vector. Then
$$
\text{dim}\tx = n+ \text{dim}II_v(T)  
$$
In particular, if $a\leq n$, then $\dtau = \tdim \tker II_v$.

ii.  Terracini ??  [GH]  Let $X^n\subset\ppp V$
be a smooth variety, let $x\in X$ be a general point.
Let  $v\in T=T_xX$ be $II$-generic, then
$$
\text{dim}\sx = \left\{\matrix
n+ \text{dim}II_v(T) &\text{ if }III^v(v,v,v)= 0\\
 n+ \text{dim}II_v(T)  + 1 &\text{ if }III^v(v,v,v)\neq
0\endmatrix\right\}.
$$

iii. $III\neq 0$ implies $\sx$ is nondegenerate.

iv.   Rank restrictions for varieties
with degenerate tangential varieties, [L5, 7.1]
 Let  $X^n\subset\pp{n+a}$ be a smooth variety with degenerate secant
variety
 that is a hypersurface.
Let $x\in X$ be a general point.
 Let $\rtau$ be the maximum rank of a quadric
in $|II|_x$ annhilating a $II$-generic vector in $T_xX$.
Then
$$
\rtau\geq n-a+2.
$$
\endproclaim

\subheading{9.4.1 Exercise} Show directly that if $a>\binom{n+1}2$,
then $\dsigma (X)=0$.

\subheading{9.4.2 Exercise} Show that the condition $\tdim II_v(T)=a_0$
is equivalent to showing the image of the rational map
$$
\align ii: \ppp T &\dashrightarrow \ppp N\\
[v] &[II(v,v)]\endalign
$$
has an $ a_0$-dimensional image, or, in other words, that the
quadrics in $\ii$ satisfy $a-a_0$ polynomial equations.

\smallpagebreak

The proof of i. is given below. 
To calculate $\tdim\sx$ one takes a general
point of
$\sx$ written as the sum of two points on $X$ and calculates the
tangent space by expanding a Taylor series around one of the
points. See [L5 \S 10] for details.

vi. is a
  stronger  rank restriction than for arbitrary varieties  in two ways:
first, the rank is higher, and second, one looks not at generic quadrics,
but   at
generic quadrics in
the subvariety of quadrics annhilating a $II$-generic
vector.
Note that in vi, $\dtau +\rtau\leq n$ with equality only in the  
the borderline case of Zak's theorem on linear normality.

To prove vi, one studies the subvariety
$X^*_{\Delta}\subset X^*$ of hyperplanes $H$ such that at some
point of   tangency, $H$ annhilates a $II$-generic vector.
One first concludes that over  a 
general point $x\in X$,
$\tdim (X^*_{\Delta}\cap \ppp N^*_xX)\geq\rtau$,
and then shows that in fact equality occurs.

Each $H\in X^*_{\Delta}$ is  
tangent to a $\dtau$-dimensional subvariety $Y_H$
of $X$ whose
tangent space at $x$ is $(q_H)_{sing}$.
In the case of the Severi varieties
 these subvarieties 
$Y_H$ are quadric hypersurfaces in a linear space. They were
exploited in Zak's classification  Severi varieties.

\subheading{9.5 Dimension of $\tx$}

 The dimension of $\tau (X)$ is the dimension
of its  tangent space at a smooth  point.
We work on $\cf^1$,
which we consider as a bundle over the smooth points of
$\tx$ by the mapping
$$
\align \cf& \ra \tx\tag 9.5.1 \\
f& \mapsto [\ee 1].\endalign
$$
 We may take $[e_1]$ as
a typical element of $\tau (X)$. Let $2\leq \rho, \sigma\leq n$.  By
(*),
$$
de_1  \equiv  \ooo  {0} {1}e_0  +\ooo  {\rho} {1}e_{\rho} +
\ooo  {\mu} {1}e_ {\mu} \ \ \text{mod}\ \{ e_1 \} \tag 9.5.2
$$
and dim$T_{[e_1]}\tx = $ $\{$ the number of independent  1-forms in (9.5.2)
$\}$.  Over $X$,
$\{\ooo 0 1 ,\ooo\rho 1\}$ are independent so we only need to know the 
number
of independent one-forms among the $\ooo{\mu}{1}$, but (up to twisting)
$$
\ooo  {\mu} {1}\ot (e_{\mu}\text{mod}\hat T)\equiv  II_{X,x}({\underline
e_1},\cdot ).
 \tag 9.5.3
$$
where $\underline e_1=\aa 0^*\ot (\aa 1\tmod \aa 0 ) \in T^*$.
(9.5.3) proves i.

\subheading{9.6 Towards an infinitesimal version of the connectedness
theorem}

Note that we can compute
$II_{\tx,[\ee 1]}$ from third order information at
$x=[\ee 0]$. Adapt further such that $II_{X, \underline
e_1}(T)=\{\ee{n+h}\}$.
Let $1\leq\lambda\leq a_0$, $a_0+1\leq\epsilon\leq n$ 
and adapt such that $\tker  (II_{\underline e_1} )=\{\ee\epsilon\}$.
Then

$$
II_{\tau(X),[\ee 1]}=
(\rr{n+h}1\lambda k\oo \lambda  
+2\qq{n+h}jk\ooo j1 )\oo k\ot\ee{n+h}\tag 9.6.1
$$

  $\underline e_{\epsilon}\in \tsingloc II_{\tau(X),[\ee 1]}$,
which proves 9.2.vi.

If $X$ is smooth, the connectedness theorem implies 
$III^v=0$. In frames, if we let, $\tsingloc (\tann (\ee 1))=
\{ \ee 1,\ee s,\ee \epsilon\}$, then $III^v=0$ implies further that
$$
\rr{n+h}1 1 k=\rr{n+h}1 s k=0.\tag 9.6.2$$  

  Griffiths and Harris refer
to (9.2.vi) as the \lq\lq  deepest result\rq\rq\ 
in their paper. They speculate that it should be related to the
connectedness
theorem. (9.4.iv) is valid for an open subset $X^0$ of
any variety $X$, but
$\overline{\tau (X^0)}\neq\tau (X)$ in general, essentially
because one uses the tangent star instead of the tangent
cone to define $\tau (X)$. See [L5] for a discussion.
Thus there is no hope for a purely local version of the
Fulton-Hansen theorem (3.6).

On the other hand, on might hope to prove
directly that if $X$ is smooth, then 
$III^v=0$
holds soley based on smoothness considerations
or {\it a priori} restrictions on $F_3$ present
when  $II$ takes
a certain form (see \S 13 ). For example, it would
be sufficient to prove that $|II|\up 1=0$.

 Such a result would
give a proof of (3.6) that is local in nature in the case $X$
is smooth. Examples of second fundamental forms
 where  the additional vanishing
(9.6.2) is automatic are given in [L9].

 \heading\S 10.  Systems of quadrics with tangential defects\endheading

The main  result 
of the second step in my study of  the infinitesimal
geometry of varieties with degenerate secant varieties is the following:

\proclaim{10.1 Theorem [L5]} 
Let $\ii\subset\ppp S^2T^*$ be a system of quadrics arising from the second
fundamental form at a general point of
a smooth variety
  $X^n\subset\pp\na$, $a\leq n$  having a 
  degenerate secant variety 
that is a hypersurface.

Let $v\in T$ be a generic tangent vector. Consider the mapping
$II_v :T\ra N$, defined by $II_v(w)=II(v,w)$. Let 
  $P\in\ii$ denote the unique quadric such that $[v]\in P_{sing}$.
Then

 i. $\tker II_v\subset P_{sing}$.

ii. If 
$\{ v, \tker II_v\}= P_{sing}$, then there exists a canonical (possibly
degenerate) quadratic form $Q_v$ on $\tker II_v$.

iii. Under the hypotheses of ii., 
$T/P_{sing}$ is a $Clifford(\tker II_v, Q_v)$ module.

In other words, in this case,   each smooth point of $\sx$, produces a
canonical
$Clifford(\bcc\dsigma )$ action on a vector space of dimension less
than $n-\dsigma$.
\endproclaim

Note that assuming $\sx$ is a hypersurface is no loss of generality
because one can always project  until $\sx$ is a hypersurface.
The conclusions of ii, iii are valid is slightly more general
contexts.

\noindent{\bf 10.2 Definition} Let
$T$ be an $n$-dimensional vector space and let
 $A\subset S^2T^*$ be an $a$-dimensional
system of quadrics, with $a\leq n$. We say $A$ has a {\it tangential
defect} of
$\delta_{\tau}$ if for all $v\in T$, there exists a 
$\dtau$-dimensional
subspace $U_v\subset  A$ such that for all $q\in U_v$,
$[v]\in q_{sing}$, i.e. $q(v,w)=0$ for all $w\in T$.

\smallpagebreak
 
Equivalently, consider the rational map:
$$
\align ii :\ppp T &\dasharrow \ppp A  \\
x&\mapsto [Q^1(x)\hd Q^a(x)]\endalign
$$
then $\dtau = \tcodim ii (\ppp T)$. To see the equivalence
of the two definitions, note that $II_v$ is  
the derivative of $ii$ at $[v]$ after taking appropriate quotients.

It is often useful  to consider a pararmetrized system of quadrics to
mimic the geometric situation. We will write $II\in S^2T^*\ot N$
where $N$ is an $a$-dimensional vector space, and $\ii\subset\ppp
S^2T^*$ denote the corresponding system.

For simplicity, we often assume $II$ corresponds to   the second
fundamental
form at a point of a variety $X$ where
$\sx$ is a hypersurface. In this case, for generic $v\in T$,
$q=\tann (v)$ is unique. Call this a {\it critical tangential defect}.
Note that for a critical tangential defect,
$\dtau = n-a+1$.

In the case of a critical tangential defect, 
for any $II$-generic $v\in T$, there is a natural 
hyperplane in $\ii$, namely $II(v,v)\upperp$.

The following proposition is analgous   to the fact that 
$\hat x\subset \tker d\gamma_x(v)$
for all $v\in T$:

\proclaim{10.3 Proposition}$\tker II_v\subset \tsingloc\tann
(v)$\endproclaim

  The proof follows from the notion of the projective
second fundamental form of a mapping, see [L5]. 

\subheading{10.3.1 Exercise} Prove (10.3) using moving frames.

\smallpagebreak

Systems 
of quadrics with a tangential defect satisfy a Bertini
type theorem. Recall that the classical
Bertini theorem says that if 
$A\subseteq S^2T^*$ is any system of quadrics (or polynomials
of any degree for that matter) and 
$q\in A$ is a general element, then
$q_{sing}\subseteq\tbaseloc (A)$. We can't hope for such a strong statment
for $q=\tann (v)$, because $v\notin\tbaseloc \ii$. But in fact the next
best
thing is true.

\proclaim{10.4 Bertini type Lemma, [L5, 6.16]} Let $II\in S^2T^*\ot N$
be a system of quadrics with a critical tangential defect. Let
$v\in T$ be $II$-generic,   let $q\in\ii$ be the annhilator of $v$. 
Then $q_{sing}\subseteq\tbaseloc (II(v,\hat q_{sing})\upperp)$.
\endproclaim

This Bertini type lemma implies
that $\ii/II(v,\hat q_{sing})\upperp$ can be considered
 as a system of quadrics
on $\hat q_{sing}$. In the case of a critical defect, we obtain
    a well defined quadric which we denote
$\tilde Q_v$.

Note that there are large linear subspaces
in $\tbaseloc \ii$ because   any
vector in the base locus
of $\tilde Q_v$ must be (as an element of $\ppp T$) in the base locus
of
$\ii$.  Moreover,   there is a positive dimensional variety
of different $\tilde Q_v$'s.

 Fix $v\in T$, $II$-generic. We have the
subspaces $\tker II_v$ and $\tsingloc (\tann (v))$ giving a flag
$$
\tker II_v
\subset\{ v,\tker II_v\}\subseteq \tsingloc (\tann (v))\subset T\tag 10.5 
$$
For future reference, we record here the  way these spaces will
be referred to in indices:
$$
\{\ee\epsilon\}
\subset\{ \ee 1, \ee\epsilon\}\subseteq
\{ \ee 1, \ee\epsilon,\ee s\}=\{\ee\xi\}\subseteq 
\{ \ee 1, \ee\epsilon,\ee s,\ee j\}=\{\ee\alpha\}\tag 10.6
$$
I.e. $\{ \ee\epsilon\}= \tker II_v$,
$\{\ee s\}=\tsingloc (\tann (v))/\{ v,\tker II_v\}$,
$\{\ee j\} =T/\tsingloc (\tann (v))$. 
 
\proclaim{10.7 Theorem, the canonical Clifford algebra structure, [L5,
6.23]}
Let $A\subset S^2T^*$ be a system of quadrics with a critical tangential
defect, and let $v\in T$ be $II$-generic. 
Assume $\ii/II(v, \hat q_{sing})\upperp$ restricted to $\tker II_v$
corresponds to a single quadratic form which we denote $Q_v$. Then
$T/(\tann (v))_{sing}$
is canonically a $Clifford(\tker II_v, Q_v)$ module.
\endproclaim

\subheading{10.8 Remark/question}In the cases
 $Q_v$  exists,
it is likely its rank is linearly bounded
from below, and thus in this situation,
the tangential defect $\dtau$ could
grow at best logrithimically with respect to $\tdim X$. This is because
$\tdim\tker II_v=\dtau$ and\linebreak
$n>\tdim (T/(\tsingloc (\tann (v)))$. Were $Q_v$ nondegnerate,
then $\tdim (T/( \tsingloc (\tann (v)))$ would have to be on the order
of $2^{\dtau -1}$.

\smallpagebreak

\demo{Proofs}
Before proving (10.4), let's warm up by proving the classical
Bertini lemma for an arbitrary system of quadrics. To keep notation
consistent
for what comes, I'll assume the system of quadrics occurs as
the second fundamental form at some point of a open subset of a variety .
Writing the system as $\{\ooo\mu \alpha\oo\alpha =\qq\mu\alpha
\beta\oo\alpha\oo\beta\}$, assume $q^{\na}$ is a generic quadric
in the system 
(further indices,
$n+1\leq\mu ,\nu\leq  \na,
 1\leq  \lambda\leq a -1$)
and that $\ppp\{\ee \xi\}= q^{\na}_{sing}$, $1\leq
\xi ,\eta\leq
\tdim q^{\na}_{sing}$. Differentiating, we have
$$
\rr\na \xi\eta\beta \oo\beta = \qq{n+\lambda} \xi\eta\ooo\na{n+\lambda}
 \tag 10.9
$$
Since we have assumed $q^{\na}$ is a generic quadric, the forms
$\ooo\na{n+\lambda}$ are all linearly independent and independent of the
semi-basic forms $\oo\beta$. Thus both sides of the equation above
have to vanish and hence the coefficients $\qq{n+\lambda} 
\xi\eta$ are zero,
proving the classical Bertini theorem for a system of quadrics.\enddemo

\demo{Proof of 10.3}
For simplicity assume we have a critical tangential defect, so
that
for a $II$-generic vector $v\in T$, $\tann (v)$ consists
of a single quadric. We adapt frames such that $\ee 1=v$ and $q^{\na}
=\tann (v)$.
Since $q^{\na}$ is no longer a generic quadric in the system,
the forms $\ooo\na{n+\lambda}$ are no longer necessarily independent or
independent of the semi-basic forms.

Let $[\ee j]$ denote a basis of $T/q^{\na}_{sing}$,
$2\leq j,k,l\leq a-1$. Since we have
a critical defect, $II_v: T\ra N$ has rank $a-1$ and
  $\trank (q^{\na})= a-2$.  We may normalize further such that
$II(v,v)=\ee\none$, i.e. that $\ee 1\in\tbaseloc\{ q^{n+k},q^{\na}\}$,
and in fact that
$II(\underline e_1,\underline e_s)=\underline e_{n+s},
II(\underline e_1,\underline e_j)=\underline e_{n+j}$. We
need to show that $q^{\na}_{sing}\subset\tbaseloc\{ q^{n+k},q^{\na}\}$,

Computing we have
$$
\align
&\rr\na 11\beta\oo\beta =  \ooo\na\none\tag 10.10\\  
&\rr\na 1j\beta\oo\beta = - \ooo\na{n+j} + \qq\na jk\ooo
k1 \tag 10.11\\
&\rr\na 1s\beta\oo\beta = - \ooo\na{n+s}  \tag 10.12\endalign
$$
Now for the essential point:
Since $\ee 1$ is generic and the $\ee k$ were taken generically,
the forms $\ooo k1$ are linearly independent and independent of
the semi-basic forms $\oo\beta$. Moreover, the matrix
$\qq\na jk$ is invertible.
 From (10.10), (10.13) we see that $\ooo\na\none, \ooo\na{n+s}$ are
semi-basic. 
(10.11)
implies we can solve for
$\ooo\na {n+j}$ in terms of the $\ooo k1$ modulo the 
semi-basic forms, so   the forms $\ooo\na{n+j}$ are independent and
independent
of the semi-basic forms. Now consider (10.9)
modulo the basic forms. We have
$$
0\equiv \qq{n+j} \xi\eta\ooo\na{n+j}\tmod\{\oo\alpha\}\tag 10.13
$$
(where we are still using $\ee \xi$ as a basis for $q^{\na}_{sing}$
although index ranges have been shifted).
  (10.3) follows as the forms  $\ooo\na{n+j}$ are all independent
and independent of the semi-basic forms.\qed\enddemo

\demo{Proof of the Clifford algebra structure}
Assume we are in the case
$\{ v,\tker II_v\} = \tsingloc (\tann (v))$.
(This is the case if $a=\frac n2+2$.)
Consider 
\def\nj{n+j}
$$
\align
&\rr\nj 11\beta\oo\beta= -\ooo\nj\none + 2\ooo j1 ,\tag 10.14
\\
&
\rr\nj 1\epsilon\beta\oo\beta = \ooo j\epsilon + \qq\nj\epsilon k \ooo k1.
\endalign
$$
which imply
$$
\ooo\nj\none\equiv 2\ooo j1\tmod \{\oo\alpha\}.\tag 10.15
$$
$$
\ooo j\epsilon \equiv -\qq\nj \epsilon k\ooo k1\tmod\{\oo\alpha \}.
$$
Computing
$$
\rr\nj\epsilon\delta\beta\oo\beta =
-\qq\none\epsilon\delta\ooo\nj\none +
\qq\nj\epsilon k\ooo k\delta + \qq\nj\delta k\ooo k\epsilon  \tag 10.16
$$
and moding out by the semi-basic forms and using (10.15), we obtain,
$$
\qq\nj\epsilon k (-\qq\nk\delta i\ooo i1 ) +
\qq\nj\delta k (-\qq\nk\epsilon i\ooo i1 )\equiv\qq\none\epsilon\delta
(2\ooo
j1)\tmod\{\oo\alpha\} .\tag 10.17
$$
I.e. that
$$
\qq\nj\epsilon k
\qq\nk\delta i + 
\qq\nj\delta k
\qq\nk\epsilon i
= -2\qq\none\epsilon\delta \delta^i_j
\ \forall \epsilon,\delta, j,k,i\tag  10.18
$$
Consider the map 
$$
\align 
\tker II_v &\ra\text{End}(T/P_{sing})\tag 10.19\\
w^{\epsilon}\ee\epsilon &\mapsto
w^{\epsilon}\qq\nk j\epsilon (\ee j)^*\ot\ee k\endalign
$$
By (10.18), the fundamental lemma of Clifford algebras applies.\qed\enddemo

\heading{\S 11  Recognizing uniruled varieties}\endheading

I now prove (0.10.2) stated in the introduction, in fact a more
general result, but first we need a technical definition:

Fix $x\in X$ and let $L$ denote an osculating linear space.
Assume the Gauss map of $X$ is nondegenerate
Let $\pi :\cf\ra X$ denote the frame bundle.
If $L$ osculates at least  to order $k+1$,
then when one restricts the differential invariants
$F_j\in \pi^*(S^jT^*_xX\ot N_xX)$   to $S^jT_xL^*\ot T_xL\upperp$, one
obtains
 maps
$$
R_j :
S^jT_xL^*\ot N \ra T_xL\ot T_xL\upperp \tmod \tim R_{j-1}. \tag 11.1.2
$$

If these maps
all have maximal rank for $2\leq j\leq m$, we will say that 
{\it the maximal rank condition
holds through level $m$}. The maximal rank condition
is a pointwise
genericity condition on the subspace of the space of tensors for 
$F_2\hd F_{m}$ having the property that there is a $k$-plane in
their common base locus.

\proclaim{11.1.3 Theorem [L7] (Expectation Theorem)} Let $(n,a,k,m)$ be  
natural numbers satisfying $m\geq 3$ and
$$
a[\binom{k+m-1}{m-1} -k -1]\geq k(n-k).
$$
Let $X^n\subset\Bbb A^{\na}$ or $X^n\subset\pp\na$ be a open subset
of a smooth (respectively analytic)
submanifold of an affine or projective space
having the properties that at each
(resp. at a general) $x\in X$ there exists a
$k$-dimensional linear space $L_x$,
disjoint from the fiber of the Gauss map,
 osculating to order $m$
and such that
the maximal rank condition holds through level $m-1$.
 Then $L_x$ is locally contained in $X$.

  If $F$ denotes
the fiber of the Gauss map at $x$ and $\tdim L_x\cap F=\lambda$, then the
same
conclusion holds as above with $k,n$ respectively replaced by
$k-\lambda, n-\lambda$.
In fact, if the Gauss map is degenerate, one
can replace $L$ by the span of $L$ and $F$.
\endproclaim

\proclaim{11.1.4 Theorem [L7]} Let $X^n\subset\Bbb C^{\na}$ or
$X^n\subset\Bbb C\pp\na$
 be a open subset  
of a complex analytic submanifold of an affine or projective space, and let
$x\in X$ be a general point.
If $n\geq 4$ and a linear space
$L^{n-2}$  osculates to order four at $x$, then locally 
$L^{n-2}\subset X$.
\endproclaim

\proclaim{11.1.5 Proposition [L7]}   There exist 
analytic open subsets of   varieties $X^n\subset\Bbb A^{n+a}$
 and $X^n\subset\pp{n+a}$ 
having a line osculating to order $\frac{n-1}a+2$ 
at every point that are not ruled.
In fact over $\Bbb C$, every open subset of a variety  has this property.  

\endproclaim

\demo{Proof of   11.1.3}
For notational simplicity, I only prove the case $\gamma$ is nondegenerate.
  Let
$(v_1\hd v_n)=(v_{\xi},v_{\rho} )$ denote a basis
of $T_xX$ adapted such that   $T_xL= \{ v_{\xi}\}  $.
The index ranges are $n+1\leq\mu,\nu\leq n+a$,
$1\leq \xi,\eta\leq k $,
$k  +1\leq\rho,\sigma\leq n $.

By hypothesis $F_j(v_{\xi_1}\hd v_{\xi_j})=0$ for all $1\leq j\leq m$.
The coefficients of $F_{i+1}$ in this range are given by the formula
([L6],2.20) which simplifies under our hypotheses to
$$
\rr \mu{\xi_1\hd}{\xi_i}\rho\oo\rho =
\frak S_{\xi_1\hd \xi_i}\rr
\mu{\xi_1\hd}{\xi_{i-1}}\sigma\ooo\sigma{\xi_i}.
 \tag 11.1.6
$$
(11.1.6) expresses the forms $\ooo\sigma\xi$ in terms of 
the forms $\oo\rho$. Geometrically, say we were in the case of a unique
linear space $L$ at each point. Then we would  have a map:
$$
\align l: X&\ra G\tag 11.1.7 \\
x&\mapsto L_x
\endalign
$$ 
where $G$ denotes the appropriate Grassmanian of
$k$ planes. In this case  the forms $\ooo\sigma\xi$
correspond to a spanning set of the cotangent space of the image.
In the general case one still has such a map, only  from the  
bundle over $X$ whose fibers parametrize the space of $L$'s through
a point. We wish to
show the mapping is constant along tangent directions to $L$, that is
the pullback of the forms $\ooo\sigma\xi$ are zero when restricted to
$T_xL$.

In terms of tensors, for each element
of $S^iT_xL^*\ot N$ we obtain a (possibly zero) element of $T_xL\upperp$
that is
identified with an element of $T_xL^*\ot T_xL\upperp$.
 Suppressing the element of $T_xL\upperp$
from the notation,
we obtain   (11.1.2).

We see that if these maps are all of largest possible rank, then at
level $i$ we will have filled a
$$
a[
\binom {k+1}2
 +\binom {{k}+2}3 + \hdots +\binom {{k}+(i-2)}{i-1} ]
 = a[\binom {{k}+i-1}{i-1} -i-1]\tag 11.1.8
$$
dimensional subspace of $T_xL^*\ot T_xL\upperp$, which is of dimension
$(n-k)
 k $.
 
Assuming the maximal rank  condition holds through level $m-1$
and $a[\binom{{k}+m-1}{m-1}-k-1]\geq  {k} (n-k)$, we see
that   $\ooo\sigma\eta\equiv 0\tmod \{\oo\rho\}$ $\forall \sigma,\eta$.
Intuitively, by the remark above, this finishes the proof. In details,
by induction, for $i>m$ we have 
$$
\rr \mu{\xi_1\hd}{\xi_i}{\xi_{i+1}}\oo{\xi_{i+1}} +
\rr \mu{\xi_1\hd}{\xi_i}\rho\oo\rho =
\frak S_{\xi_1\hd
\xi_i}\rr{\mu}{\xi_1\hd}{\xi_{i-1}}\sigma\ooo\sigma{\xi_i},
\tag 11.1.9
$$
which implies
$\rr {\mu}{\xi_1\hd}{\xi_i}{\xi_{i+1}}=0$
for all $i$, and thus the linear space
osculates to infinte order.
\qed\enddemo

To prove theorems 1 and 5, one observes that if the maximal rank condition
fails in these cases, the relevant invariants are forced to be zero anyway.

\heading{\S 12 Quadrics containing $X$}\endheading

 In this section I derive the generalized Monge system
  described in (0.6).

$X$ is locally the intersection of quadrics if $N^*_xX$ is spanned
by the differentials of quadratic equations.

In order that $N_x^*$ be spanned by differentials of quadratic
polynomials, it is necessary that
$$
\{ dP_x | P\in \tker \bii k\vtwox \} = N_x^*\tag 12.1.k
$$
for all $k$.
(We supress reference to the base
point $x$ in what follows.)
 For $k\leq 2$, (12.1.k) automatically holds; for
$k=3$ (12.1.3) will hold if and only if

$$
\ff {3\mu} = 3a^{\mu}_{\nu\gamma}\oo\gamma\ff {2\nu}
\tag 12.2
$$
for some constants $a^{\mu}_{\nu\gamma}\in \Bbb C$.
Notice that if
$\rr\mu\alpha\beta\gamma =  \frak S_{\alpha\beta\gamma}
a^{\mu}_{\nu\gamma} \qq\nu\alpha\beta  $
in some frame, it holds in any choice of frame (with different
constants $a^{\mu}_{\nu\gamma}$), so the expression (12.2) has intrinsic
meaning. If (12.2) holds, then
 $$
\tker\bii 3\vtwox =
\{\xx\mu\xx 0 -
\qq\mu\alpha\beta\xx\alpha\xx\beta-a^{\mu}_{\nu\beta} \xx\nu
\xx\beta,
\xx\mu\xx\nu \}.\tag 12.3
$$
Continuing in the same fashion, we uncover the following conditions:
\def\al#1#2#3{a^{#1}_{{#2}{#3}}}
\def\fff#1#2{\Bbb F\Bbb F^{{#1} {#2}}}
%temp def
$$
\align &
 \fff 3\mu = 3\al\mu\nu\gamma\oo\gamma\fff 2\nu \tag 12.4\\
&\fff 4\mu = 4\al\mu\nu\alpha\oo\alpha\fff 3\nu + 3 \bb\mu\nu\tau
\fff 2\nu\fff 2\tau\\
&\fff 5\mu = 5\al\mu\nu\gamma\oo\gamma\fff 4\nu + 10\bb\mu\nu\tau\fff
3\nu\fff 2\tau
\endalign
$$
%endtemp def zone
\def\fff#1#2{\Bbb F\Bbb F^{#1}_{#2}}
  where $ \al\mu\nu\alpha , \bb\mu\nu\tau=\bb\mu\tau\nu\in \Bbb C$.
Moreover, if there are no linear syzygies among the quadrics in $\ii$,
as explained in \S 14, then $\ff 6=0$ and thus
 $N_x^*$ is  spanned by the differentials of 
quadrics and these  uadrics are smooth along $X$ so
they generate $I(X)$.    In this
case,   we will call (12.4) the
{\it generalized Monge system for quadrics}.

%back to ususal def
\def\al#1{\alpha_{#1}}

 In summary:

\proclaim{12.5 Theorem [L6]}
Let $X\subset\bppp V$ be a variety and $x\in X$ a general point.
Assume $III_{Xx}=0$   and that there are
no linear syzygies in $ \ii_x$.
Then
$$
\align &
\tdim \{ \text{quadrics osculating to order three at } x \}\leq
a
+\binom{a+1}2-1  \tag 12.6\\
&
\tdim \{ \text {quadrics osculating to order four at } x \}\leq
a-1.
\endalign
$$

If  the generalized Monge system (12.4) holds, then
$$
 I_2=\tker\bii 4{\vtwox x} 
$$
Equality occurs in the first (respectively second)  line of  (12.6)
 if and only if    the first
(resp. second) line of (12.4) holds at $x$.    If the
generalized Monge system  does not hold, then $I_X$ is not generated
 by quadrics.
 \endproclaim

 If one assumes appropriate genericity conditions, there exist analogous
Monge equations for $I_d$ of order $2d+1$ in small codimension, see [L6].

\heading{\S 13 Recognizing homogeneous varieties}\endheading

Before studying homogeneous varieties, let's consider a general question:

If one fixes the dimension and codimension of $X$, there is an integer
$k_0$ such that if $F_2\hd F_{k_0}$ satisfy some mild genericity
hypotheses, then all $F_l$'s are determined by the 
$F_2\hd F_{k_0}$'s and their derivatives. Since $X$ is analytic,
this means that $X$ is entirely determined by $k_0$ derivatives
on any open set. For a hypersurface, with $n\geq 2$, Jensen and
Musso [JM] proved that $k_0=3$ (the case $n=2$ is due to Cartan
and Fubini).
For a plane curve, $k_0=6$. In general, the function 
$k_0(n,a)$ appears to be unknown.

\subheading{ 13.1 Problem}   What, if any,
values of $a$ have $k_0(n,a)=2$?

\smallpagebreak

Note that the Codazzi equation in Euclidean geometry implies
that if $III^{\Bbb E}=0$, then $k_0^{\Bbb E}=2$.

\subheading{13.2 Locally uniform varieties}

Now let's specialize our question to the case of {\it locally
uniform} varieties. 

\subheading{ 13.2.1 Definition} A variety $X\subset\ppp V$ is
{\it locally uniform to order $k$} in the neighborhood of a point $x\in X$
if there exists a local framing (section of $\cf^1$) in which the
coefficients of $F_l$ are constant for $2\leq l\leq k$.

For example, all hypersurfaces are locally uniform to
order two in   neighborhoods of general points.

\subheading{13.2.3 Remark} R. Bryant has proven that varieties
that are locally uniform to order $\infty$ are in fact
locally homogeneous (personal communication).

\subheading{ 13.2.2 Questions} How does (13.1) simplify for
locally uniform varieties? Are locally uniform varieties
always rational (assuming $a\neq 1$ and perhaps some additional
conditions)? How can one characterize the varieties of codimension
two that are locally uniform?

\smallpagebreak

Recall that  the coefficients of $F_3$ are given by:
$$
\rr\mu\a\b\gamma\oo\gamma =
- d\qq\mu\a\b -\qq\mu\a\b\ooo 0 0 -\qq\nu\a\b\ooo\nu\mu
+\qq\mu\a\c\ooo\c\b + \qq\mu\c\b\ooo\c\a
$$

\subheading{13.2.3 Exercise}
We can always take a local framing such that $\ooo 00 =0$.

\smallpagebreak

On any section of $\cf$ we have
$\ooo\mu\nu= a^{\mu}_{\nu\a}\oo\a$, $\ooo\a\b= b^{\a}_{\b\c}\oo\c$
for some functions $a^{\mu}_{\nu\a},b^{\a}_{\b\c}$. 
Thus, assuming $d\qq\mu\a\b=0$, i.e. local uniformity to order two,
 we obtain restrictions
on $F_3$. Namely, at each point, let
$$
A\in N\ot N^*\ot T^*\ \ \ B\in T^*\ot T\ot T^* \tag 13.2.3
$$
then
$$
F_3\in A\cdot II + S(B\cdot II) \tag 13.2.4
$$
where $A\cdot II$, $B\cdot II$ are the natural contractions and $S$ is
symmetrization in {\it two} factors.
 
 $F_3$ is symmetric in all
three factors which places significant restrictions
on the admissible $A$'s and $B$'s.   Moreover,
$F_3$ can be modified by $T^*\circ II + <II, T^*\ot N, II>$ so
one might hope that using all these conditions to show
$F_3$ must be zero in certain situations, or at least severly limited.

\subheading{13.2.6 Problem} Give a clean cohomological description of the
restrictions described above.

\smallpagebreak

In the example below, these conditions are {\it not} enough to conclude 
$F_3=0$, higher order considerations must be taken into account.

Further differentiation
  places systems of partial differential equations on the  functions  
$a^{\mu}_{\nu\a},b^{\a}_{\b\c}$, and it is these overdetermined
systems of pde that finally allow one to conclude $F_3=0$ in
the calculations below.

\subheading{13.3 Homogeneous varieties}

Homogeneous varieties usually have rather special differential
invariants. For example, $II$ of a homogeneous variety $X=G/P$ must
be invariant under the action of the semi-simple part of $P$ on
$S^2T^*\ot N$.

We have seen in the introduction that even when $II$ is invariant
under a large group, it may not determine $X$ (e.g. the adjoint
varieties).

\subheading{13.4 Idea of the proofs}

The idea of the proofs is as follows: given
any variety $X\subset\ppp V$,  consider
the first order adapted
frame bundle, 
  $\pi :\cf^1_X\ra X$.
Each $f\in\cf^1$ determines a splitting of the
flag $\hat X\subset\hat T\subset V$ which we denote $\hat x + T + N$.

Write the pullback of the Maurer-Cartan form of $GL(V)$ to
$\cf^1$
as
$$
\Omega = \pmatrix \oo 0&\ooo 0\beta&\ooo 0\nu\\
\oo\alpha&\ooo\alpha\beta&\ooo\alpha\nu\\
0&\ooo\mu\beta&\ooo\mu\nu\endpmatrix
$$
with index ranges $1\leq\alpha,\beta\leq\tdim X$,
$\tdim X+1\leq\mu,\nu\leq \tdim\ppp V$.

If $X=G/P$, one can reduce
$\cf^1$ until it is isomorphic to $G$ (with fiber
isomorphic to $P$). In that case one obtains the
Maurer-Cartan form symbolically as:
$$
\Omega_G = \pmatrix \oo 0&\ooo 0\beta&0\\
\oo\alpha&\ooo\alpha\beta=\rho_T(\frak h) &\ooo\alpha\nu
=A_2(\ooo 0\beta ) \\
0&\ooo\mu\beta=
A_1(\oo \alpha ) &\ooo\mu\nu=\rho_N(\frak h) \endpmatrix
$$
where 
$H$ is the semi-simple part of $P$,
$T=T_xX$, $N=N_xX$ are $H$-modules with representations $\rho_T$, $\rho_N$,
and
$A_1,A_2$ are $H$-equivariant maps.
  The zero in the upper
right hand block indicates that  any infinitesimal change in
the splitting statisfies the
\lq\lq transversality\rq\rq \
condition that $dN\subseteq \{ T + N\}$.
The dependence of the
$\ooo\alpha\nu$ block   on the forms $\ooo 0\beta$
  indicates that if one changes the choice of $T$,
there is a corresponding change in choice of $N$ mandated.

If $X$ is a variety with the same second fundamental form as $G/P$,
by restricting bases  
we can reduce $\cf^1_X$ to a bundle $\cf^2_X$
where the pullback of the the Maurer-Cartan
form looks like:
$$
\Omega  = \pmatrix \oo 0&\ooo 0\beta&\ooo 0\nu\\
\oo\alpha&\ooo\alpha\beta=\rho_T(\frak h)+ w_1 &\ooo\alpha\nu
  \\
0&\ooo\mu\beta  &\ooo\mu\nu=\rho_N(\frak h)+ w_2 \endpmatrix
$$
where $w_1,w_2$ are linear combinations of the other
forms appearing in the Maurer-Cartan form.
The proofs proceed  by showing that there are  
reductions of $\cf^2_X$ to $G$.

In practice, the rigidity proofs proceed by showing the
  invariants $F_k\in \pi^*(S^kT^*X\ot NX)$  are zero for $k>2$.

\subheading{13.5 Proof of the Grassmanian case}

 $z=0$ because $\ii\up 1=0$.

\def\dee#1#2{e_{{#1}{#2}}}
 
Let $V$ have basis $\{\ee 0,\dee 1j,\dee 2j,\dee jk\}$, where 
$3\leq j,k,l\leq n+2$, $\{\alpha\}=\{ 1j, 2j\}$.
Normalize such that
$
II= (\oo{1j}\oo{2k}-\oo{1k}\oo{2j})\ot\dee jk, \ j<k
$.
Note that   the forms $\ooo {1i}{1j},\ooo {2i}{2j},
\ooo {1i}{1i},\ooo {2i}{2i},\ooo {2i}{1i},
\ooo {1i}{2i} $
are all independent and independent of the semi-basic forms because
they represent infinitesimal motions that preserve our normalization
of $II$.
We have
$$
\align &
\rr{ij}{(1k)}{(1l)}\beta = 0\ \forall i,j,k,l\text{ distinct  and }
\forall\beta\tag 13.5.1\\
 &
\rr{ij}{(2k)}{(2l)}\beta = 0\ \forall i,j,k,l\text{ distinct  and }
\forall\beta\tag 13.5.2\endalign
$$
(these equations imply that the refined third fundamental
form is zero).
 From now on, assume all indices are distinct.
Using (13.5.1), (13.5.2), we have
$$
\align &
\rr{(ij)}{(1i)}{(1k)}{(1i)}\oo{1i}+
\rr{(ij)}{(1i)}{(1k)}{(1j)}\oo{1j}+
\rr{(ij)}{(1i)}{(1k)}{(2i)}\oo{2i}+
\rr{(ij)}{(1i)}{(1k)}{(2j)}\oo{2j}+
\rr{(ij)}{(1i)}{(1k)}{(2l)}\oo{2l}\tag 13.5.3\\
&\ \ \ =\ooo{(2j)}{(1k)} \endalign
$$
 The right hand side of (13.5.3) is independent of $i$,
so comparing with the same expression using $m$ instead of $i$,
(here we use that $n\geq 4$)
we obtain:
$$
\align
&\rr{(ij)}{(1i)}{(1k)}{(1i)}=0\tag 13.5.4\\
&\rr{(ij)}{(1i)}{(1k)}{(1j)}=\rr{(mj)}{(1m)}{(1k)}{(1j)}
  \tag 13.5.5\\
&\rr{(ij)}{(1i)}{(1k)}{(2i)}=\rr{(mj)}{(1m)}{(1k)}{(2i)}
  \tag 13.5.6\\
&\rr{(ij)}{(1i)}{(1k)}{(2j)}=\rr{(mj)}{(1m)}{(1k)}{(2j)}
  \tag 13.5.7\\
&\rr{(ij)}{(1i)}{(1k)}{(2l)}=\rr{(mj)}{(1m)}{(1k)}{(2l)} .
  \tag 13.5.8\endalign
$$
Now
$$
\align
&
\Delta \rr{(ij)}{(1i)}{(1k)}{(2l)}= g^{(2j)}_{(kl)}
\\
&
\Delta \rr{(ij)}{(1i)}{(1k)}{(2j)}= g^{(2j)}_{(jk)} 
+g^{0}_{(1k)} .
\endalign
$$
Using these equations and the corresponding
equations with the role of $1$ and $2$ reversed, we reduce to frames where
$ 
 \rr{(ij)}{(1i)}{(1k)}{(2l)}= 0, 
  \rr{(ij)}{(1i)}{(1k)}{(2j)}=0,
 \rr{(ij)}{(2i)}{(2k)}{(1l)}= 0,
  \rr{(ij)}{(2i)}{(2k)}{(1j)}=0$. 
In these frames, $\ooo{1k}{2j},\ooo{2k}{1j}=0$ hence
$$
0=\rr {(ij)}{(1i)}{(1i)}{\beta}\oo{\beta}=-2\ooo{2j}{1i}\tag 13.5.9
$$
and similarly with the role of $1$ and $2$ reversed.
Thus the only nonzero terms left in $F_3$ are
$ \rr{(ij)}{(1i)}{(1j)}{(1k)},\rr{(ij)}{(2i)}{(2j)}{(2k)}$.
Consider
$$
\rr {(ij)}{(1i)}{(1j)}{(1i)\beta}\oo{\beta}=2
 \rr{(ij)}{(1i)}{(1j)}{(1k)}\ooo{1k}{1i}\tag 13.5.10
$$
Both sides of (13.5.10) must be zero because the forms
$\ooo{1k}{1i}$ are all independent and independent of
the semi-basic forms. The analogous equation holds
with $2$'s. Hence we see $F_3=0$.

To have a nonzero coefficient of $F_4$, $\rr{(ij)}\alpha
\beta{\gamma\delta}$, in the lower indicies there must
be two $1$'s and two $2$'s, and at least two of the $k$-indices
must be $i$ or $j$. Consider
$$
\align
&\rr {(ij)}{(1i)}{(1k)}{(2l)(2j)}\oo{2j}=\ooo{2j}{kl}\tag 13.5.11\\
&\rr {(ij)}{(1i)}{(1i)}{(2l)(2j)}\oo{2j}=2\ooo{2j}{il}.\tag
13.5.12\endalign
$$
Since the right hand side of (13.5.11) is independent of
$i$, we conclude (after switching the roles of $i$ and $j$)
that $\rr {(ij)}{(1i)}{(1k)}{(2l)(2j)}$ is independent of $i,j$
(with neither $k,l$ equal to $i$ or $j$, but $k=l$ is possible). Using
$$
\Delta \rr {(ij)}{(1i)}{(1k)}{(2l)(2j)}= g^0_{kl}
$$
we normalize all these terms to zero. This implies $\ooo{2j}{il}=0
$
and thus $\rr {(ij)}{(1i)}{(1i)}{(2l)(2j)}=0$ for all $i,j,l$ distinct
as well, and similarly with the role of $1$ and $2$ reversed.
Thus the remaining nonzero terms in $F_4$ are
$\rr {(ij)}{(1i)}{(1i)}{(2j)(2j)},\rr {(ij)}{(1i)}{(2i)}{(1j)(2j)}$.
Consider
$$
\align 
&\rr {(ij)}{(1i)}{(2j)}{(1k)(2l)(1i)}\oo{1i}+
\rr {(ij)}{(1i)}{(2j)}{(1k)(2l)(1j)}\oo{1j}\tag 13.5.13\\
&+
\rr {(ij)}{(1i)}{(2j)}{(1k)(2l)(2i)}\oo{2i}+
\rr {(ij)}{(1i)}{(2j)}{(1k)(2l)(2j)}\oo{2j}\\ 
& = -\ooo 0{kl}.\endalign
$$
Since the right hand side of (13.5.13) is independent of $i,j$, we
conclude $\ooo 0{kl}=0$ and hence the left hand side is
zero as well. Now it is easy to see the rest of the terms in
$F_5$ are zero and all higher forms are zero.\qed

\subheading{13.6 Problem} Determine rigidity for the case $G(2,5)$.

\heading \S 14 Complete intersections\endheading

While the other pathologies we have so far studied could be
calculated locally, failing to be a complete intersection is
a global issue. For example at any $x\in\xsm$, $X$ is locally
a complete intersection.
To attempt to recognize whether or not a variety is a complete
intersection from computable local information would be futile.
In this section I will discuss computable
local conditions that insure $X$ is a complete intersection.

\smallpagebreak

 Let $V=\bcc{n+a+1}$ and let $X^n\subset\bppp V=\bcc{}\bpp{n+a}$ be a
variety
of dimension $n$.
Let $\xsm$ denote the smooth points of $X$.
 Let $I_X\subset S^{\bullet}V^*$ denote
the ideal of $X$
and let $I_{X,d}=I_d = S^dV^*\cap I_X$ denote
the $d$-th graded piece of $I_X$. Fixing a smooth point  $x\in X$,
there is a distinguished subspace of  $I_d$, namely the
hypersurfaces of degree $d$ that are singular at $x$, i.e. $P\in I_d$ such
that
$(dP)_x=0$, where $dP$ denotes the
exterior derivative of the polynomial  $P$.

\smallpagebreak

The following definition is due to Lvovsky [Lv]:

\subheading{14.1  Definition [L6]}. Let $X\subset\ppp V$ be a variety. Let
$P\in I_d$ and let $Z=Z_P\subset\ppp V$ be the corresponding hypersurface.
We will say
$Z$ {\it trivially contains} $X$ if
$P= l^1P_1+\hdots l^mP_m$ with $P_1\hd P_m\in I_{d-1}$ and $l^1\hd l^m
\in V^*$, and   otherwise
that  $Z$ {\it essentially contains} $X$.

\smallpagebreak

[L6, 1.1] stated in the introduction generalizes to the following
statement:

 \proclaim{14.2 Proposition [L6,1.6],
A local characterization of complete intersections} Let $X\subset\bppp V$
be a
variety. The following are equivalent:

1. $X$ is a complete intersection.

2. Every hypersurface
essentially containing $X$ is smooth at all $x\in\xsm$.

3.  Let $x\in \xsm$. Every
 hypersurface
essentially containing $X$ is smooth at $x$.
 \endproclaim

(14.2) localises the study of complete intersections to
a point, and further, filters the conormal bundle at that point
to enable us to study one degree at a time. Unfortunately, to determine
if a hypersurface essentially contains $X$, one might need to take
an arbitrarily high number of derivatives. To have computable
conditions, we will work  with osculating
hypersurfaces rather than the
hypersurfaces
 containing $X$. The advantage will be that we will only
need to study a fixed number of derivatives for each fixed degree of
hypersurface;  the disadvantage is that we will only obtain
sufficient conditions to be a complete intersection.

By [L4, 3.16, 3.17] stated in the introduction, we see that at
best one could prove there are no singular hypersurfaces of
degree $d$  osculating to order $2d+2$ at $x$; and that the first 
restrictions one could hope for are at order $d+1$.

\smallpagebreak

I now specialize to the case $d=2$.

Looking at (4.6.2), we see that
  $\tker\fff 3{\vtwox,x}$ is as small as possible if there are no
{\it linear syzygies} among the quadrics in $II_{X,x}$.

(If $A\subset S^2T^*$, a linear syzygy among the quadrics
in $A$ is a relation of the form $l_j\circ Q^j=0$, where
$l_j\in T^*$ and $Q^j\in A$. More invariantly, consider
the symmetrization map $\Cal S: T^*\ot S^2T^*\ra S^3T^*$ and its
restriction to $A$, $\Cal S|_A : T^*\ot A\ra S^3T^*$.
Let $A^{[1]}= \tker (\Cal S|_A)$. Then $A^{[1]}$ is the space
of linear syzygies of $A$.)

 We have
the following linear algebra  lemma:

\proclaim{14.3 Lemma [L4, 6.19]}Let $A^p\subset S^2T^*$ be an
$p$-dimensional
system of quadrics on an $n$ dimensional vector space.  Say there
is a linear syzygy
$$
l^1Q_1 +\hdots + l^pQ_p = 0
$$
where both
$l^i\in T^*$ and $Q_i\in A$ are  independent sets of vectors.
Then $\forall Q\in A$,
$$
\text{rank }Q\leq 2(p-1).
$$
\endproclaim

For the proof, see [L6].
If one now compares [L4, 6.19] with the rank restriction theorem,
one sees that  if $a<\frac{n-(b+1)+3}3$ then there are no linear
syzygies in $\ii$. Combined with  the generalized Monge
system, we obtain:

\proclaim{14.4 Theorem [L6, 6.26] }
Let  $X^n\subset\bpp{n+a}$ be a variety and $x\in X$ a general point.
  Let $b=\tdim X_{sing}$. (Set $b=-1$ if $X$ is smooth.)
If
$a<\frac{n-(b+1)+3}{3}
$
then
$$
\align &
\tdim \{ \text {quadrics osculating to order three at } x \}\leq
a
+\binom{a+1}2-1 \tag 12.5\\
&
\tdim \{ \text {quadrics osculating to order four at } x \}\leq
a-1.
\endalign
$$
 Equality occurs in the first (respectively second)  line of  (14.5)
 if and only if the 
generalized Monge system  holds to order three (respectively four) at
$x$.
  If  the generalized Monge system
   holds, then  $X$ is a complete intersection of the $a-1$
dimensional family of quadrics osculating to order four.
\endproclaim

 [L4, 6.28] stated in the introduction follows immediately
in the following stronger form:

\proclaim{14.6 Corollary [L6, 6.28]}Let  $X^n\subset\bpp{n+a}$ be a variety
and
 $x\in X$ a general point.
  Let $b=\tdim X_{sing}$. (Set $b=-1$ if $X$ is smooth.)
If
$a<\frac{n-(b+1)+3}{3}
$
then any quadric osculating to order four at $x$ is smooth at $x$ and
any quadric osculating to order five at $x$ contains $X$.
\endproclaim

  While the higher order Monge equations are more
complicated to write down, in principle they are no more
difficult to understand. Thus the problem of determining the 
subbundle of the conormal bundle consisisting of trivially containing
hypersurfaces is in principle resolved. The dengeneracy conditions on
the differential invariants of $X$,   while more complicated
to write down, is also, in principle resolved. The first condition
corresponds to the existence of (excess) linear syzygies. The first
real problem 
in attempting to generalize  (14.4)
is that for polynomials of higher
degrees, the presence of linear syzygies is not a serious pathology, in
particular, it has little relation to the \lq\lq cone locus\rq\rq 
(maximal multiplicity
 locus) of the polynomials. Moreover, while there are indications that
there may be higher order rank restriction theorems, none have
yet  been proven, and any I would be willing to conjecture
would relate to
restricting the cone locus.

\heading \S 15 Errata and clarifications\endheading

1. In [L3] p 315, it says that the Grassmanians $G(2,m)$ will have
degenerate $\gamma_k$. It should say the Segre's $Seg(\pp l\times\pp m)$
will have degenerate $\gamma_k$.

2. In [L6]   I give an example of a variety with $F_3=0$ but
$F_4\neq 0$ at a  point, but I do not show that that point
is a general point. An example  where the phenomena does occur
at general points is a cone over a curve in $\pp 2$. 
 
3. In [L5] 12.1,  I state that 
when $X$ is smooth and $\tau (X)$ is
degenerate, the fibers of the Gauss map
of $\tau (X)$ are of dimension at least $\dtau +2$, but what is proven
is that $\tau (X) $ has dual defect at least $\dtau +2$. The statement as
announced is false. Moreover, the observation of the defect of $\tau (X)^*$
was already made by Zak [Z].

\enddocument